\documentclass[10pt,a4paper]{article}
\usepackage{exscale, makeidx, amsmath, amssymb, 
comment,
epsfig, bbm,  mathrsfs,psfrag,float}
\usepackage{times}

\small
\normalsize
\usepackage[latin1]{inputenc}
\usepackage[T1]{fontenc}
\usepackage{amsfonts}
\usepackage{amsbsy}
\usepackage{amscd}
\usepackage{theorem}

\usepackage{epsf}
\usepackage{psfrag}
\usepackage{graphicx}

\usepackage{multicol}
\allowdisplaybreaks

\def\cqfd{ \hfill $\blacksquare$ }
\def\cq{ \hfill $\square$ }
\newcommand{\beq}{\begin{eqnarray*}}
\newcommand{\eeq}{\end{eqnarray*}}

\newcommand{\ben}{\begin{enumerate}}
\newcommand{\een}{\end{enumerate}}

\newcommand{\beqs}{\begin{eqnarray*}&\displaystyle}
\newcommand{\eeqs}{&\end{eqnarray*}}

\newcommand{\n}{{\rm n}}
\newcommand{\Br}{{\rm Br}}

\newtheorem{theorem}{Theorem}[section]
\newtheorem{lemma}[theorem]{Lemma}

\newtheorem{corollary}[theorem]{Corollary}
{\theorembodyfont{\rmfamily}\newtheorem{definition}[theorem]{Definition}}
{\theorembodyfont{\rmfamily}\newtheorem{remark}[theorem]{Remark}}
{\theorembodyfont{\rmfamily}\newtheorem{example}[theorem]{Example}}

\newenvironment{pf}{\begin{trivlist}\item[]\bf Proof. \rm}
                     {\hspace*{\fill} \cqfd\end{trivlist}}

\def\noi{\noindent}
\newcommand{\fsigma}{{\boldsymbol{\sigma}}}
\newcommand{\ftau}{{\boldsymbol{\tau}}}

\newcommand{\fTau}{{\boldsymbol{\cT}}}

\newcommand{\fFau}{{\boldsymbol{\cF}}}

\newcommand{\bbx}{{\bf x}}
\newcommand{\bx}{{\rm x}}

\newcommand{\ft}{\mathbf{t}}
\newcommand{\fT}{\mathbf{T}}
\newcommand{\fF}{\mathbf{F}}
\newcommand{\bD}{\mathbb{D}}

\newcommand{\bE}{{\bf E}}

\newcommand{\bN}{\mathbb{N}}
\newcommand{\bQ}{\mathbb{Q}}

\newcommand{\bP}{{\bf P}}
\newcommand{\bR}{\mathbb{R}}

\newcommand{\bT}{\mathbb{T}}
\newcommand{\bU}{\mathbb{U}}
\newcommand{\bZ}{\mathbb{Z}}

\newcommand{\cE}{\mathcal{E}}
\newcommand{\cB}{\mathcal{B}}
\newcommand{\cC}{\mathcal{C}}
\newcommand{\cF}{\mathcal{F}}
\newcommand{\cG}{\mathcal{G}}
\newcommand{\cH}{\mathcal{H}}

\newcommand{\cM}{\mathcal{M}}

\newcommand{\cD}{\mathcal{D}}
\newcommand{\cK}{\mathcal{K}}
\newcommand{\dl}{{\bf d}}
\newcommand{\dgh}{{\boldsymbol{\delta}}}
\newcommand{\Ell}{{\boldsymbol{\ell}}}
\newcommand{\btl}{\bT_{\Ell_1}}

\newcommand{\ccM}{\mathscr{M}}

\newcommand{\paste}{\mathrm{Paste}\, } 

\newcommand{\abv}{\mathrm{Abv} \, } 

\newcommand{\blw}{\mathrm{Blw}\, } 

\newcommand{\pnt}{\Upsilon}

\newcommand{\TR}{\, \widetilde{\textsc{T}}\textsc{ree}\, } 

\newcommand{\Trr}{\textsc{Tree}\, } 

\newcommand{\LevP}{\mathrm{P}} 

\newcommand{\cP}{\mathcal{P}}

\newcommand{\cT}{\mathcal{T}}

\newcommand{\lgeo}{[\![}
\newcommand{\rgeo}{]\!]}
\newcommand{\s}{\sigma}
\newcommand{\un}{{\bf 1}}
\newcommand{\vre}{\varepsilon}

\begin{document}
\title{
HEREDITARY TREE GROWTH AND L\'EVY FORESTS.
\thanks{This research was in part supported by
    EPSRC grant GR/T26368/01 and ANR A3 Projet BLAN.}  
}             
\author{Thomas \textsc{Duquesne}\thanks{Universit\'e P.~et M.~Curie, 
Laboratoire de Probabilit\'es et Mod\`eles Al\'eatoires, Bo\^ite courrier 188, 4 Place Jussieu, 
75252 Paris Cedex, France.} \and Matthias \textsc{Winkel}\thanks{University of
  Oxford, Department of Statistics, 1 South Parks Road, Oxford OX1 3TG, UK; 
        email winkel@stats.ox.ac.uk}\vspace{-0.1cm}}

\maketitle

\vspace{-0.7cm}

\begin{abstract} We introduce the notion of a hereditary property for rooted real trees and we also consider reduction of trees by a given hereditary property. Leaf-length erasure, also called trimming, is included as a special case of hereditary reduction. We only consider the metric structure of trees, and our framework is the space $\bT$ of pointed isometry classes of locally compact rooted real trees equipped with the Gromov-Hausdorff distance. Some of the main results of the paper are a general tightness criterion in $\bT$ and limit theorems for growing families of trees. We apply these results to Galton-Watson trees with exponentially distributed edge lengths. This class is preserved by hereditary reduction. Then we consider families of such Galton-Watson trees that are consistent under hereditary reduction and that we call growth processes. We prove that the associated families of offspring distributions are completely characterised by the branching mechanism of a continuous-state branching process. We also prove that such growth processes converge to L\'evy forests. As a by-product of this convergence, we obtain a characterisation of the laws of L\'evy forests in terms of leaf-length erasure and we obtain invariance principles for discrete Galton-Watson trees, including the super-critical cases. 

\smallskip 


\end{abstract}

\section{Introduction}
 This paper concerns general results on continuum trees and convergence of random trees. Here we view trees as certain metric spaces called real trees and we are using and developing a framework initiated by Aldous \cite{Al1, Al2} and Evans, Pitman and Winter \cite{EvPitWin} who first considered in our probabilistic context the space of compact real trees equipped with the Gromov-Hausdorff distance. The convergence results of our paper are applied to a large class of growth processes of Galton-Watson forests. This class of tree-growth processes contains the important example of forests consistent under leaf-length erasure 
(see Neveu \cite{Ne2}, and Le Jan \cite{LJ91} in the context of super-processes), and it is also closely related to two specific models considered by Geiger and Kauffmann \cite{GeiKau} and by the present authors in \cite{DuWi1}. 
We prove that in some sense, any way of growing Galton-Watson trees yields, in the limit, L\'evy trees, which are continuum random trees introduced by Le Gall and Le Jan \cite{LGLJ1} that have been further studied in \cite{DuLG} and also by Abraham, Delmas \cite{Abr-07} and Weill \cite{Weill}.

Let us briefly review in this introduction the main results of the paper. First we 
recall a few definitions on real trees and the space of trees we consider. A real tree is a path-connected metric space $(T,d)$ with the following property: any two points $\sigma,\sigma^\prime\in T$ are connected by a unique injective path denoted by $\lgeo \s, \s^\prime \rgeo$, 
which furthermore is isometric to the interval $[0, d(\s, \s^\prime)]$ of the real line. Informally, real trees are obtained by gluing together, without creating loops, intervals of $\bR$ equipped with the usual metric. 
However, note that real trees may have a complicated local structure, like Aldous's (Brownian) Continuum Random Tree, which is a compact real tree whose set of leaves is uncountable and dense.   
In each real tree $(T, d)$, we distinguish a point $\rho \in T $ that is viewed as the root. So we speak of $(T, d, \rho)$ as a rooted real tree. 

We shall focus on complete locally compact rooted real trees (CLCR real trees for short). We then say that two CLCR real trees are equivalent if there exists a root-preserving isometry from one tree onto the other. We simply denote by $\widetilde{T}$ the pointed isometry class of a given CLCR real tree $(T, d, \rho)$. We denote by $\bT$ the set of pointed isometry classes of CLCR real trees. We equip 
$\bT$ with the pointed Gromov-Hausdorff distance denoted by $\dgh$ (see Section \ref{realtrees} for a definition). Then $(\bT, \dgh)$ is a Polish space. This result is due to Gromov \cite{Gro} for compact 
metric spaces and to Evans Pitman and Winter \cite{EvPitWin} for compact real trees 
(see also \cite{DuWi1} for the standard adaptation to pointed CLCR real trees). 
The main results of this paper (tightness criterion, limit theorems and invariance principles) take place in the space $(\bT, \dgh)$.  
  
  Let us briefly explain the notion of a hereditary property in this context. Let $(T, d , \rho)$ be a CLCR real 
tree. Then we define for every $\sigma\in T$ the \em subtree above $\sigma$ \em as $\theta_\sigma T=\{\sigma^\prime\in T\colon\sigma\in\lgeo\rho,\sigma^\prime\rgeo\}$. Note that $(\theta_\sigma T,d,\sigma)$ is also a CLCR real tree. We denote by $\widetilde{\theta}_\sigma T$ its pointed isometry class in $\bT$.    
A \textit{hereditary property} is a Borel subset $A \subset \bT$ such that for every CLCR real tree 
$(T,d,\rho)$ and for every $\s \in T$, if $\widetilde{\theta}_\s T \in A$, then $\widetilde{T} \in A$. 
In order to rephrase this definition informally, let us view $T$ as a continuum of individuals whose progenitor is the root $\rho$: if   
an individual $\sigma\in T$ ``has the hereditary property $A$'', namely if $\widetilde{\theta}_\sigma T\in A$, then the progenitor $\rho$ also ``has the property $A$''; implicitly, a hereditary property may be lost on the ancestral lineage between the progenitor and an individual, and an individual can only inherit a hereditary property if all his ancestors had it.

  We then define the \textit{$A$-reduced subtree of $T$} as $R_A(T)$ where
$$\mbox{$R_A (T)$ is the closure in $T$ of the subset 
$ \{ \rho \} \cup \{ \sigma \in T\colon \widetilde{\theta}_\sigma T \in A   \} $.}$$ 
Then, $(R_A (T), d, \rho)$ is a CLCR real tree and its pointed isometry class only depends on the isometry class of $T$. Hence, there is an induced function from $\bT$ to $\bT$ that we simply denote by $R_A$. 
Hereditary properties can be \textit{composed} in the following sense: let $A, A^\prime \subset \bT$ be two hereditary properties, we then set $A^\prime \circ A = \{ \widetilde{T} \in \bT\colon\!\!$ $R_A (\widetilde{T}) \in A^\prime \}$ and Lemma \ref{comphereprop} asserts that $A^\prime \circ A$ is hereditary and moreover $R_{A^\prime \circ A} = R_{A^\prime} \circ R_A$. 

   The most important example of hereditary reduction is the \textit{leaf-length erasure} (also called trimming) that is defined as follows. For any CLCR real tree $(T, d, \rho)$, denote by $\Gamma (T)= \sup_{\s \in T} d(\rho, \s)$ its \textit{total height} (that is possibly infinite). Since it only depends on $\widetilde{T}$, it induces a function on $\bT$, that is also denoted by $\Gamma$ and that is $\dgh$-continuous. Then for any $h \in [0, \infty)$, we set $A_h = \{ \widetilde{T} \in \bT\colon \Gamma (\widetilde{T} ) \geq h \}$, which is clearly hereditary. We shall simply write $R_{A_h}$ as $R_h$ and refer to $R_h$ as the $h$-leaf-length erasure. Note that for any $h, h^\prime \in [0, \infty)$,  $A_{h^\prime} \circ A_h= A_{h+h^\prime}$ and thus, $R_{h^\prime } \circ R_h = R_{h^\prime + h}$.        

Leaf-length erasure was first considered by Kesten \cite{Kes86} for discrete trees. Then it was studied
by Neveu \cite{Ne2}, by Neveu and Pitman \cite{NP} to approximate the Brownian tree and also by Le Gall \cite{LG86}; later Le Jan \cite{LJ91} used it to construct superprocesses with a stable branching
mechanism. In the context of compact real trees, leaf-length erasure was more systematically used by Evans, Pitman and Winter \cite{EvPitWin}. They proved in particular that $R_h$ is $\dgh$-continuous.  

One important fact to note is that for any CLCR real tree $T$, $R_h( T)$ is a 
\textit{real tree with edge lengths}, that is $R_h(T)$ has a discrete branching structure: the set of 
branch points has no accumulation points and all branch points have finite degree. For every $a \in [0, \infty)$ and every $h\in (0, \infty)$, we set 
$$ Z_a^{(h)}(T)=\# \, \{ \,  \sigma\in R_h(T)\colon d(\rho,\sigma) = a \, \} .$$
Since $R_h(T)$ is a real tree with edge lengths, $a \mapsto  Z_a^{(h)}(T)$ is 
an $\bN$-valued function that is left-continuous with right limits. We call the process $a\mapsto Z_a^{(h)} (T)$ the \textit{$h$-erased profile}. Since $Z_a^{(h)}(T)$ only depends on the pointed isometry class of $T$, it induces a function on $\bT$ that is denoted in the same way.       

   As proved in Lemma \ref{keycompact}, the erased profiles allow to control the covering numbers of balls of CLCR real trees, which is the key argument in the following general tightness criterion in $(\bT , \dgh)$: 
let $(\widetilde{\cT}_j)_{j \in J}$ be a family of $\bT$-valued random trees. Their laws are 
$\dgh$-tight 
if and only if for every fixed $h, a \in (0, \infty)$, the laws of the $\bN$-valued random variables 
$(Z_a^{(h)}(\widetilde{\cT}_j))_{j\in J}$ are tight. This is Theorem \ref{tightness}, which is one of the main results of this paper. This tightness criterion is used to obtain, among other results, invariance principles for Galton-Watson trees 
(see Theorem \ref{invprincth} in Section \ref{Invprincsec}). 
More generally, it is well-adapted to any model of random trees that allows a certain control on the erased profiles.      
    
    In this paper we consider tree-valued processes that grow. 
 One simple way to understand the growth is to view the trees in a certain ambient metric space and 
 to say that the trees grow with respect to the inclusion partial order. However, we only consider the metric 
 structure and we want to consider neither the ambient metric spaces that allow such 
 constructions nor the details concerning the many ways the tree-growth process can be embedded in a 
 given ambient space. This is why we introduce the following intrinsic definition of growth: let $\widetilde{T}, \widetilde{T}^\prime \in \bT$; 
 we say that $\widetilde{T}$ can be ``embedded'' in $\widetilde{T}^\prime$, which is denoted by $\widetilde{T}\preceq \widetilde{T}^\prime$, if we can find representatives $(T, d, \rho)$ and $(T^\prime, d^\prime, \rho^\prime)$ of $\widetilde{T}$ and $\widetilde{T}^\prime$, and an isometrical embedding $f\colon  T\hookrightarrow T^\prime$ such that $f(\rho)= \rho^\prime$.       
    
   Note that $\preceq$ is a partial order on $\bT$. Moreover, for any hereditary property $A$, if $\widetilde{T} \! \preceq \! \widetilde{T}^\prime$ then $R_A(\widetilde{T})\! \preceq \! R_A(\widetilde{T}^\prime)$ and in particular for any $h \!>\! 0$, $R_h(\widetilde{T}) \! \preceq \! R_h 
(\widetilde{T}^\prime)$, which implies for any $a \in [0, \infty)$, $Z^{(h)}_a (  \widetilde{T}) \! \leq \!  Z^{(h)}_a (  \widetilde{T}^\prime)$. This yields the following convergence criterion stated in Theorem \ref{tightness2}: let $( \widetilde{\cT}_n)_{n \in \bN}$ be a $\bT$-valued sequence of random trees such that for all $n \in \bN$, $ \widetilde{\cT}_n \! \preceq \! \widetilde{\cT}_{n+1}$ almost surely; furthermore assume that  for every fixed $a, h\in (0, \infty)$, the family of laws of the $\bN$-valued random variables 
$Z_a^{(h)}(\widetilde{\cT}_n)$, $n\! \in\! \bN$, is tight; then, there exists a random tree $\widetilde{\cT}$ in $\bT$ such that $ \lim_{n \rightarrow \infty} \dgh (\widetilde{\cT}_n , \widetilde{\cT})= 0$ a.s. 
  
   This result is used to prove that, when convergent, any growing family of Galton-Watson forests tends to either a Galton-Watson forest or to a L\'evy forest. Before explaining this result, let us informally discuss  
the model of Galton-Watson real forests that are considered here (they are introduced precisely in Section \ref{GWrealtreesec}). 
Their laws are characterised by the following three parameters: 
the offspring distribution $\xi= (\xi(k))_{k\in \bN}$, the lifetime parameter $c \in (0, \infty)$, and the initial distribution $\mu= (\mu(k))_{k\in \bN}$. We view a Galton-Watson real forest with parameters $\xi$, $c$ and $\mu$ (a \textit{{\rm GW}($\xi,c;\mu$)-real forest} for short) as the forest of genealogical trees of a population that evolves as follows: at generation $0$, the population has $N$ independent progenitors, where $N$ has law $\mu$; the lifetimes of individuals are independent and exponentially distributed with mean $1/c\, $; when they die, individuals independently give birth to a random number of children distributed according to 
$\xi$. Here, we shall assume that $\xi$ is proper and non-trivial, namely that 
$\xi (1)=0$ and $\xi (0) <1$. Moreover we also assume that the associated continuous-time 
$\bN$-valued branching process is conservative (see (\ref{nonexpdiscr}) in Section \ref{GWrealtreesec} for more details).

Theorems \ref{realforreduc} and \ref{GWrealreduc} 
assert that the class of GW-laws is preserved by hereditary reduction. More precisely, 
let $\widetilde{\cF}$ be a GW($\xi, c; \mu$)-real forest, let $A\subset \bT$ be a hereditary property and denote by  
$\alpha$ the probability that the $A$-reduced tree of a single GW($\xi, c$)-real tree is just a point. Assume that $\alpha \in (0, 1)$. Then $R_A (\widetilde{\cF})$ is a 
GW($\xi^\prime,c^\prime;\mu^\prime$)-real forest, where ($\xi^\prime,c^\prime;\mu^\prime$) is given in terms of $\alpha$ and ($\xi, c; \mu$) as follows: if we denote by $\varphi_\nu$ the generating function of a probability measure $\nu$ on $\bN$, then 
\begin{equation}\label{xiprime}  \varphi_{\xi^{\prime}} (r)=  r + \frac{\varphi_\xi(  \alpha \!+\! (1\!- \!\alpha)r )\!- \!\alpha \!- \! (1\!- \!\alpha)r}{(1\!-\! \alpha)(1\!- \!\varphi_\xi^\prime(\alpha))} , \quad c^{\prime} \!=\!  (1\!-\! \varphi^\prime_\xi (\alpha))\, c  \! \quad  \textrm{and} \! \quad
 \varphi_{\mu^{\prime}} (r) \!= \! \varphi_\mu (\alpha + (1-\alpha)r). 
\end{equation}
This important property of hereditary reduction naturally leads us to consider families of GW-real forests that are consistent under hereditary reduction. Namely, we call a family $(\widetilde{\cF}_\lambda)_{\lambda\in[0,\infty)}$ of Galton-Watson real forests a \em growth process \em if for all $\lambda^\prime\ge\lambda$ there exist hereditary $A_{\lambda,\lambda^\prime}\subset\bT$, such that almost surely 
$$ \textrm{$\forall  \lambda^\prime, \lambda \in [0, \infty)$ such that $\lambda^\prime\ge\lambda$,}
\qquad \widetilde{\cF}_\lambda\!=\!R_{A_{\lambda,\lambda^\prime}}(\widetilde{\cF}_{\lambda^\prime})\; .$$
Let us say that $\widetilde{\cF}_\lambda$ is a GW($\xi_\lambda, c_\lambda; \mu_\lambda$)-real forest and assume that $\mu_\lambda$ tends to the Dirac mass at infinity, as $\lambda \rightarrow \infty$. This assumption implies that each offspring distribution $\xi_\lambda$ appears as a \textit{reduced law} as in (\ref{xiprime}) for all $\alpha$ sufficiently close to 1 (such offspring distributions are called infinitely extensible \cite{DuWi1}).   
Then Theorem \ref{structuregrowth} shows that they are quite specific. Namely, 
the laws ($\xi_\lambda, c_\lambda, \mu_\lambda$) are entirely governed by a triplet 
($\psi, \beta, \varrho$) defined as follows. 
\begin{itemize}

\item $\psi\colon [0, \infty) \rightarrow \bR$ is the \textit{branching mechanism} of a continuous-state branching process. Namely, $\psi$ is the Laplace exponent of a spectrally positive L\'evy process and it is therefore of the L\'evy-Khintchine form 
$$ \psi (\lambda) = {\bf a} \lambda + \frac{_1}{^2} {\bf b}
\lambda^2 + \int_{(0, \infty)} \left( e^{-\lambda x} -1+\lambda x \un_{\{ x<1\}} \right) \pi (dx) \; , $$
where ${\bf a} \in \bR$, ${\bf b}\in [0, \infty)$ and $\pi$ is a Borel measure on $(0, \infty)$ 
such that $\int_{(0, \infty)}(1\wedge x^2) \pi (dx) < \infty$. Moreover, $\psi$ has to satisfy two additional conditions: it takes 
positive values eventually and $\int_{0+} dr /|\psi (r)| =\infty$ (the latter assumption is equivalent to assuming that the continuous-state branching process governed by $\psi$ is conservative).

\item  $\beta\colon [0, \infty) \rightarrow [0, \infty)$ is non-decreasing and such that $\lim_{\lambda \rightarrow \infty} \beta (\lambda) = \infty$. 

\item $\varrho$ is a probability measure on $[0, \infty)$ that is distinct from the Dirac mass at zero. 
 
\end{itemize}

\noi  
Then, for any $\lambda \in [0, \infty)$, ($\xi_\lambda, c_\lambda, \mu_\lambda$) is derived from ($\psi, \beta, \varrho$) as follows: 
$$ \varphi_{\xi_\lambda} (r)= r+ \frac{\psi ( (1-r)\beta(\lambda))}{\beta(\lambda)\,  \psi^\prime \! (\beta(\lambda) )} \, , \quad c_\lambda= \psi^\prime \! (\beta(\lambda) ) \quad \! \textrm{and} \quad \! \varphi_{\mu_\lambda} (r)= \int_{(0, \infty)} e^{-(1-r)y \beta(\lambda)} \, \varrho(dy) . $$
These offspring distributions already appeared in the more specific context of \cite{DuLG,DuWi1}. 
Of particular interest are cases where the offspring distributions $\xi_\lambda$ are all equal to a certain $\xi$ not depending on $\lambda\ge 0$. We easily see that this exactly corresponds to the stable cases where $\psi(\lambda)= \lambda^\gamma$, for a certain $\gamma \in (1, 2]$ and thus $\varphi_\xi (r)= r+ \frac{1}{\gamma} (1-r)^\gamma$. We call these laws the $\gamma$-stable offspring distributions. The Brownian case corresponds to the critical binary offspring distribution ($\gamma=2$).  They appear in previous work \cite{LJ91, Ne2} (and in a slightly different form in \cite{DuLG, Va, Ya}). 
  
Observe that $\lambda \mapsto \Gamma (\widetilde{\cF}_\lambda)$ is non-decreasing almost surely, since the growth process is 
$\preceq$-non-decreasing. Then $ \lim_{\lambda \rightarrow \infty}    \Gamma (\widetilde{\cF}_\lambda)$ exists in $[0, \infty]$ and standard branching process arguments give that 
$\bP(  \lim_{\lambda \rightarrow \infty}  \Gamma (\widetilde{\cF}_\lambda) < \infty) >0$ iff  
$\psi$ satisfies $\int^\infty dr / \psi (r) < \infty$. This is a necessary and sufficient condition for the growth process to converge in $(\bT,\dgh)$, almost surely. Namely, Theorem \ref{asconvgrowth} asserts that   
there exists a random CLCR real tree $\widetilde{\cF}$ such that 
$$\lim_{\lambda \rightarrow \infty} \dgh \big( \widetilde{\cF}_\lambda, \widetilde{\cF}\big)= 0  \qquad\textrm{almost surely}.$$
The limiting tree $\widetilde{\cF}$ is a $\psi$-L\'evy forest. Moreover, the branching processes associated to $\widetilde{\cF}_\lambda$, $\lambda\ge 0$, also converge almost surely. Namely, for any $a\in [0, \infty)$, we get 
$$ \frac{1}{{\beta (\lambda)}} \, \# \, \{ \s \in \cF_\lambda : d(\rho, \s) = a \}\;   \underset{\lambda \rightarrow \infty}{-\!\!\!-\!\!\!-\!\!\! \longrightarrow} \;  Z_a\qquad\mbox{almost surely,}$$   
where $(Z_a)_{a \in [0,  \infty)}$ is a continuous-state branching process with branching mechanism 
$\psi$ and whose initial value $Z_0$ has law $\varrho$ (see Section \ref{Levytreeconvsec} for more details). As an application of these results, Theorem \ref{Levyfochar} provides a nice characterisation of L\'evy forests via leaf-length erasure. This result is used in the proof of the invariance principles for discrete Galton-Watson trees that follow in Section \ref{Invprincsec}. These limit theorems are general and they notably include the super-critical cases. Invariance principles for critical and sub-critical trees were obtained in \cite{DuLG} using different arguments and a different formalism.

  This paper is organised as follows. In Section \ref{realtrees}, we recall basic definitions concerning real trees and the Gromov-Hausdorff metric. In Section \ref{funcsubsec}, we 
discuss various operations on real trees. 
The technical details concerning the measurability of such operations are postponed 
to an appendix. Section \ref{convcritsec} is devoted to the statement and the proof of the main 
general tightness results and convergence theorems for real trees. 
Section \ref{GWrealtreesec} discusses an intrinsic definition of Galton-Watson real trees, and hereditary reduction is studied in Section \ref{realhereditsec}. In Section \ref{defgrowthsec}, we define and study growth processes of Galton-Watson forests,  whose limits are discussed in Section \ref{Levytreeconvsec}. Section \ref{Invprincsec} is devoted to invariance principles for rescaled discrete Galton-Watson trees with unit edge lengths.

\section{Preliminary results.}
\label{prelim}
 \subsection{Real trees.} 
\label{realtrees}
 Real trees form a class of loop-free length spaces, which turn out to be limits of many discrete trees. 
More precisely, we say that a metric space $(T,d)$ is 
{\it a real tree} if it satisfies the following conditions: 
\begin{itemize} 
\item  for all $\sigma , \sigma' \in T$, there exists a unique  
    isometry $f_{\sigma , \sigma'}\colon[0,d(\sigma,\sigma')]\rightarrow T$ such that 
    $f_{\s,\s'}(0)=\s$ and $f_{\sigma,\sigma'}(d(\sigma,\sigma'))=\s'$; we set $\lgeo \s, \s'\rgeo :=f_{\s,\s'}([0, d(\s,\s')]) $, which is the geodesic joining $\sigma$ to $\sigma^\prime$;

\item if $q\colon[0,1]\rightarrow T$ is continuous injective, we have $q([0,1])= \lgeo q(0), q(1) \rgeo$.
\end{itemize}

Let $\rho$ be a distinguished point of $T$, which is viewed as the {\it root} of $T$. Then $(T,d,\rho)$ is called a {\it rooted real tree}. We also denote by 
$\rgeo \s,\s' \rgeo $, $\lgeo \s,\s'\lgeo$ and 
$\rgeo \s,\s'\lgeo$ the respective images of $(0,d(\s, \s')] $, $[0,d(\s, \s')) $ and 
$(0,d(\s, \s')) $ under $f_{\s,\s'}$. We view the tree $T$ as the family tree of a population whose progenitor is $\rho$. For any $\sigma , \s^\prime \in T$, their most recent common ancestor is then the unique point denoted by $\sigma \wedge \sigma^\prime$ such that $\lgeo \rho , \sigma \rgeo \cap \lgeo \rho , \sigma^\prime \rgeo = \lgeo \rho , \sigma \wedge \sigma^\prime \rgeo $. Thus, 
\begin{equation}
\label{branchheight}\textstyle 
\forall \s , \s^\prime  \in T \, , \quad d( \rho , \s \wedge \s^\prime ) = \frac{1}{2} \left(d( \rho , \s)+d( \rho , \s^\prime)-d( \s , \s^\prime)  \right) \; .
\end{equation}
There is a nice characterization 
of real trees in terms of the 
{\it four points condition}: a connected metric space $(X,d)$ is a real tree iff for all
$\s_1,\s_2,\s_3,\s_4\in X$, 
\begin{equation}
\label{fourpoint}
d(\s_1, \s_2) +d(\s_3, \s_4) \leq 
(d(\s_1, \s_3) +d(\s_2, \s_4))\vee (d(\s_3, \s_2) +d(\s_1, \s_4)). 
\end{equation}
We refer to \cite{Dress84,DMT96,DT96} for general
results concerning real trees, to \cite{Pau88, Pau89} for applications of real 
trees to group theory and to \cite{DuLG, DuLG2, Ev00, EvPitWin, EvWin} 
and \cite{LyoHam} 
for probabilistic use of real 
trees.  

\paragraph{Gromov-Hausdorff distance on the space of complete locally compact real trees.} We say that two pointed metric spaces 
$(X_1,d_1,\rho_1)$ and $(X_2,d_2,\rho_2)$ 
are {\it equivalent} iff there exists a \em pointed isometry\em, i.e.\ an isometry 
$f$ from $X_1$ {\it onto} $X_2$ such that $f(\rho_1)=\rho_2$. 
The Gromov-Hausdorff distance of 
two {\it pointed compact} metric spaces 
$(X_1,d_1,\rho_1)$ and $(X_2,d_2,\rho_2)$ is given by the following.  
\begin{equation}
\label{dghcpct}
 \dgh_{{\rm cpct}} (X_1,X_2)= \inf 
\left\{ \, d_{{\rm Haus}}(f_1(X_1),f_2(X_2))\vee d (f_1(\rho_1),
f_2(\rho_2))  \, \right\} \; .
\end{equation}
Here the infimum is taken over all $(f_1, f_2, (E, d))$, where $(E, d)$ is a metric space, where $f_i\colon  X_i \rightarrow E$, $i= 1 ,2$, are isometric embeddings and where $d_{{\rm Haus}}$ stands for the Hausdorff distance 
on the set of 
compact subsets of $E$. Observe that $ \dgh_{{\rm cpct}}$ 
only depends on the isometry classes of the $X_i$. It induces a metric on the set of isometry classes 
of all pointed compact metric spaces (see \cite{Gro}).

  We then denote by $\bT_{{\rm cpct}}$ the set of pointed isometry classes of compact rooted real trees. 
Evans, Pitman and Winter \cite{EvPitWin} 
showed that $\bT_{{\rm cpct}}$ equipped with the Gromov-Hausdorff 
distance $ \dgh_{{\rm cpct}}$ is a complete and separable 
metric space.

We then denote by $\bT$ {\it the set of pointed isometry classes of complete locally compact rooted real trees}. The Gromov-Haudorff distance is extended to $\bT$ in the following standard way.  
Let $(T_1, d_1, \rho_1)$ and $(T_2, d_2, \rho_2)$ be two complete locally compact rooted real trees. 
Recall that the Hopf-Rinow theorem (see \cite[Chapter 1]{Gro}) implies that the closed balls of $T_1$ and $T_2$ 
are compact sets (note that this entails that $T_1$ and $T_2$ are separable). We then set 
\beqs \dgh(T_1,T_2)=\sum_{k\ge 1}2^{-k}\dgh_{\rm cpct}(\overline{B}_{T_1}(\rho_1,k),\overline{B}_{T_2}(\rho_2,k)),
\eeqs
where for any $i\in \{ 1, 2\}$, $\overline{B}_{T_i}(\rho_i,k)$ stands for the closed (compact) ball with center $\rho_i$ and radius $k$ in the locally compact rooted real tree $(T_i,d_i,\rho_i)$. 
Clearly, $\dgh$ only depends on the isometry classes of $T_1$ and $T_2$. It defines a metric on 
$\bT$ and $(\bT, \dgh)$ is Polish (see \cite[Proposition 3.4]{DuWi1}).  

\medskip

\noi
{\it Notation and convention.} We shorten {\it Complete Locally Compact Rooted} real tree to {\it CLCR real tree}. Let $(T, d , \rho)$ be a CLCR real tree. We denote by $\widetilde{T} \in \bT$ its pointed isometry class. We shall denote by $\pnt$ the pointed isometry class of a point tree 
$(\{\rho\},d,\rho)$.

\paragraph{Isometrical embeddings of CLCR real trees in pointed Polish spaces.}
\label{ellunemb}
In this paper, we deal with growing families of real trees that may be embedded into a given space for technical reasons. More precisely, let $(E, d, \rho)$ be a pointed Polish space. We introduce the following set. 
\begin{equation}
\label{TEdef}
\bT_E= \left\{ T \subseteq E \colon  \textrm{$\rho \in T$ and $(T, d, \rho)$ is a CLCR real tree} \right\}\; .
\end{equation}
Let us denote by $d_{{\rm Haus}}$ the Hausdorff distance on compact
subsets of $E$. Then, for any 
$T, T^\prime \in \bT_E $, we define 
$$ {\bf d}_E (T, T^\prime)=\sum_{k\ge 1} 2^{-k} d_{{\rm Haus}} \left(\overline{B}_T(\rho , k) , \overline{B}_{T^\prime}(\rho, k)  \right) . $$
Clearly, ${\bf d}_E$ is a distance on $\bT_E$ and we have 
\begin{equation}
\label{dcontrol}
\dgh (T, T^\prime) \leq {\bf d}_E (T, T^\prime)\; .
\end{equation}
Note that $(E, d)$ may be ``small'' and that the set $\{ \widetilde{T} \, ; T \in \bT_E\}$ may be strictly included in $\bT$. However, following Aldous's idea (see \cite{Al1}), it is possible to 
embed all CLCR trees in the vector space 
$\Ell_1$ of summable real-valued sequences equipped with the $\lVert \cdot  \rVert_1$-norm and where $0$ is the distinguished point. \cite[Proposition 3.7]{DuWi1} shows that every element of $\bT$ has a representative in $\btl$. Namely $\bT = \{ \widetilde{T}\, ; T \in \btl \}$. 
Moreover, \cite[Proposition 3.6]{DuWi1} asserts that $(\btl , \dl_{\Ell_1})$ is a Polish space. The arguments of the proof can be directly adapted to the more general case to prove that 
$(\bT_E , {\bf d}_E)$ is a Polish space.

\paragraph{Approximation by real trees with edge lengths.} Let $(T, d,\rho)$ be a rooted real tree. 
For all $\s\in T$ we denote by ${\rm n}(\s,T)$ the {\it degree} of $\s$, namely the 
(possibly infinite) 
number of connected components of $T\setminus\{\s\}$. We also introduce the following 
$${\rm Lf}(T)=\{ \s\in T\setminus\{ \rho \}\colon  {\rm n}(\s,T)=1\} \quad 
{\rm and } 
\quad {\rm Br}(T)=\{ \s\in T\setminus\{ \rho \}\colon  {\rm n}(\s,T)\geq 3\}$$
that are respectively the set of the {\it leaves} of $T$ and 
the set of {\it branch points} of $T$. Note that the root is neither considered as a leaf nor as a branch point. Recall from \cite{DuWi1} that the set of branch points of a CLCR real tree is at most countable. 
\begin{definition}
\label{realedgetrees} 
A rooted real tree $(T,d,\rho)$ is called a {\it real tree with edge lengths}
if it is complete and if
\begin{equation}
\label{degsumfin}
\forall R \in (0, \infty) \, , \quad  {\rm n}(\rho,T)+ \sum  {\rm n}(\s,T) < \infty \; ,
\end{equation}
where the sum is over all the branch points $\s \in {\rm Br} (T)$ such that $d(\rho, \s) \leq R$. \cq 
\end{definition}
An equivalent definition is the following: 
a complete real tree $(T,d,\rho)$ is real tree with edge lengths iff 
\begin{enumerate}
\item[(a)] the degree of branch points is bounded in every ball of finite radius; 
\item[(b)] every ball contains a finite number of branch points and a finite number of leaves. 
\end{enumerate}
Indeed, the only non-trivial point to check is that (\ref{degsumfin}) implies that 
every ball contains a finite number of leaves. We argue by contradiction: fix $R \in (0, \infty)$ and suppose that there exists a sequence 
$(\s_n)_{n \in \bN}$ of distinct leaves such that $d(\rho, \s_n ) \leq R$. Let $B$ be the set of points 
$\gamma $ such that $\gamma = \s_i \wedge \s_j $ for infinitely many pairs of integers $i<j$. Clearly, 
$B \subseteq {\rm Br} (T) \cap B_T (\rho, R)$ and (\ref{degsumfin}) implies that $B$ is finite and non-empty.   
Let $\beta \in B$ be such that $d(\rho, \beta)= \max_{\gamma \in B} d(\rho, \gamma)$.   
By (\ref{degsumfin}), $\beta $ has finite degree (say $n+1$); let $T_1, \ldots , T_n$ be the open connected components of $T\backslash \{ \beta \}$ that do not contain the root. One of these connected components $T_i$ has to contain infinitely many terms of the sequence $(\s_n)_{n \in \bN}$. Then there is 
$\beta^\prime \in T_i \cap B$ such that $d(\rho, \beta^\prime) > d(\rho, \beta)$, which is not possible. 

\medskip

Let $(T, d, \rho)$ be a real tree with edge lengths. Denote the connected components of $T\setminus(\Br(T)\cup\{ \rho\})$ by $C_i$, $i\in I$. The components $C_i$, $i\in I$, are called the {\it edges} of $T$. They all are isometric to intervals of the real line. Note that their endpoints are leaves, branch points or the root. 
By (a) and (b), for all $R \in (0, \infty)$, only finitely many of such edges have points at distance 
less than $R$ from the root. This implies that the closed ball with center $\rho$ and radius $R$ is compact. Thus, $(T, d,\rho)$ is locally compact. Let us denote by $\bT_{\rm edge}$ the set of pointed isometry classes of the real trees with edge lengths. Then, we get 
$$ \bT_{{\rm edge}} \subset \bT \; .$$

We next introduce {\it leaf-length erasure} (also called trimming) that is used throughout the paper.  
Let $h\in (0, \infty)$ and let $(T,d,\rho)$ be a CLCR real tree. We set 
\begin{equation} 
\label{hleafera}
R_h(T)=\{\rho\}\cup\{\sigma\in T\colon \exists\sigma^\prime\in T\mbox{ such that } \sigma\in\lgeo\rho,\sigma^\prime\rgeo\mbox{ and } 
                          d(\sigma,\sigma^\prime)\ge h\}.
\end{equation}
Clearly, $R_h(T)$ is path-connected and it is easy to prove that it is a closed subset of $T$. Then, 
the four points condition implies that 
$(R_h(T),d,\rho)$ is a CLCR real tree. 
We call it the 
{\it $h$-leaf-length erased tree associated with $T$}. Note that its isometry class only depends on that of $T$ so $R_h$ can be defined from $\bT$ to $\bT$. 

  Leaf-length erasure was first introduced by Kesten \cite{Kes86} for discrete trees 
and further studied and applied by many others \cite{Ne2,NP,LG86,LJ91,EvPitWin}. 
In the following lemma we sum up the various properties of the leaf-length erasure operator $R_h$ that we shall use in this paper; it is a straightforward extension of the same result for compact trees that is due to Evans, Pitman and Winter \cite[Lemma 2.6]{EvPitWin}.  
\begin{lemma}
\label{eraseprop} 

\begin{enumerate}
\item[$\rm(i)$] For every $h \in (0, \infty)$ and every CLCR real tree $(T, d, \rho)$, 
$(R_h(T),d,\rho)$ is a real tree with edge lengths. 
\item[$\rm(ii)$] For every $h \in [0, \infty) $, $R_h$ is $\dgh$-continuous. 
\item[$\rm(iii)$] $R_{h+h^\prime}=R_h\circ R_{h^\prime}$, $h,h^\prime \in [0, \infty)$. 
\item[$\rm(iv)$] For every $h\in [0, \infty)$ and for 
every CLCR real tree $(T,d,\rho)$, we have $\dgh (R_h(T),T)\le h$. 
\end{enumerate}
\end{lemma}
\begin{pf} Since $T$ is locally compact, we immediately get that $ \n(\sigma , R_h(T)) < \infty$, for every $\sigma \in R_h (T)$. To prove $\rm(i)$, it remains to prove that $\Br(R_h(T))$ has no limit points. Suppose that there exists a sequence $(\s_n)_{n \in \bN}$ of distinct branching
points of $R_h(T)$ converging to $\sigma\in R_h(T)$. Then it is easy to see that the closed ball 
$\overline{B}_T(\sigma,2h)$ is not locally compact, 
which is a contradiction. Points $\rm(iii)$ and $\rm(iv)$ are easy consequences of 
\cite[Lemma 2.6]{EvPitWin}. Let us briefly explain $\rm(ii)$: for every CLCR real tree $(T,d,\rho)$ and for every $r>0$, we set $T(h,r)=R_h(\overline{B}_T(\rho,r+2h))$. Then observe that $ \overline{B}_{T(h,r)}(\rho,r)=R_h(T)\cap \overline{B}_T(\rho,r)$. This entails $\rm(ii)$ by \cite[Lemma 2.6 \rm(i)]{EvPitWin}, which asserts that $R_h$ is $\dgh_{\rm cpct}$-continuous on $\bT_{\rm cpct}$.
\end{pf}

\subsection{Specific functions on real trees.}
\label{funcsubsec}

We introduce functions of real trees such as the total height, the profile and various procedures 
that allow to split or to paste trees. Continuity or measurability of such functions 
is quite expected, however some of the proofs are technical. Thus, to ease the reading, we postpone them to Appendix \ref{appendixA}. 

\paragraph{The total height.} Let $(T,d,\rho)$ be a CLCR real tree. 
The total height of $T$ is given by 
$$\Gamma (T)= \sup_{\sigma \in T}
 d(\rho, \sigma) \in [0, \infty]\; .$$ 
Note that $\Gamma (T)$ only depends on the pointed isometry class $\widetilde{T}$ of $(T,d,\rho)$: it induces a function on $\bT$ that is denoted in the same way. It is easy to check that $\Gamma$ is 
$\dgh$-continuous and note that 
\begin{equation}
\label{heighterase}
 \forall h \in [0, \infty) \; , \quad  \Gamma \circ R_h = (\Gamma -h)_+ \; .
\end{equation}
\paragraph{The $h$-erased profile.} For every $a \in [0, \infty)$ and every $h\in (0, \infty)$, we set 
\begin{equation}
\label{Zhdef}
 Z_a^{(h)}(T)=\# \, \{ \,  \sigma\in R_h(T)\colon d(\rho,\sigma) = a \, \} .
\end{equation}
Since $R_h(T)$ is a real tree with edge lengths, it is easy to check that $a \mapsto  Z_a^{(h)}(T)$ is 
an $\bN$-valued function that is left-continuous with right limits. We call the process $a\mapsto Z_a^{(h)} (T)$ the $h$-erased profile. We shall denote by $Z^{(h)}_{a+} (T)$ its right limit at $a$.  
Note that $h \mapsto Z_a^{(h)}(T)$ is non-increasing. Since $Z_a^{(h)}(T)$ only depends on the pointed isometry class of $T$, it induces a function on $\bT$ that is denoted in the same way. 

\paragraph{Splitting measures.} Let us recall that a point measure on $\bT$ is a measure 
of the form $\sum_{i\in I} \delta_{\widetilde{T}_i}$. We introduce the following set of point measures
$$\textstyle \ccM(\bT)= \left\{ \sum_{i \in I} \delta_{\widetilde{T}_i} \textrm{ point measure on $\bT$}\colon   
\forall h \in(0, \infty), \;   \# \{ i\in  I\colon \Gamma (\widetilde{T}_i) > h \} < \infty  \right\}  . $$
Note that $\ccM (\bT)$ contains the null measure (for which $I$ is taken empty). 
We use the following standard notation on point measures: for every $M = \sum_{{i \in I}} \delta_{\widetilde{T}_i}$ in $ \ccM (\bT)$ and for every function $G\colon \bT \rightarrow [0, \infty )$, we set 
$ \langle M, G \rangle =  
\sum_{{i \in I}} G( \widetilde{T}_i)  $. 
We shall also denote by $\langle M\rangle $ the total mass of $M$: we then check that 
$\langle M\rangle= \langle M, \un_{\bT} \rangle=  \#I$. 

  Since $ (\bT, \dgh)$ is not locally compact, the vague topology on $\ccM(\bT)$ is not useful. We rather equip $\ccM(\bT)$ with the sigma-field $\cG_{\ccM (\bT)}$ generated by the functions $M\in \ccM (\bT)  \mapsto \langle M, G\rangle$, for
$G\colon  \bT\rightarrow [0, \infty)$ Borel-measurable.

We also define the {\it trees above level $a$} as follows. 
Let $(T, d, \rho)$ be a CLCR real tree. 
Denote by $T^\circ_{i} (a)$, $i \in I(a)$, the connected components of the (possibly empty) open set 
$T \backslash \overline{B}_T (\rho, a)$; for every $i \in  I(a)$, denote by $T_i (a)$ the closure of $T^o_i (a)$ in $T$ and denote by $\s_i (a)$ the unique vertex such that $T_i (a)= T^\circ_i (a)  \cup \{ \s_i (a) \}$. Observe that 
$d(\rho , \sigma_i (a))= a$. Clearly, the trees 
$(T_i (a) , d , \s_i (a))$, $ i \in I(a)$, are CLCR real trees 
whose pointed isometry classes are denoted by $\widetilde{T}_i (a)$, $i\in I(a)$. 
We then set 
\begin{equation}
\label{asplitdef}
\cM_a  (T)= \sum_{i\in I(a)} \delta_{\widetilde{T}_i (a)}  \in \ccM (\bT) 
\quad {\rm and} \quad Z^+_a (T) = \langle \cM_a (T) \rangle = \# I (a)  \in \bN \cup \{ \infty \}\; .
\end{equation}
The fact that $\cM_a (T) \in \ccM (\bT)$ is an easy consequence of the local compactness of $T$. The measure $\cM_a (T)$ is called the {\it splitting measure of $T$ at level $a$} and $Z^+_a (T)$ is called the 
{\it right profile of $T$}. Note that these functions only depend on $\widetilde{T}$: they induce functions on $\bT$ that are denoted in the same way. It is easy to check that for every $a \in [0, \infty)$ and every $h \in (0, \infty)$, 
\begin{equation} 
\label{profconn}
Z^{(h)}_{a+} = Z^+_a \circ R_h \quad {\rm and} \quad Z^+_a = \lim_{h \rightarrow 0+} Z^{(h)}_{a+} \; .
\end{equation}

\paragraph{Grafting real trees.} Let us explain how to graft real trees on the vertices of another real tree.
Let $(T, d, \rho)$ be a rooted real tree, let 
$(T_i, d_i, \rho_i)$, $i \in I$, be 
a family of real trees and let $\s_i \in T $, $ i\in I$, be a
collection of vertices of $T$. 
We then set 
$$ T^\prime =T \; \sqcup   \bigsqcup_{i\in I} T_i\setminus{\{ \rho_i\}} \; ,$$ 
where $\sqcup$ stands for the disjoint union. We next define a metric $d^\prime$ on $T^\prime \times T^\prime$ as follows:   
$d^\prime$ coincides with $d$ on $T\times T$; assume that 
$\s\in  T_i\setminus{\{ \rho_i\}}$; if $\s^\prime \in T$, 
then we set $d'(\s,\s^\prime)=d_i(\s,\rho_i) + d(\s_i,\s^\prime) $; if $\s^\prime \in T_j\setminus{\{ \rho_j\}}$ with $j\neq i$, then we set $d^\prime(\s,\s^\prime)=  d_i(\s,\rho_i)+d(\s_i,\s_j)+d_j(\s^\prime,\rho_j)$; if $\s^\prime \in  T_i\setminus{\{ \rho_i\}}$, then we set
$ d'(\s,\s^\prime)=d_i(\s,\s^\prime)$. 
It is easy to prove that $(T^\prime,d^\prime,\rho^\prime :=\rho)$ is a rooted real tree and we use the  notation 
\begin{equation}
\label{pastnota}
 (T^\prime, d^\prime,\rho^\prime)= T \circledast_{i\in I} (\s_i, T_i) \; .
 \end{equation}
We shall extensively use the following special case: we assume that $T$ reduces to its root $T= \{ \rho \}$: in this case $\s_i = \rho$ and we use the specific notation 
$$\circledast_{i\in I} T_i:=  \{ \rho \} \circledast_{i\in I} (\rho, T_i)   \; , $$ 
with the convention that $  \circledast_{i\in I} T_i= \{ \rho \}$ if $I$ is empty. In words,  
$\circledast_{i\in I} T_i $ is a tree obtained by pasting the trees $T_i$, $i\in I$, at their roots. 

  We now assume that for every $i\in I$, $(T_i , d_i, \rho_i)$ is a CLCR real tree and we assume that for every $h \in (0, \infty)$, $\# \{ i\in I \colon  \Gamma (T_i) > h \} < \infty$. This implies that $\{ i\in I \colon  T_i \neq \{ \rho_i \} \}$ is countable.  
It is easy to check first that  
$\circledast_{i\in I} T_i $ is a CLCR real tree. Note that its pointed isometry class only depends on the pointed isometry classes of $T_i$, $i\in I$. We denote by $\circledast_{i\in I}  \widetilde{T}_i$ the pointed isometry class of $\circledast_{i\in I} T_i $.

  This induces a function $\paste\colon \ccM (\bT) \rightarrow \bT$, that is defined as follows: for every 
$M = \sum_{i\in I} \delta_{\widetilde{T}_i} \in \ccM (\bT)$, we set 
\begin{equation}
\label{pastedef}
\paste( M )= \circledast_{i\in I}  \widetilde{T}_i \; , 
\end{equation}
with the convention that $\paste (M)= \pnt$ if $M$ is the null measure. It is easy to check that 
\begin{equation}
\label{splitpaste}
 \forall \widetilde{T} \in \bT\, , \; \paste (\cM_0 (\widetilde{T} )) = \widetilde{T} \quad {\rm and} \quad \forall M=\sum_{i\in I}\delta_{\widetilde{T}_i} \in \ccM (\bT)\, , \; \cM_0 ( \paste (M))= 
 \sum_{i\in I\colon\widetilde{T}_i\neq\pnt}\delta_{\widetilde{T}_i}.
 \end{equation}
\paragraph{The tree above a given level.} Let $(T, d, \rho)$ be a CLCR real tree and let $a\in [0, \infty)$.  
Recall the definition $(T_i (a), d, \sigma_i (a))$, $i\in I(a)$, of the subtrees of $T$ above $a$. We then set 
\begin{equation}
\label{abovedef}
 \abv (a, T) =  \circledast_{i\in I(a)} T_i (a) \; .
\end{equation} 
The tree  $ \abv (a, T)$ is called the {\it tree above level $a$}. 
By definition, the pointed isometry class of  $\abv (a, T)$ is $\circledast_{i\in I(a)}\widetilde{ T}_i (a)$ and it only depends on $\widetilde{T}$. This induces a function from $\bT$ to $\bT$ that is denoted in the same way and we check that 
\begin{equation}
 \cM_0 \circ \abv (a, \cdot) = \cM_a \quad {\rm and} \quad \abv (a, \cdot) 
 = \paste \circ \cM_a \; .
\end{equation} 
We also denote by $\blw (a , T)$ the tree $(\overline{B}_T (\rho , a) , d , \rho)$, that is the {\it tree below the level $a$}. Its pointed isometry class only depends on $\widetilde{T}$: this induces a function 
from $\bT$ to $\bT$ that is denoted in the same way. 
\begin{lemma}
\label{abovebelow} The map $(a, \widetilde{T}) \in [0, \infty) \times \bT \longmapsto \abv (a, \widetilde{T})$ is jointly continuous. The same holds true for $\blw$. 
\end{lemma}
\begin{pf} See Appendix \ref{AbvBlwpfsec}. 
\end{pf}
\begin{lemma}
\label{measuraZht} Fix $a\in [0, \infty)$ and $h\in (0, \infty)$. The following assertions hold true. 
\begin{enumerate}
\item[\rm(i)] $\cM_a\colon  \bT \rightarrow \ccM (\bT)$ is measurable and so are $Z^+_a\colon  \bT \rightarrow \bN \cup \{ \infty \} $ and $Z_a^{(h)} \colon  \bT \rightarrow \bN$. 
\item[\rm(ii)]  $\paste\colon  \ccM  (\bT) \rightarrow \bT$ is measurable. 
\end{enumerate}
\end{lemma}
\begin{pf} See Appendix \ref{mesusec}.
\end{pf}
\paragraph{The functions $D$, $\vartheta$ and ${\bf k}$.} Let $(T, d , \rho)$ be a CLCR real tree. 
Recall that $\rho \notin  {\rm Br} (T) \cup {\rm Lf} (T)$, by convention. We next set 
\begin{equation}
\label{distbranch}
D (T)= \inf \left\{ d(\rho , \sigma) \; ; \;  \sigma \in {\rm Br} (T) \cup {\rm Lf} (T) \, \right\} \end{equation}
with the convention that $\inf \emptyset = \infty$. Namely, it is important to note that 
$D(T)= \infty$ iff either $T$ is reduced to a point or $T$ 
is equivalent to a finite number of half lines pasted at their finite endpoint. Recall from (\ref{asplitdef}) the definition of $Z^+_a $ and from (\ref{abovedef}) that of $\abv (a, \cdot) $. If $D(T) < \infty$, we set  
\begin{equation}
\label{firstabove}
\vartheta T: = \abv (D(T)\, , T) \quad {\rm and} \quad {\bf k} (T)= Z^+_0 ( \vartheta T) \; . 
\end{equation}
If $D(T)= \infty$, we set $\vartheta T = \{ \rho \}$ and ${\bf k} (T) = 0$. 
Observe that $D (T)$, ${\bf k} (T)$ and the pointed isometry class of $\vartheta T$ only depend on the pointed isometry class of $T$. 
So they can be viewed as functions on $\bT$. With a slight abuse of language, $D (T)$ is viewed as the distance from the first branch point, $\vartheta T$ as the tree above the first branch point 
and ${\bf k} (T) +1$ as the degree of the first branch point.   
\begin{lemma}
\label{measDthetak} We have $D=\lim_{h\rightarrow 0}D\circ R_h$. Moreover, the functions 
$D\colon  \bT \rightarrow [0, \infty]$, ${\bf k}\colon  \bT \rightarrow \bN \cup \{ \infty \}$ and $\vartheta \colon  \bT \rightarrow \bT$ are measurable. 
\end{lemma}
\begin{pf} See Appendix \ref{measDthetakpf}.
\end{pf}
The functions $D$, $\vartheta$ and ${\bf k}$ are useful only when applied to real trees with edge 
lengths and they allow to characterise $\bT_{\rm edge}$. More precisely, for every $n \in \bN$, we 
recursively define $\vartheta_n $ and $D_n $ by setting $ \vartheta_{n+1} = \vartheta \circ \vartheta_n$ (where $\vartheta_0$ is taken as the identity map on $\bT$) 
and $D_{n} = D \circ \vartheta_n$. Thus, $D=D_0$ and $\vartheta = \vartheta_1$. 
\begin{lemma}
\label{edgechar} 
We have that $\bT_{\rm edge}=\left\{ \widetilde{T}\in\bT\colon \sum_{n\in \bN } D_n(\widetilde{T})=\infty \right\}$ is a Borel subset of $\bT$.
\end{lemma} 
\begin{pf} See Appendix \ref{edgecharpf}.
\end{pf}

\subsection{Convergence criteria} 
\label{convcritsec}
This section is devoted to new general 
convergence theorems for complete locally compact rooted real trees (CLCR real trees): we first state a general 
tightness criterion involving the $h$-erased profiles and we also state 
an almost sure convergence criterion 
for sequences of trees that ``grow'' in a broad sense that we will make precise later.

We use the following notation. Let $(T, d, \rho)$ be a rooted real tree, not necessarily locally compact. Then, for every $h \in (0, \infty)$ the definition (\ref{hleafera}) of $h$-leaf-erased tree $(R_h(T), d, \rho)$ makes sense and it defines a path-connected subset of $T$. It is therefore a 
rooted real tree. For every $a\in [0,\infty)$, we define $Z^{(h)}_a (T) $ as in (\ref{Zhdef}). Note that 
this quantity may be infinite if $T$ is not locally compact. We denote by $N_T(h,r)$ 
the (possibly infinite) minimal number of open balls with radius $h$ 
that are necessary to cover $\overline{B}_T(\rho,r)$. 
\begin{lemma}
\label{keycompact} Let $(T, d, \rho)$ be a complete rooted real tree (not necessarily locally compact). 
Then, 
$$ \forall a, h \in (0, \infty)\, , \; \forall r \in (a+h , \infty ) \, , \quad  Z^{(h)}_a (T) \leq N_T(h,r) \; .$$
Moreover, let $\cD \subseteq (0, \infty)$ be a dense subset. Then, for every $h, r \in (0, \infty)$, there is a finite subset $\cD_{h, r} \subset \cD \cap [0, r]$ such that 
\begin{equation}
\label{Zhentro}
N_T (h,r) \leq 1+ \sum_{a\in \cD_{h,r}} Z^{(h/3)}_a (T) \; .
\end{equation}
Consequently, $(T, d, \rho)$ is locally compact iff there exists a dense subset $\cD \subseteq (0, \infty)$  and a sequence $(h_p)_{p\in \bN }$ decreasing to $0$ such that 
\begin{equation}
\label{loccompcri}
\forall a\in \cD \, , \; \forall p \in \bN \, , \quad Z^{(h_p)}_a (T) < \infty. 
\end{equation}
\end{lemma}
\begin{pf} 
Let $a, h \in (0, \infty)$ and $r > a+h$. We fix $L= \{ \sigma\in R_h(T)\colon d(\rho,\sigma)= a\}$ so that $Z^{(h)}_a (T)= \# L$. Note that $N_T (h, r) \geq 1$. Thus, (\ref{Zhentro}) holds if $\# L \leq 1$. Now assume that $\# L \geq 2$. With every $\sigma\in L$, we can associate $\zeta(\sigma)\in T$ such that $d(\rho,\zeta(\sigma))=a+h$ and $\sigma\in\lgeo\rho,\zeta(\sigma)\rgeo$. Note that $\zeta (\sigma) \in \overline{B}_T (\rho, r)$. Suppose that 
$\sigma$ and $\sigma^\prime$ are two distinct points of $L$. Recall 
the notation $\sigma \wedge \sigma^\prime$ for the most recent common ancestor of $\sigma$ and $\sigma^\prime$. Then, $d(\zeta(\sigma),\zeta(\sigma^\prime))>2h$ because $\zeta (\sigma) \wedge \zeta(\sigma^\prime) = \sigma\wedge \sigma^\prime$ is below level 
$a$ in $T$. Thus, $ Z^{(h)}_a (T)= \# L \leq N_T(h,r)$.  

Let $h \in (0, \infty)$ and let $\cD\subset(0,\infty)$ be dense. We can find an increasing sequence $a(k)\in \cD$, $k\in \bN$, such that $a(0) < h/2$ and $h \leq a(k+1)-a(k) < 3h/2$. For every $k\in\bN$, we then set 
$L_k=\{\sigma\in R_h(T)\colon d(\rho,\sigma)= a(k)\}$. Let $\sigma\in T$ be such that $d(\rho,\sigma) \geq 2h$. There exists $k\in\bN$ such that $a(k+1)\le d(\rho,\sigma)<a(k+2)$.
Thus, $\lgeo\rho,\sigma\rgeo\cap L_k=\{\sigma^\prime\}$ and $d(\sigma,\sigma^\prime)< a(k+2)-a(k) <3h $. This implies $N_T(3h,r)\le 1+\sum_{0\leq k \leq \lfloor r/h\rfloor}Z_{a(k)}^{(h)}(T)$, 
where $\lfloor r/h\rfloor$ stands for the integer part of $r/h$. This easily entails (\ref{Zhentro}).

 If $T$ is locally compact, we already noticed that $Z^{(h)}_a (T) < \infty$, for every $a, h\in (0, \infty)$. Conversely, suppose that (\ref{loccompcri}) holds true. Note that $h \mapsto Z^{(h)}_a (T)$ is non-increasing. Then, (\ref{Zhentro}) applies and for every $r $, and $(\overline{B}_T(\rho,r), d)$ is a uniformly bounded complete metric space. It is therefore compact. This implies that $T$ is locally compact. 
\end{pf}
Note that $N_T(h,r)$ only depends on the isometry class $\widetilde{T}$ of $T$, which justifies the notation $N_{\widetilde{T}}(h,r)$.  Let $\cC \subseteq \bT$. 
As a consequence of a general result on compactness with respect to the pointed Gromov-Hausdorff metric $\dgh_{{\rm cpct}}$ (see \cite{Gro} or \cite[Theorem 8.1.10]{BuBu}), we have the following result.  
\begin{equation}
\label{Grogene}
\textrm{The $\dgh$-closure of $\cC$ is compact} \; \Longleftrightarrow \; \forall h,r\in (0, \infty) \, , \quad 
\sup_{\widetilde{T}\in\cC} N_{\widetilde{T}}(h,r)<\infty .
\end{equation}
Thus, Lemma \ref{keycompact} immediately entails  the following theorem.  
\begin{theorem}
\label{prop8} Let $\cC\subseteq \bT$. Then, the $\dgh$-closure of $\cC$ is compact iff there exists a sequence $(h_p)_{p\in \bN}$ decreasing to $0$ and a countable dense subset $\cD\subset (0, \infty)$ such that 
\begin{equation}
\label{tightcriter}
\forall p  \in \bN  \, , \; \forall a \in \cD \, , \quad  \quad  \sup_{\widetilde{T}\in\cC} \;  Z_a^{(h_p)}(\widetilde{T})\; <\infty \; \, .
\end{equation}
The same statement holds true when $Z^{(h_p)}_a $ is replaced by $Z^{(h_p)}_{a+}$.  
 \end{theorem}  
The following tightness criterion for random trees is a consequence of Theorem \ref{prop8}. 
\begin{theorem}
\label{tightness} 
Let $(\widetilde{\cT}_j)_{j \in J}$ be family of $\bT$-valued random trees. Their laws are 
$\dgh$-tight 
if and only if for every fixed $h, a \in (0, \infty)$, the laws of the $\bN$-valued random variables 
$(Z_a^{(h)}(\widetilde{\cT}_j))_{j\in J}$ are tight. 
The same result holds true when $(Z_a^{(h)}(\widetilde{\cT}_j))_{j\in J}$ is replaced by $(Z_{a+}^{(h)}(\widetilde{\cT}_j))_{j\in J}$. 
\end{theorem}
\begin{pf} Suppose that for every fixed $h, a \in (0, \infty)$ the laws of $(Z_a^{(h)}(\widetilde{\cT}_j))_{j\in J}$ are tight. Let $\cD= \{ a_q ; q\in \bN\}$ be some dense subset of $(0, \infty)$ and let $(h_p)_{p\in \bN}$ be a sequence decreasing to $0$. Let $\varepsilon \in (0, 1)$. For all $p, q \in \bN$, there exists 
$K_{p,q}\in (0, \infty)$ such that 
$$ \sup_{j \in J} \bP ( Z^{(h_p)}_{a_q} (\widetilde{\cT}_j)  >K_{p,q} ) < \varepsilon 2^{-p-q-2} \; .$$
We then set $\cC= \{ \widetilde{T} \in \bT\colon  \forall p,q\in \bN, \, Z^{(h_p)}_{a_q} (\widetilde{T}) \leq K_{p,q}  \}$. Theorem \ref{prop8} implies that its $\dgh$-closure $\overline{\cC}$ is compact and we 
easily prove that for every $j\in J$ 
$$  \bP (  \widetilde{\cT}_j \notin \overline{\cC} ) \leq \bP ( \exists   p, q \in \bN \colon  Z^{(h_p)}_{a_q} (\widetilde{\cT}_j)  >K_{p,q} ) \leq \sum_{p,q \in \bN}  \bP ( Z^{(h_p)}_{a_q} (\widetilde{\cT}_j)  >K_{p,q} ) <\varepsilon \; , $$
which entails the tightness of the laws of $\widetilde{\cT}_j$, $j\in J$. Conversely, let $\varepsilon \in (0, 1)$. There exists a $\dgh$-compact subset $\cK \subset \bT$ such that  $\inf_{j\in J} \bP (\widetilde{\cT}_j \in \cK) \geq 1-\varepsilon $. Fix $a, h \in (0, \infty)$, $r > a+h$, and set $K= \sup_{\widetilde{T} \in \cK} N_{\widetilde{T}} (h,r) $. By (\ref{Grogene}), 
$ K< \infty$ and Lemma \ref{keycompact} implies that $\sup_{\widetilde{T} \in \cK} Z_a^{(h)} ( \widetilde{T}) \leq K$. Thus, $\inf_{j\in J} \bP (Z_a^{(h)} (\widetilde{\cT}_j ) \leq K ) \geq 1-\varepsilon $, which implies the tightness of the laws of $Z_a^{(h)}(\widetilde{\cT}_j)$, $j\in J$, for each $a,h\in (0, \infty)$.  
\end{pf}
\begin{corollary}
\label{convcriter} Let $(\widetilde{\cT}_n)_{n\in\bN}$ be a sequence of $\bT$-valued random variables. We suppose that for every sufficiently small $h \in (0, \infty)$, the laws of $R_h (\widetilde{\cT}_n)$, 
$n \in \bN$, converge weakly in $(\bT, \dgh)$ as $n \to \infty$. Then, the laws of the 
$\widetilde{\cT}_n $, $n \in \bN$, converge 
weakly in $(\bT, \dgh)$ as $n \to \infty$. 
\end{corollary}
\begin{pf} Fix $a, h\in (0, \infty)$. Since the laws of $R_{h/2} (\widetilde{\cT}_n)$, $n \in \bN$,   converge weakly, Theorem \ref{tightness} implies that the laws of 
$Z_a^{(h/2)} (R_{h/2} (\widetilde{\cT}_n))= Z^{(h)}_a (\widetilde{\cT}_n)$, $n \in \bN$, are tight. Theorem \ref{tightness} then entails that the laws of $\widetilde{\cT}_n$, $n \in \bN$, are 
$\dgh$-tight. 
Denote by $\Lambda_1$ and $\Lambda_2$ two limit laws of $(\widetilde{\cT}_n)_{n\in\bN}$. Then, the laws of $R_h$ under $\Lambda_1$ and $\Lambda_2$ coincide, since $R_h$ is $\dgh$-continuous and since we assume that the laws of $R_h (\widetilde{\cT}_n)$, $n\in \bN$, converge weakly. 
Thus, $\Lambda_1=\Lambda_2$ since $R_h$ converges uniformly to the identity on $\bT$, as $h$ decreases to $0$.
\end{pf}

  We now discuss a stronger convergence for sequences of CLCR real trees that grow in a the following weak sense. 
\begin{definition}  
\label{weakgrowthdef} Let $\widetilde{T}, \widetilde{T}^\prime \in \bT$. We say that 
$\widetilde{T}$ can be embedded in $\widetilde{T}^\prime$, which is denoted by $\widetilde{T}\preceq \widetilde{T}^\prime$ if we can find representatives $(T, d, \rho)$ and $(T^\prime, d^\prime, \rho^\prime)$ of $\widetilde{T}$ and $\widetilde{T}^\prime$, and an isometrical embedding $f\colon  T\hookrightarrow T^\prime$ such that $f(\rho)= \rho^\prime$. \cq 
\end{definition}
We easily check that $\preceq$ is a partial order on $\bT$ (anti-symmetry property is the only non-trivial point to check). Note that if 
$\widetilde{T}\preceq \widetilde{T}^\prime$, then for all $a\in [0, \infty)$ and for all $h \in (0, \infty)$,   
\begin{equation}
\label{herasegrow}
 Z^+_a (\widetilde{T} ) \leq Z^+_a (\widetilde{T}^\prime ) \, , \quad 
 R_h (\widetilde{T} ) \preceq R_h (\widetilde{T}^\prime) \quad {\rm and} \quad  Z^{(h)}_a (\widetilde{T}) \leq Z^{(h)}_a (\widetilde{T}^\prime) \; .
\end{equation}
We first prove the following convergence criterion for a growing sequence of trees. 
\begin{theorem}
\label{critgrowdeter}  Let $(\widetilde{T}_n)_{n \in \bN}$ be a $\bT$-valued sequence such that $\widetilde{T}_n \preceq \widetilde{T}_{n+1}$, $n \in \bN$. 
Let $\cD \subseteq (0, \infty)$ be dense and let $(h_p)_{p\in  \bN}$ be a sequence decreasing to $0$. 
We assume the following:
\begin{equation}
\label{herasecontrol}
\forall p \in \bN \, , \; \forall a \in \cD \, , \qquad \sup_{n \geq 0} Z^{(h_p)}_a (\widetilde{T}_n) \;  < \;  \infty \; .
\end{equation}
Then, there exists $\widetilde{T}\in \bT$ such that 
$\lim_{n \rightarrow \infty}  \dgh (\widetilde{T}_n , \widetilde{T})= 0$. 
\end{theorem}
The proof of the theorem is in several steps. We first state the following representation lemma. 
\begin{lemma}
\label{represgrow} Let $(\widetilde{T}_n)_{n \in \bN}$ be a $\bT$-valued sequence such that $\widetilde{T}_n \preceq \widetilde{T}_{n+1}$, $n \in \bN$. Then, there exists a pointed Polish space $(E, d, \rho)$ and a sequence 
$T_n \subseteq E$ such that for all $n \in \bN$ 
$$ \rho \in T_n \, , \quad  T_n \subseteq T_{n+1} \quad {\rm and} \quad \textrm{$(T_n , d, \rho)$ is a representative of $\,\widetilde{T}_n$} \; .$$
\end{lemma}
\begin{pf} We first prove recursively that for all $n\in \bN$,  
there is  a representative $(T^*_n , d_n , \rho_n)$ of $\widetilde{T}_n$ and 
an isometrical embedding $f_n \colon  T_n^* \hookrightarrow T_{n+1}^*$ such that $f_n (\rho_n)= \rho_{n+1}$. Indeed, assume that $(T_0^*, d_0, \rho_0), \ldots, (T_n^*, d_n, \rho_n)$ and $f_0, \ldots, f_{n-1}$ exist and satisfy the above mentioned conditions. 
Since $\widetilde{T}_n \preceq \widetilde{T}_{n+1}$, there are
$(T^\prime_n , d^\prime_n , \rho_n^\prime)$ and $(T^*_{n+1} , d_{n+1} , \rho_{n+1} )$ 
that are representatives of $\widetilde{T}_n$ and $\widetilde{T}_{n+1}$, respectively, and there is  
an isometrical embedding $f^\prime_n \colon  T^\prime_n \hookrightarrow T^*_{n+1}$ 
such that $f^\prime_n (\rho_n^\prime)= \rho_{n+1}$. There also exists an 
isometry $\phi\colon T^*_n \rightarrow T^\prime_n$ such that $\phi (\rho_n)= \rho_n^\prime$. 
Then, we set $f_n = f^\prime_n \circ \phi$ that satisfies the desired property, which completes the recurrence.  
  
We next construct the desired spaces in a projective way. To that end, we first set 
$$ \forall n \in \bN \, , \quad S_n := \left\{ \fsigma = (\sigma_k)_{k\in \bN} \; \in \prod_{{k\in \bN}} T^*_{n+k}\colon  
 \sigma_{k+1} = f_{n+k} (\sigma_k) \right\} \quad {\rm and} \quad S_\infty = \bigsqcup_{{n \in \bN}} S_n \; , $$  
where $\sqcup$ stands for the disjoint union. For all $\sigma\in T_n^*$, we also 
define the sequence $\jmath_n (\sigma)= (\sigma_k)_{k \in \bN}$ such that $\sigma_0 = \sigma $ and 
$\sigma_{k+1}= f_{n+k} (\sigma_k)$, $k\in \bN$. 
Note that $\jmath_n $ is one-to-one from $T_n^*$ onto $S_n$. Let $n, n^\prime \in \bN$. 
Observe that for all $\fsigma = (\s_k)_{k\in \bN} \in S_n $ and all  
$\fsigma^\prime  = (\s^\prime_k)_{k\in \bN} \in S_{n^\prime} $,
$$ \forall k \geq \lvert n^\prime -n \rvert \, , \quad d_{n+ k} (\sigma_k , \sigma^\prime_{n-n^\prime+k})= d_{n^\prime +k} (\s_{n^\prime-n+k} , \s^\prime_k) =: d(\fsigma, \fsigma^\prime)\; .$$    
This induces a pseudo-metric $d$ on $S_\infty$ and we see that $d(\fsigma, \fsigma^\prime)= 0$ iff a certain shift of the sequence $\fsigma$ is equal to a certain shift of the sequence 
$\fsigma^\prime$. We now introduce the usual equivalence relation $\equiv$ by specifying that $\fsigma \equiv \fsigma^\prime$ iff $d(\fsigma, \fsigma^\prime)= 0$. 
We denote by $E_\infty = S_\infty / \!\! \equiv$ 
the quotient space and we denote by 
$\pi \colon  S_\infty \rightarrow E_\infty$ the canonical projection that associates with a point $\fsigma$ its 
equivalence class. Then, $d$ induces a true metric on $E_\infty$ that is denoted in the same way (to simplify notation). We next set: 
$$ \rho= \pi ( \jmath_0 (\rho_0)) \, \quad {\rm and} \quad  T_n = \pi (S_n) \, , \quad n \in \bN \; .$$   
First note that $\jmath_n (\rho_n)\equiv \jmath_0 (\rho_0)$. Thus, $\rho \in T_n$. 
We next check that $T_n \subseteq T_{n+1}$. Indeed, let $\sigma \in T^*_n$ and set $\jmath_n (\sigma)= (\sigma_k )_{k\in \bN}$. Then, $\sigma_1 \in T_{n+1}^*$ and $\jmath_{n+1} (\sigma_1)= (\sigma_{k+1})_{k\in \bN} \in S_{n+1}$; thus 
$\jmath_n (\sigma) \equiv \jmath_{n+1} (\sigma_1)$ and $\pi (\jmath_n (\sigma) )\in T_{n+1}$, which 
entails $T_n \subseteq T_{n+1}$. 

For all $n \in \bN$, we then set $ \phi_n = \pi \circ \jmath_n \colon  T^*_n 
\rightarrow E_\infty $. It is easy to check that $\phi_n (\rho_n)= \rho$ and that $\phi_n$ is an isometry from $T^*_n $ onto $T_n$. Therefore $(T_n, d, \rho)$ is a CLCR real tree whose pointed isometry class is $\widetilde{T}_n$. Observe that $E_\infty = \bigcup_{n \in \bN} T_n$, which proves that 
$(E_\infty, d_\infty)$ is separable. We get the desired result by taking $E$ as a completion of $(E_\infty , d)$ such that $E_\infty \subseteq E$.\end{pf}

 Recall from Section \ref{realtrees} the definition of $\bT_E$. We next prove the following criterion for convergence in $\bT_E$. 
\begin{lemma}
\label{deterconv} Let $(E, d, \rho)$ be a pointed Polish space. 
Let $(T_n)_{n\in \bN}$ be a $\bT_E$-valued sequence such that $T_n \subseteq T_{n + 1}$. 
Let $\cD \subseteq (0, \infty)$ be dense and let $(h_p)_{p\in  \bN}$ be a sequence decreasing to $0$. 
We assume the following. 
\begin{equation}
\label{herasecontrolE}	
 \forall p \in \bN \, , \; \forall a \in \cD \, , \quad \sup_{n \geq 0} Z^{(h_p)}_a (T_n) < \infty \; .
 \end{equation}
Let $T$ be the closure of $\, \bigcup_{n \in \bN} T_n $ in $(E, d)$. Then,  
$$ T \in \bT_E \quad {\rm and} \quad  \lim_{n\rightarrow \infty} {\bf d}_E (T_n , T)= 0 \; .$$ 
\end{lemma}
\begin{pf} We first set $T_\infty= \bigcup_{n\in \bN} T_n$, so that $T$ is the closure of 
$T_\infty$ is $E$. The restriction of $d$ to $T_\infty$ satisfies the four points condition. 
We check that the same holds true for the restriction of $d$ to $T$. Moreover, every point of 
$T_\infty$ belongs to a certain $T_n$ and 
is connected to $\rho$ by a geodesic. This easily implies that 
$(T, d)$ is a connected space. Thus, $(T, d, \rho)$ is a Polish rooted real tree. 

  We fix $h, a \in (0, \infty)$ and we  set 
$L_a = \{ \sigma \in R_{2h} (T)\colon  d(\rho, \sigma)= a\}$. Thus, $\# L_a = Z^{(2h)}_a (T)$. Let  $\sigma \in L_a$. 
There exists  $s \in T $ such that $ \sigma \in \lgeo \rho , s\rgeo$ and $d(\sigma, s) \geq 2h$. Since $T$ is the closure of $T_\infty$, there exists $n (\sigma ) \in \bN$ and $\gamma \in T_{n (\sigma )}$ such that 
$d (s, \gamma ) < h$. This implies $\sigma \in \lgeo \rho , \gamma \rgeo$ and $d( \sigma, \gamma ) > h$. Namely, this implies that $\sigma \in R_h (T_{n (\sigma )})$ and thus $\sigma \in R_h (T_{n})$, for all $n \geq n (\sigma )$. We then get 
$$ Z^{(2h)}_a (T)= \sup_{n \geq 0} \;  \# \{ \sigma \in L_a   \colon    n(\sigma ) \leq n \}  \leq \sup_{ n \geq 0} Z^{(h)}_a (T_n)  \; . $$
Note that this inequality holds true for all $h, a \in (0, \infty)$. By (\ref{herasecontrolE}), 
Lemma \ref{keycompact} applies and $(T, d, \rho)$ is a CLCR real tree, which shows that $T\in \bT_E$. 

  Let us fix $r, \varepsilon \in (0, \infty)$. Since $\overline{B}_T (\rho, r)$ is compact and since $T_\infty$ is dense in $T$, we easily see that 
there exist $n_\varepsilon \in \bN$ and $\sigma_1, \ldots, \sigma_p \in \overline{B}_{T_{n_\varepsilon} }(\rho , r)$ such that for all $\sigma\in \overline{B}_T (\rho , r)$, $d( \{ \sigma _1, \ldots , \sigma_p\} , \sigma ) < \varepsilon $. For all $n \geq n_\varepsilon$, since $ \overline{B}_{T_{n_\varepsilon}} (\rho , r) \subseteq \overline{B}_{T_n} (\rho, r) \subseteq   \overline{B}_T (\rho, r) $, we get $d_{{\rm Haus}} (\overline{B}_{T_n} (\rho, r), \overline{B}_T (\rho, r)) < \varepsilon$, which easily completes the proof of the lemma. 
\end{pf}

\noi
{\bf Proof of Theorem \ref{critgrowdeter}.} Let $(\widetilde{T}_n)_{n \in \bN}$ be a $\bT $-valued sequence that satisfies $\widetilde{T}_n 
\preceq \widetilde{T}_{n+1}$ and the assumption (\ref{herasecontrol}). 
By Lemma \ref{represgrow}, there exist 
a pointed Polish space $(E, d, \rho)$ and $T_n \in \bT_E$, $n \in \bN$, that satisfy $T_n \subseteq T_{n+1}$, (\ref{herasecontrolE}) and such that $(T_n , d, \rho)$ is a representative of $\widetilde{T}_n$. Then,  Lemma \ref{deterconv} applies and there exists $T\in \bT_E$ such that $\lim_{n\rightarrow \infty} {\bf d}_E (T_n , T)= 0 $, which implies $\lim_{n\rightarrow \infty}  \dgh (\widetilde{T}_n , \widetilde{T})= 0$, by (\ref{dcontrol}). \cqfd 

\bigskip

\noi
Theorem \ref{critgrowdeter} entails the following criterion for growing sequences of random trees. 
\begin{theorem}
\label{tightness2} 
Let $( \widetilde{\cT}_n)_{n \in \bN}$ be a $\bT$-valued sequence of random trees such that for all $n \in \bN$, 
$$ \textrm{$\bP$-a.s.} \quad \widetilde{\cT}_n \preceq \widetilde{\cT}_{n+1} \; .$$

\noindent We assume that for every fixed $a, h\in (0, \infty)$, the family of laws of $\bN$-valued random variables 
$Z_a^{(h)}(\widetilde{\cT}_n)$, $n\in \bN$, is tight. Then, there exists a random tree $\widetilde{\cT}$ in $\bT$ such that 
$$\textrm{$\bP$-a.s.} \quad \lim_{n \rightarrow \infty} \dgh (\widetilde{\cT}_n , \widetilde{\cT})= 0 \; .$$
\end{theorem}
\begin{pf}We fix $a, h \in (0, \infty)$. By (\ref{herasegrow}), a.s.$\,Z_a^{(h)}(\widetilde{\cT}_n) \leq Z_a^{(h)}(\widetilde{\cT}_{n+1}) $, for all $n \in \bN$. Then for all $K \in (0, \infty)$, 
$$\bP \left( \sup_{{n \in \bN}} Z_a^{(h)}(\widetilde{\cT}_n) > K \right) = \sup_{{n \in \bN}} \bP \left(Z_a^{(h)}(\widetilde{\cT}_n) >K \right) \; .$$
Since the family of laws of $Z_a^{(h)}(\widetilde{\cT}_n)$, $n\!\in\!\bN$, is tight, we get 
$\bP (\sup_{n \in \bN} Z_a^{(h)}(\widetilde{\cT}_n)\!<\!\infty)\!=\!1$, for each $a, h \in (0, \infty)$, and Theorem \ref{critgrowdeter} easily completes the proof. 
\end{pf}

Lemma \ref{deterconv} also entails the following almost sure convergence criterion for random trees in $\bT_E$. The proof, quite similar to that of Theorem \ref{tightness2}, is left to the reader.  
\begin{theorem}
\label{tightness3} 
Let $(E, d , \rho)$ be a pointed Polish space. 
Let $( \cT_n)_{n \in \bN}$ be a $\bT_E$-valued sequence of random trees such that a.s.\ $\cT_n \subseteq \cT_{n+1}$, $n \in \bN$. We assume that for each fixed $a, h\in (0, \infty)$, the family of laws of $\bN$-valued random variables 
$Z_a^{(h)}(\cT_n)$, $n\in \bN$, is tight. Let $\cT$ be  the closure of $\bigcup_{n \in \bN} \cT_n $ in $E$. Then, 
$$  \textrm{$\bP$-a.s.} \quad \cT \in \bT_E \quad {\rm and} \quad \lim_{n\rightarrow \infty} {\bf d}_E (\cT_n , \cT)= 0 \; .$$
 \end{theorem}

\section{Galton-Watson trees and reduction by hereditary properties}
\label{GWsec}

\subsection{Galton-Watson trees as random real trees.}
\label{GWrealtreesec}
In this section, we give an intrinsic definition of Galton-Watson trees as $\bT$-valued random variables. Informally, given $\xi= (\xi (k))_{k\in \bN}$, a probability measure on $\bN$, and $c \in (0, \infty)$, a Galton-Watson tree with offspring distribution $\xi$ and lifetime parameter $c$ is the genealogical tree of a population that has a single progenitor and evolves as follows: the lifetimes of the individuals are independent and exponentially distributed with mean $1/c\, $; when they die, individuals independently give birth to a 
random number of children distributed according to $\xi$. 
We denote by $\varphi_\xi$ the generating function of $\xi$. 
Namely, 
\begin{equation}
\label{varphiandmeannota}
\forall\,  r \in [0, 1]\colon \sum_{k \geq 0} r^k \xi (k) = \varphi_\xi  (r) \; .
\end{equation}
We make the following assumptions on the offspring distribution $\xi$: 
\begin{equation}
\label{nonexpdiscr}
\textsl{non-trivial: $\; \xi (0)+ \xi (1) < 1 , \quad$ } \; \textsl{proper: $\; \xi (1)= 0 , \quad $}  \; \textsl{conservative:} \; \int^{1-}  \frac{dr}{ ( \, \varphi_\xi (r) -r )_- } \, =\, \infty \, , 
\end{equation}
where $(\, \cdot \, )_-$ stands for the negative part. We assume that $\xi$ is 
proper because we are only interested in the underlying geometrical tree. 
The conservativity assumption comes from a standard criterion that asserts that the right profile of a Galton-Watson tree is an $\bN$-Markov process that is conservative iff the last condition of (\ref{nonexpdiscr}) is satisfied 
(see \cite[Section III.3]{AthNey} for more details).

Recall from Section \ref{funcsubsec} that for all $\widetilde{T}_1, \widetilde{T}_2 \in \bT$, $\widetilde{T}_1 \circledast \widetilde{T}_2$ stands for the pointed isometry class of the tree obtained by pasting two representatives of $\widetilde{T}_1$ and $\widetilde{T}_2$ at their roots. We let the reader check that $(\widetilde{T}_1 ,  \widetilde{T}_2) \mapsto \widetilde{T}_1 \circledast \widetilde{T}_2$ is a continuous, symmetric and associative operation on $\bT$. For all Borel probability measures $Q_1$ and $Q_2$ on $\bT$, we define $Q_1 \circledast Q_2$ as the image measure of the product measure $Q_1 \otimes Q_2$ under $\circledast$. Namely, for all measurable $G\colon \bT \rightarrow [0, \infty)$, 
$$ \int_{\bT}  (Q_1 \circledast Q_2 )(d \widetilde{T}) \,  G ( \widetilde{T}) = \int_\bT\int_\bT Q_1 ( d\widetilde{T}_1)Q_2 (d \widetilde{T}_2) \,  G( \widetilde{T}_1 \circledast \widetilde{T}_2) \; .$$
We easily check that $(Q_1 , Q_2)\mapsto Q_1 \circledast Q_2$ is a weakly continuous, symmetric and associative operation. For all Borel probability measures $Q$ on $\bT$, and all $n \in \bN$, we recursively define $Q^{\circledast n}$ by setting $Q^{\circledast(n+1)}= Q\circledast Q^{\circledast n}$ where $Q^{\circledast 0} $ is the Dirac mass at the point tree $\pnt$. 
Recall from Section \ref{funcsubsec} the definition of $D$, $\vartheta$ and ${\bf k}$. 
\begin{definition} 
\label{GWrealtreedef} 
Let $\xi$ be a proper and conservative offspring distribution. Let $c \in (0, \infty)$. 
\begin{enumerate}
\item[(a)] {\it A Galton-Watson real tree with offspring distribution $\xi$ and edge parameter $c$} 
(a GW($\xi, c$)-{\it real tree} for short) is a  $\bT$-valued random variable whose distribution $Q$ satisfies the following. For all $n \in \bN$, for all measurable functions $g\colon [0, \infty) \rightarrow [0, \infty)$ and $G\colon \bT \rightarrow [0, \infty)$, 
\begin{equation}
\label{defiGWcalcul}
Q\left[  \un_{\{{\bf k} = n \} } G (\vartheta) g( D)   \right]= \xi (n) Q^{\circledast n} [ \, G \, ]  \int_0^\infty g(t) ce^{-ct} dt \; .
\end{equation}
\item[(b)] {\it A Galton-Watson real forest with offspring distribution $\xi$, edge parameter $c$ and initial distribution $\mu$} (a GW($\xi, c; \mu$)-{\it real forest} for short) is a $\bT$-valued random variable whose distribution is $ P= \sum_{n \geq 0} \mu(n) Q^{\circledast n} $, where $Q$ stands for the distribution of a GW($\xi, c)$-real tree.  \cq
\end{enumerate}
\end{definition}
\begin{lemma}\label{uniqGW} 
For each proper conservative offspring distribution $\xi$ and each $c\in(0,\infty)$, there exists a unique distribution $Q_{\xi, c}$ on $\bT$ that satisfies 
{\rm(\ref{defiGWcalcul})} in Definition \ref{GWrealtreedef} {\rm(a)}. Moreover, $Q_{\xi, c}(\bT_{\rm edge})= 1$ and 
$(\xi,  c )\mapsto  Q_{\xi, c }$ is injective. 
\end{lemma}
\begin{pf} See Appendix \ref{uniqGWpf}.
\end{pf}

\noi
{\bf Notation.} For all proper and conservative distributions $\xi$, for all $c \in (0, \infty)$ and for all probability measures $\mu$ on $\bN$, we set $P_{\xi, c , \mu} =  \sum_{n \geq 0} \mu(n) Q_{\xi, c}^{\circledast n} $ as the law of a GW($\xi, c ;\mu$)-real forest. \cq 

\medskip 

In the non-conservative case $Q_{\xi , c}$ does not exist as a distribution on $\bT$. In fact, explosive GW-real trees would not be locally compact. Intuitively, the trees would be well-behaved up to and including the first explosion height, but there would be infinitely many subtrees above this height destroying local compactness. We do not intend to develop this more formally here.

 GW-real trees are characterized by their {\it regenerative branching property} as proved by 
Weill \cite{Weill} in the compact case. 
More precisely, for every $a \in [0, \infty)$, we denote by $\cB_a $ the sigma-field on $\bT$ generated by $\blw (a, \cdot )$. Note that $(\cB_a )_{ a \in [0, \infty)}$ is a filtration on $\bT$. We denote by $(\cB_{a+})_{a \in [0, \infty)}$ the associated right-continuous filtration. Recall from (\ref{asplitdef}) the definition of the right profile at height $a$ that is denoted by $Z^+_a$. We easily check that $Z^+_{a}$ is $\cB_{a+}$-measurable. 
\begin{lemma}
\label{regenene} Let $Q$ be a Borel probability measure on $\bT$. We first assume that 
$Q(0 < D< \infty)>0 $. We also assume that $Q$-a.s.\ the process $a\mapsto Z^+_a $ is $\bN$-valued and cadlag.  
Then, the following assertions are equivalent.
\begin{enumerate}
\item[$\rm(i)$] For every $a \in [0, \infty)$, the conditional distribution given $Z^+_a$ of $\abv (a, \cdot)$ under $Q$ is $Q^{\circledast Z^+_{a} } $.   
\item[$\rm(ii)$] For every $a \in [0, \infty)$, the conditional distribution given $\cB_{a+}$ of $\abv (a, \cdot)$ under $Q$ is $Q^{\circledast Z^+_{a} } $. 
\item[$\rm(iii)$] $Q$ is the distribution of a {\rm GW($\xi, c$)}-real tree, for a certain $c \in (0, \infty)$ and a certain proper conservative offspring distribution $\xi$. 
\end{enumerate}
\end{lemma}
\begin{pf}\rm  In the compact case, this statement is close to \cite[Theorem 1.2]{Weill}. We briefly prove Lemma \ref{regenene} in Appendix \ref{regenenepf} using different arguments. We also mention that the cadlag assumption can be dropped, as we show in Theorem \ref{vrairege}.
\end{pf}
Let us recall basic results on the branching processes associated with GW-trees. Let $c\in (0, \infty)$ and let $\xi$ be a non-trivial proper conservative offspring distribution. The regenerative branching property entails that $(Z^+_{a})_{ a \in [0, \infty)}$ under $Q_{\xi,  c}$ is an $\bN$-valued Markov process whose 
matrix-generator $(q_{i,j})_{ i,j \in \bN}$ is given by $ q_{i,j} = c\, i \, \xi (j-i+1)$, if $j>i$ or $j=i-1$, $q_{i,i}= -c\, i$ and $q_{i,j}=0$ if $j<i-1$. Let us set 
\begin{equation}
\label{laplvdef}
w(a, \theta):=Q_{\xi, c} \left[ e^{-\theta Z^+_{a}  } \right]  \;   , \quad  a, \theta \in [0, \infty) . 
\end{equation}
The definition of GW-real trees implies that for all $\theta \in [0, \infty)$, $a \mapsto w(a, \theta)$ is the unique nonnegative solution of the differential equation $\partial w (a, \theta)  /\partial a= c \left( \varphi_\xi (w(a, \theta)) -w(a, \theta) \right)$, with $w(0, \theta) = e^{-\theta} $. 
A simple change of variables allows to rewrite the previous equation in the following form: 
\begin{equation}
\label{inteqharris}
\forall \theta  \in (0, \infty) \backslash  \{-\log q_\xi\}  \; ,  \; \forall a \in [ 0, \infty) \; ,\quad 
\int_{e^{-\theta}}^{w(a,\theta)} \frac{dr}{\varphi_\xi (r) -r} \, =\,
c   a  \; , 
\end{equation}
where $q_\xi$ is the smallest solution of the equation $\varphi_\xi (r)= r$ and where we agree on $-\log (q_\xi)= \infty$ if $q_\xi= 0$ (see e.g.$\;$\cite[Section III.3]{AthNey} for more details). 
Note that $Q_{\xi, c}$-a.s.$\; \Gamma= 
\inf\{ a \in [0, \infty)\colon Z^+_a= 0\}$. Thus,
\begin{equation}
\label{heightharris}
\textrm{$w(a):=Q_{\xi, c} (\Gamma \leq a) =\lim_{\theta \rightarrow \infty} \downarrow w(a,\theta )$ satisfies} \quad  \int_{0}^{w(a) }\frac{dr}{\varphi_\xi (r) -r} \, =\,
c a. 
\end{equation}
Let $\mu$ be a probability measure on $\bN$ that is distinct from $\delta_0$. We easily derive from Lemma \ref{regenene} that for all $a\in [0, \infty)$,  
\begin{equation}
\label{branchfor}
\textrm{the conditional distribution given $Z^+_a$ of $\abv (a, \cdot)$ under $P_{\xi, c,\mu}$ is 
$Q_{\xi, c}^{\circledast Z^+_a}$.}
\end{equation}
Then, under $P_{\xi, c, \mu}$, $(Z^+_a)_{a\in [0, \infty)}$ is also a Markovian branching process with offspring distribution $\xi$, with lifetime parameter $c$ and with initial distribution $\mu$. Let $\varphi_\mu$ be the generating function of $\mu$. We thus get 
\begin{equation}
\label{forestharris}
 \forall a, \theta \in [0, \infty) \, , \quad P_{\xi, c, \mu} \left[e^{-\theta Z^+_a} \right]= \varphi_\mu \left( w(a, \theta) \right) \; .
\end{equation}

\begin{lemma}
\label{poubelle} Let $(\xi_n)_{n \geq 0}$, be a sequence of proper conservative offspring distributions that converges weakly to the proper conservative offspring distribution $\xi_{\infty}$. Let $(c_n)_{n \geq 0}$ be a sequence of positive real numbers that converges to $c_\infty \in (0, \infty)$. Let 
$(\mu_n)_{n \geq 0}$ be a sequence of probability measures on $\bN$ that converges weakly to the probability measure $\mu_\infty$ on $\bN$. Then, the following convergences hold weakly on $\bT${\rm:}   
$$  \lim_{n\rightarrow \infty} Q_{\xi_n , c_n}= Q_{\xi_\infty , c_\infty} \quad {\rm and}  \quad  \lim_{n\rightarrow \infty} P_{\xi_n , c_n, \mu_n}= P_{\xi_\infty , c_\infty, \mu_\infty} \; .$$
Moreover, for all $a\in [0, \infty )$, as $n \rightarrow\infty$, 
the law of $Z^+_a$ under $Q_{\xi_n , c_n }$ (resp.\ under 
$P_{\xi_n , c_n, \mu_n}$) converges weakly to the law of $Z^+_a$ under $Q_{\xi_\infty , c_\infty}$ 
(resp.\ under $P_{\xi_\infty , c_\infty , \mu_\infty}$).   
\end{lemma}
\begin{pf} See Appendix \ref{poubellepf}. 
\end{pf}

\subsection{Hereditary properties and their reduction procedures.} 
\label{realhereditsec}
Here, we generalise $h$-erasure $R_h(T)$ of a CLCR real tree $(T,d,\rho)$ to reduction for more general hereditary properties as explained in the introduction. 
More precisely, let $(T, d, \rho)$ be a CLCR real tree and let $\sigma \in T$. Recall that  {\it the subtree above $\sigma$} is given by 
\begin{equation} 
\label{abvsubtree}
\theta_\sigma T = \left\{ \sigma^\prime \in T\colon \sigma \in \lgeo \rho , \sigma^\prime \rgeo  \right\}\; .
\end{equation}
Note that $(\theta_\sigma T , d , \sigma) $ is a CLCR real tree. We simply denote by $\widetilde{\theta}_\s T$ its pointed isometry class. Observe that 
\begin{equation}
\label{shiftconst}
\forall \s^\prime \in \lgeo \rho , \sigma \rgeo \, , \quad \theta_\s T= \theta_\s (\theta_{\s^\prime} T)\; .  
\end{equation}

\begin{definition}
\label{heredrealdef} 
\begin{enumerate}
\item[(a)] Let $A \subset \bT$ be a Borel subset of $\bT$ such that 
$\pnt \notin A$. 
The set $A$ is {\it hereditary} if it satisfies the following: for every CLCR real tree 
$(T,d,\rho)$ and for every $\s \in T$, if $\widetilde{\theta}_\s T \in A$, then $\widetilde{T} \in A$. 
\item[(b)] Let $(T,d, \rho)$ be a CLCR real tree and let $A \subset \bT$ be hereditary. 
We denote by $R_A (T)$ the closure in $T$ of the subset  
\begin{equation}
\label{RAset}
 \{ \rho \} \cup \left\{ \sigma \in T\colon \widetilde{\theta}_\sigma T \in A    \right\} \; .
\end{equation} 
Then, $(R_A (T), d, \rho)$ is a CLCR real tree that we call the {\it $A$-reduced tree of $T$}. \cq 
\end{enumerate}
\end{definition}
Let us briefly explain why $(R_A (T), d, \rho)$ is a CLCR real tree. 
Since $R_A(T)$ is a closed subset of $T$, we only need to prove that it is connected. 
Denote by $S$ the set given by (\ref{RAset}). Let $\s \in R_A(T) \backslash \{ \rho \}$ (if any).
Then, there exists a sequence $\sigma_n \in S$ converging to $\s$ and such that $\widetilde{\theta}_{\s_n} T \in A$. Let $\s^\prime \in \lgeo \rho , \s \lgeo \, $. 
For all sufficiently large $n$, $\s^\prime \in \lgeo \rho , \s_n \lgeo \, $, and 
(\ref{shiftconst}) implies that $\widetilde{\theta}_{\s^\prime} T \in A$. This proves that $\lgeo \rho, \s \rgeo \subseteq R_A (T)$, for all $\s \in R_A (T)$. 

\medskip

Clearly, the pointed isometry class of $(R_A (T), d, \rho)$ only depends on $\widetilde{T}$ and $A$. This then induces a function on $\bT$ 
that is denoted in the same way by $R_A\colon  \bT \rightarrow \bT$. Note that $R_A (\pnt)= \pnt$ and that 
\begin{equation}
\label{Aredembed}
\forall \widetilde{T} \in \bT \, , \quad R_A (\widetilde{T} ) \preceq \widetilde{T} \; .
\end{equation}
\begin{lemma}
\label{measAred} Let $A \subseteq \bT$ be a hereditary property. Then $R_A\colon  \bT \rightarrow \bT$ is Borel-measurable.  
\end{lemma}
\begin{pf}
See Appendix \ref{RAproofsec}. 
\end{pf}
\begin{lemma}
\label{commulem}
Let $(T,d, \rho)$ be a CLCR real tree and let $A\subseteq \bT$ be hereditary. Then, the following holds true: 
\begin{equation}
\label{commushift}
\forall \s \in R_A (T) \, , \quad \theta_\s R_A (T) = R_A (\theta_\s T) \; .
\end{equation}
\end{lemma}
\begin{pf} Set $S= \{ \s^\prime \in \theta_\s T \colon  \widetilde{\theta}_{\s^\prime} T \in A  \}$. Then, $R_A (\theta_\s T)$ is the closure of $\{ \s \} \cup S$. Note that if $\s \in R_A (T)$, then $R_A (T) \cap \theta_\s T =   \theta_\s R_A (T) $. Thus, $\{ \sigma \} \cup S \subseteq \theta_\s R_A (T) $, 
which implies that $R_A (\theta_\s T)\subseteq  \theta_\s R_A (T)$. 
Conversely, let $\s^\prime \in  R_A (T) \cap \theta_\s T $ be distinct from $\s$ (if any). 
By definition of $R_A (T)$, there exists a sequence $\s_n $ converging to $\s^\prime$ such that 
$\widetilde{\theta}_{\s_n} T\in A $. Since $\s \in \lgeo \rho , \s^\prime \lgeo\, $, 
then for all sufficiently large $n$ we get $\s_n \in \theta_\s T$, and thus $\s_n \in S$. This proves that 
 $R_A (T) \cap \theta_\s T$ is contained in the closure of $\{ \s\} \cup S$, which completes the proof of 
(\ref{commushift}).  
\end{pf}
Let us now discuss how the tree reduction 
behaves with respect to the functions $\cM_a$, $\abv (a, \cdot) $ and $Z^+_a$. To that end, for every  hereditary property $A \subseteq \bT$ we introduce   
\begin{equation}
\label{Aminus}
A^- = \left\{ \widetilde{T} \in \bT\colon R_A (\widetilde{T} ) \neq \pnt \right\} \; ,
\end{equation}
that can be viewed as a regular version of $A$.  
\begin{example}
\label{heraexex}
For all $h \in [0, \infty)$, consider the hereditary property $A_h = \{ \Gamma \geq h\}$. 
Note that $R_{A_h} $ is the $h$-leaf length erasure function 
as defined by (\ref{hleafera}). Namely $R_{A_h}= R_h$. Thus, $A_h^-= \{ \Gamma >h \}$.  \cq 
\end{example}
\begin{lemma}
\label{Amoinsprop} Let $A\subseteq \bT$ be hereditary. Then, the following holds true. 
\begin{enumerate}
\item[\rm(i)] $A^- \subseteq A$ and $A^-$ is hereditary. Moreover, $R_{A^-} (T)= R_A(T)$, for every CLCR real tree $(T, d, \rho)$. This implies that $(A^-)^-= A^-$. 
\item[\rm(ii)] Let $a\in [0, \infty)$ and let $\widetilde{T} \in \bT$. If $\cM_a (\widetilde{T})= \sum_{i\in I} \delta_{\widetilde{T}_i}$, then we get 
\begin{equation}
\label{asplitAred}
\cM_a (R_A (\widetilde{T}))= \sum_{i\in I\colon  \widetilde{T}_i\in A^-}
\delta_{R_A (\widetilde{T}_i)} \quad {\rm and} \quad Z^+_a (R_A (\widetilde{T}))=\langle \cM_a (\widetilde{T}) , \un_{A^-} \rangle \; .
\end{equation}
Moreover, 
\begin{equation}
\label{AbovAred}
\!\!\!\!  \underset{\underset{\scriptstyle\; \widetilde{T}_i\in A^-}{i\in I\colon }}{\!\! \circledast} \!\! \! R_A (\widetilde{T}_i) =   \abv (a, R_A (\widetilde{T}))= R_A (\abv (a, \widetilde{T}))  \; .
\end{equation}
\end{enumerate}
\end{lemma}
\begin{pf} We first prove $\rm(i)$. 
By 
Lemma \ref{measAred}, $A^-$ is a Borel subset of $\bT$. We next prove that $A^- \subseteq A$: indeed, let 
$(T, d , \rho)$ be a CLCR real tree such that $\widetilde{T} \in A^-$; then, 
there is $\sigma \in T \backslash \{ \rho \}$, such that $\widetilde{\theta}_\s T \in A $, which implies $\widetilde{T} \in A$. Let us prove that $A^-$ is hereditary.  Assume that $\widetilde{\theta}_\s T \in A^-$, then there exists $\s^\prime \in \theta_\s T 
\backslash \{ \s \} $ such that $\widetilde{\theta}_{\s^\prime} T \in A$, which implies that $R_A (\widetilde{T})\neq \pnt$, and $\widetilde{T} \in A^-$, which entails that $A^-$ is hereditary. 

Let us prove that $R_{A^-} (T)= R_A(T)$. Since 
$A^-\subseteq A$, $R_{A^-} (T) \subseteq  R_A(T)$. 
Denote by $S$ the set defined in (\ref{RAset}) and take $\s \in S \backslash \{ \rho\}$ (if any). Then, for all $\s^\prime \in \lgeo \rho, \s \lgeo \, $, $\s \in R_A (\theta_{\s^\prime} T)$, which implies that 
$\widetilde{\theta}_{\s^\prime} T \in A^-$ and therefore $\s^\prime \in R_{A^-} (T)$. 
Thus, $S$ is in $R_{A^-} (T)$, which implies the desired result.

We next prove $\rm(ii)$. Let $(T, d, \rho)$ be a CLCR real tree. We denote by $T^\circ_i $, $i \in I$, the connected components of the open set $\{ \s \in T \colon  d( \rho ,\s) >a\}$. For every $i\in I$, denote by $\s_i$ the unique point of $T$ such that $d(\rho , \s_i)= a$ and such that $T_i := \{\s_i \} \cup T^\circ_i$ is the closure 
of $T^\circ_i$. Recall that the CLCR real trees $(T_i, d, \s_i )$ are the subtrees above $a$, namely: 
$\cM_a (T)= \sum_{i\in I} \delta_{\widetilde{T}_i}$. Let $i\in I$ be such that $\widetilde{T}_i \in A^-$. 
Thus,  $R_A (T_i)$ does not reduce to $\{ \s_i \}$ and $R_A (T_i) \backslash \{ \s_i\}$ is connected and non-empty. Namely, it is a connected component of $\{ \s \in R_A (T) \colon  d(\rho, \s) >a \}$. Conversely, every connected component of $\{ \s \in R_A (T) \colon  d(\rho, \s) >a \}$ is of this form. This easily completes the proof of $\rm(ii)$. 
\end{pf}

Recall that for leaf-length erasure $R_{h+h^\prime}= R_h \circ R_{h^\prime}$. We extend these relations by introducing the composition of hereditary properties. 
\begin{definition}
\label{comphere} Let $A,A^\prime \subseteq \bT$ be two hereditary properties. The {\it composition of $A$ by $A^\prime$} is defined as the set 
$ A^\prime \circ A := \{  \widetilde{T} \in \bT\colon R_A (\widetilde{T}) \in A^\prime \}$. \cq 
\end{definition}
\begin{lemma}
\label{comphereprop} Let $A, A^\prime \subseteq \bT$ be hereditary. Then,  
$A^\prime \circ A$ is hereditary and for every CLCR real tree $(T, d, \rho)$, we get 
$ R_{A^\prime} (R_{A} (T))= R_{A^\prime \circ A} (T) $. 
\end{lemma}
\begin{pf} Lemma \ref{measAred} implies that $A^\prime \circ A$ is a Borel subset of $\bT$. 
Let $(T, d, \rho)$ be a CLCR 
real tree and let $\s \in T$. By Lemma \ref{commulem}, 
\begin{equation}
\label{equicirc}
 \widetilde{\theta}_\s T \in A^\prime \circ A \quad \Longleftrightarrow \quad 
\s \in R_A (T) \quad {\rm and} \quad  \widetilde{\theta}_\s R_A(T)  \in A^\prime  \; .
\end{equation} 
Note that if $ \widetilde{\theta}_\s R_A(T)  \in A^\prime$, then $R_A (\widetilde{T})\in A^\prime$, namely   $\widetilde{T} \in A^\prime \circ A$. Then, ``$\Rightarrow $'' in (\ref{equicirc}) implies that $A^\prime \circ A$ is hereditary. Moreover, (\ref{equicirc}) immediately 
implies $R_{A^\prime} (R_{A} (T))= R_{A^\prime \circ A} (T)$.  
\end{pf}

The following theorem, which is the main result of this section, shows that the class of GW-laws is stable under reduction by a hereditary property. 
Recall that $Q_{\xi, c}$ stands for the law of a GW($\xi, c$)-real tree, 
that $\varphi_\xi$ is the generating function of $\xi$ and that $q_\xi$ is the smallest root of $\varphi_\xi (r)= r$. 
\begin{theorem}
\label{GWrealreduc} Let $\xi $ be a non-trivial proper conservative offspring distribution and let $c\in (0, \infty)$. Let $A \subseteq \bT$ be hereditary. We assume that 
\begin{equation}
\label{assuptAred}
\alpha:= Q_{\xi, c} (R_A = \pnt)  \; \in (0,1) \; .
\end{equation}
Then, $\xi (0) >0$, $\alpha \in (0, q_\xi]$ and thus $ \varphi^\prime_\xi (\alpha ) \in (0, 1)$. 
Moreover, $R_A$ under the conditioned law 
$Q_{\xi , c} ( \, \cdot \,  | \, R_A \neq \pnt)$ is distributed as a {\rm GW($\xi^{(\alpha)} , c^{(\alpha)}$)}-real tree where 
\begin{equation}
\label{reduclaw}
c^{(\alpha)}  =  c\,(1-  \varphi_\xi^\prime (\alpha) ) \quad {\rm and} \quad \varphi_{\xi^{(\alpha)}} (r)
=  r + \frac{\varphi_\xi(  \alpha + (1- \alpha)r )- \alpha - (1-\alpha)r}{(1-\alpha)(1- \varphi_\xi^\prime(\alpha))}  \, , \quad r\in [0, 1] \; .
\end{equation}
\end{theorem}
\begin{pf} First recall from (\ref{Aminus}) that $\{ R_A \neq \pnt\}= A^- $ and thus $\alpha=Q_{\xi,c}(\bT\setminus A^-)$. Recall from (\ref{distbranch}) and (\ref{firstabove}) the definition of the functions $D$, $\vartheta$ and ${\bf k}$. Since $A^-$ is hereditary, we get 
$\un_{\bT \backslash A^-} \leq \un_{\{ {\bf k} = 0\}} +  \un_{\{ {\bf k} \geq 1 \}} \un_{\{ \langle \cM_D (\cdot) , \un_{A^-}   \rangle = 0 \}} $. We take the expectation under $Q_{\xi , c}$: since
 $\cM_D = \cM_0 \circ \vartheta $,
(\ref{defiGWcalcul}) in the definition of $Q_{\xi , c} $ entails that 
$$ 0 < \alpha \leq \xi (0) + \sum_{n \geq 1} \xi (n) Q_{\xi , c}^{\circledast n} \left(\,  \langle  \cM_0 , \un_{A^-} \rangle \!=\! 0\,  \right)= \sum_{n \geq 0} \xi (n) \alpha^n = \varphi_\xi (\alpha) \; , $$
which entails that $\alpha \in (0, q_\xi]$ and that $\varphi^\prime_\xi (\alpha ) \in (0, 1)$. 

We denote by $Q$ the law of $R_A$ under $Q_{\xi , c} ( \, \cdot \,  | \,  A^- )$ and we 
want to apply Lemma \ref{regenene} to $Q$. To that end, we first note that $a\mapsto Z_a^+$ is $\bN$-valued and cadlag $Q$-a.s.~and we then claim that also 
\begin{equation}
\label{claimdist}
 Q ( 0<D < \infty) = Q_{\xi , c} \left(\left. 0< D \circ R_A < \infty \, \right| \, A^-  \right) > 0 \; .
\end{equation}
Indeed, let $\widetilde{T} \in \bT_{{\rm edge}}$ be such that $Z^+_0 (\widetilde{T})= 1$ and $D(\widetilde{T}) \in (0, \infty)$. If we have $\langle \cM_{D(\widetilde{T})} (\widetilde{T}) , \un_{A^-} \rangle  \geq 2$, then $D (R_A (\widetilde{T}))= D(\widetilde{T})$ and since $A^-$ is hereditary, we also get 
$\widetilde{T}\in A^-$. Consequently, 
\begin{equation}
\label{stepclaimdist}
x:= Q_{\xi , c} \left( \langle \cM_{D} ,  \un_{A^-} \rangle  \geq 2 \right) \leq Q_{\xi, c} \left( A^- \cap \{ D\circ R_A \in (0, \infty) \} \right) = (1-\alpha) Q (0<D<\infty) \; .
\end{equation}
Now by (\ref{defiGWcalcul}) in the definition of $Q_{\xi , c} $, we get 
\begin{eqnarray*}
x &\!\!\!=\!\!\!& \sum_{n \geq 2} \xi (n)  Q_{\xi , c}^{\circledast n} \left( \langle  \cM_0 , \un_{A^-} \rangle \geq  2 \right)
=  \sum_{n \geq 2} \xi (n)  \left( 1- Q_{\xi , c}^{\circledast n} \left( \, \langle  \cM_0 , \un_{A^-} \rangle \!=\! 0\,  \right)-Q_{\xi , c}^{\circledast n} \left( \, \langle  \cM_0 , \un_{A^-} \rangle \!  =\! 1 \, \right) \right) \\
&\!\!\!=\!\!\!& 
 \sum_{n \geq 2} \xi (n)  \left( 1- \alpha^n -n(1-\alpha) \alpha^{n-1} \right)= 1-\varphi_{\xi} (\alpha) -(1-\alpha) \varphi^\prime_{\xi} (\alpha) \; >0 \, ,   
\end{eqnarray*} 
since $\varphi_\xi$ is strictly convex (because $\xi$ is proper and non-trivial). 
Then, (\ref{stepclaimdist}) implies (\ref{claimdist}).

Recall that under $Q_{\xi, c}$ and conditionally given $Z^+_a$, $\abv (a, \cdot)$ has law $Q_{\xi, c}^{\circledast Z^+_a}$. Recall (\ref{asplitAred}) and (\ref{AbovAred}) and observe that 
for all $\widetilde{T} \in \bT_{{\rm edge}}$ and for all $a\in [0, \infty)$, 
$Z^+_a (R_A (\widetilde{T})) \geq 1$, implies $\widetilde{T}\in A^-$. Thus, for all  $n \geq 1$, and for  
all measurable functions $G\colon \bT \rightarrow [0, \infty)$, we get  
\begin{eqnarray}
\label{regenAred}
 B_n  &\!\!\!: =\!\!\! & \textstyle Q \left( \un_{\{ Z^+_a = n \} } G (\abv (a, \cdot) )\right)= \frac{1}{1-\alpha} 
 Q_{\xi, c} 
\left( \,  \un_{\{ Z^+_a \circ R_A= n \} } G (\abv (a, R_A (\cdot) ) ) \right) \nonumber \\
& \!\!\!=\!\!\!&\textstyle   \frac{1}{1-\alpha} Q_{\xi, c} 
\left( \un_{\{ \langle \cM_0 ( \abv (a, \cdot ) ), \un_{A^-}  \rangle = n \} } G ( R_A (\abv (a, \cdot ) ) ) \right)   \nonumber \\
& \!\!\!=\!\!\!&\frac{1}{1-\alpha}  \sum_{N \geq n} Q_{\xi, c} (Z^+_a \! = \!  N) \! 
\int_{\{ \sum_{1\leq i\leq N} \un_{A^-} (\widetilde{T}_i)= n \}}  
Q_{\xi, c}^{\otimes N} (d\widetilde{T}_1 \ldots d\widetilde{T}_N) \;  G \left(\underset{1\leq i\leq N\colon  \widetilde{T}_i\in A^-}{\circledast} 
R_A(\widetilde{T}_i) \right) \nonumber \\
& \!\!\!=\!\!\!& \frac{1}{1-\alpha}  \sum_{N \geq n} Q_{\xi, c} (Z^+_a \! = \!  N) \binom{N}{n} \alpha^{N-n} (1-\alpha)^n  
\int_{\bT^n}   Q^{\otimes n} (d\widetilde{T}^\prime_1 \ldots d\widetilde{T}^\prime_n) \, 
G\left( \widetilde{T}^\prime_1 \circledast \ldots \circledast \widetilde{T}^\prime_n \right) \nonumber \\
& \!\!\!=\!\!\!& \frac{1}{1-\alpha}  \sum_{N \geq n} Q_{\xi, c} (Z^+_a \! = \!  N) \binom{N}{n} \alpha^{N-n} (1-\alpha)^n   Q^{\circledast n} \left[ G \right] \; .
\end{eqnarray}
This proves that under $Q$, the law of $\abv (a, \cdot)$ conditionally given 
$Z^+_a$ is $Q^{\circledast Z^+_a}$. By Lemma \ref{regenene}, $Q$ is 
the law of a GW($\xi_*, c_*$)-real tree where $\xi_*$ is a proper conservative offspring distribution. 
We next set 
$w_*(a, \theta)= Q[\exp (-\theta Z^+_a) ]$ and $w(a, \theta) = Q_{\xi, c} [\exp (-\theta Z^+_a) ]$. By 
summing (\ref{regenAred}) over $n$ with $G= 1-e^{-\theta Z^+_a}$, we get  
\begin{equation}
\label{wwstar}
w_* (a, \theta)= 1- \frac{1-w\left( \, a \, , -\log (\alpha + (1-\alpha)e^{-\theta} ) \right) }{1-\alpha} \; .
\end{equation}
and by differentiating (\ref{inteqharris}) with respect to $a$ and to $\theta$, 
we get 
$$  -\frac{\partial_a w_* (a, \theta ) }{\partial_\theta w_* (a, \theta)}= c_* e^{\theta} \left( \varphi_{\xi_*} (e^{-\theta}) -e^{-\theta}  \right). $$
On the other hand, (\ref{wwstar}) implies that 
$$ -\frac{\partial_a w_*\!(a, \theta ) }{\partial_\theta w_*\!(a, \theta)}  =   -\frac{\alpha \!+\!  (1\!-\! \alpha)e^{-\theta}}{(1\!-\! \alpha)e^{-\theta}} \cdot 
\frac{\partial_a w(a,\!-\!\log (\alpha \!+\! (1\!-\! \alpha)e^{-\theta} )) }{\partial_\theta w(a, 
\!-\!\log (\alpha\! +\!  (1\!- \! \alpha)e^{-\theta} ))} \nonumber 
  =    \frac{ce^{\theta}}{1\!-\!\alpha}\!\left( \varphi_\xi (\alpha \!+\!  (1\!-\! \alpha)e^{-\theta})\!-\!\alpha \!-\!  (1\!-\! \alpha)e^{-\theta} \right)\!. 
$$
Thus, we get
 \begin{equation}
\label{harrisAred}
c_*  \left( \varphi_{\xi_*} (e^{-\theta}) -e^{-\theta}  \right)= \frac{c}{1-\alpha} \left( \varphi_\xi (\alpha +  (1\!-\! \alpha)e^{-\theta})-\alpha \!-\!  (1\!-\! \alpha)e^{-\theta} \right). 
\end{equation}
Since $\xi_*$ is proper, $\xi_* (1)= \varphi_{\xi_*}^\prime (0)= 0$. Thus, by differentiating 
(\ref{harrisAred}) and by letting $\theta$ go to $\infty$, we find
$c_*= c^{(\alpha)}$, as defined in (\ref{reduclaw}). This, combined (\ref{harrisAred}), implies 
$\xi_* = \xi^{(\alpha)}$. 
\end{pf}

The previous result and (\ref{AbovAred}) in Lemma \ref{Amoinsprop} applied to $a= 0$, immediately  implies the following statement for forests. Recall that 
$P_{\xi, c, \mu}= \sum_{n \in \bN} \mu (n) Q^{\circledast n}_{\xi, c}$. 
\begin{theorem}
\label{realforreduc} Let $\xi $ be a proper non-trivial conservative offspring distribution, let $c\in (0, \infty)$ and let $\mu$ be a probability distribution on $\bN$ such that $\mu (0) <1$.  
Let $A \subseteq \bT$ be hereditary.  
We set $\alpha= Q_{\xi, c} (R_A=\pnt)$ and we assume that $\alpha \in (0, 1)$. 
 Then, $R_A$ under $P_{\xi, c , \mu}$ has law $P_{\xi^{(\alpha)} , c^{(\alpha)}, \mu^{\{ \alpha \}}}$, where $\xi^{(\alpha)}$ and $ c^{(\alpha)}$ are given by {\rm(\ref{reduclaw})} and $\mu^{\{ \alpha\} }$ is given by 
\begin{equation}
\label{mureduclaw}
\varphi_{\mu^{\{ \alpha \}}} (r)= \varphi_\mu (\alpha + (1-\alpha) r) \, , \quad r\in [0, 1] \; .
\end{equation}
\end{theorem}

\section{Growth Processes.}
\label{Growthsec}

\subsection{Definition and characterisation.}
\label{defgrowthsec}

\subsubsection{Infinitely extensible offspring distributions.} 

We first briefly study the transform on offspring distributions that appears in (\ref{reduclaw}) in Theorem \ref{GWrealreduc}. Let $\xi $ be a {\it proper} offspring distribution, let $c \in (0, \infty)$ and let $\mu$ be another probability on $\bN$. 
Recall that $\varphi_\xi$ stands for its generating function and that $q_\xi$ is the smallest root of $\varphi_\xi (r)= r$. We introduce the following subset 
\begin{equation}
\label{defDxi}
D_\xi = \{ \alpha \in [0, 1)\colon   \alpha \leq \varphi_\xi (\alpha) \} \; .
\end{equation}
We easily check that for all $\alpha \in D_\xi$, $\varphi^\prime_{\xi} (\alpha) <1$, and we define $(\xi^{(\alpha)} , c^{(\alpha)}, \mu^{\{\alpha \}}) $ by setting for all $r\in [0, 1]$,  
\begin{equation}
\label{reddlaww}
c^{(\alpha)} \!=\! c (1- \varphi^\prime_\xi (\alpha))  , \! \quad  \varphi_{\xi^{(\alpha)}} (r)=  r + \frac{\varphi_\xi(  \alpha \!+\! (1\!- \!\alpha)r )\!- \!\alpha \!- \! (1\!- \!\alpha)r}{(1\!-\! \alpha)(1\!- \!\varphi_\xi^\prime(\alpha))} , \quad 
 \varphi_{\mu^{\{\alpha \}}} (r) \!= \! \varphi_\mu (\alpha + (1-\alpha)r).
\end{equation} 
Note that $\xi^{(\alpha)}$ is proper. Then $(\xi, c, \mu)$ can be recovered from $(\xi^{(\alpha)} , c^{(\alpha)}, \mu^{\{\alpha \}}) $ and $\alpha$ as follows. 
First note that $\varphi_{\xi^{(\alpha)}}$ can be analytically extended to  
$( \frac{-\alpha}{ 1-\alpha} \, , \, 1)$. Observe that the right limits of $\varphi_{\xi^{(\alpha)}}$ and $\varphi^\prime_{\xi^{(\alpha)}}$ at $ \frac{-\alpha}{ 1-\alpha}$ exist and are finite. Since $\xi$ is proper we easily get $(1-\varphi^\prime_\xi (\alpha))^{-1}= 1- \varphi^\prime_{\xi^{(\alpha)}} (\frac{-\alpha}{ 1-\alpha})  >1$. 
We therefore get for all $r\in [ 0\, , \, 1]$, 
  \begin{equation}
\label{propreinverse}
{\textstyle c= c^{(\alpha)}   (1- \varphi^\prime_{\xi^{(\alpha)}} (\frac{{-\alpha}}{{1-\alpha}}))} , \quad 
\varphi_\xi (r)= r+ \frac{ (1\!- \!\alpha) \, \varphi_{\xi^{(\alpha)}} (\frac{r-\alpha}{1-\alpha}) -(r\!-\!\alpha)}{1-\varphi^\prime_{\xi^{(\alpha)}}(\frac{-\alpha}{1-\alpha} )} , \quad \varphi_{\mu} (r)= \varphi_{\mu^{\{ \alpha \}}} {\textstyle(\frac{{r-\alpha}}{{1-\alpha}})}. 
\end{equation}
Next note that if  $\beta \in D_{\xi^{(\alpha)}}$, then 
\begin{equation}  
\label{compaaa}
 \left(  (\xi^{(\alpha)})^{(\beta )} , (c^{(\alpha)})^{(\beta )}, (\mu^{\{\alpha \}})^{\{ \beta \}} \right)= (\xi^{(\gamma)} , c^{(\gamma)}, \mu^{\{\gamma \}}) \quad \textrm{with $\gamma\in D_\xi$ and $1-\gamma = (1-\alpha)(1-\beta)$.}
\end{equation}
We now introduce the definition of infinitely extensible offspring distributions. 
\begin{definition}
\label{infdivoffdef} 
An offspring distribution $\xi$ 
is said to be {\it infinitely extensible} if there exists a sequence of proper offspring distributions $(\xi_n)_{ n \in \bN}$ and $\alpha_n \in D_{\xi_n}$, $n \in \bN$, such that $\xi= \xi_{n}^{{(\alpha_n)}}$ for all $n \in \bN$, 
and $\lim_{n \rightarrow \infty} \alpha_n = 1$.\cq 
\end{definition}
   For instance, if $\xi= \xi^{(\alpha)}$, for all $\alpha\in(0, 1)$, then, we easily prove that 
$\xi $ is a stable offspring distribution, namely, there exists $\gamma\in (1, 2]$ such that 
$  \varphi_\xi (r)= r+ \frac{1}{\gamma} (1-r)^\gamma $. Observe that stable offspring distributions  
are critical (i.e.$\;$their mean is equal to $1$), but their variance is infinite except in the binary case 
$\gamma = 2$. Stable offspring distributions occur in many contexts: see for instance Le Jan \cite{LJ91} and Neveu \cite{Ne2} for leaf-length erasure. See also \cite{DuLG}, Vatutin \cite{Va}
and Jakymiv \cite{Ya} for reduced trees. Note that there are offspring distributions $\xi$ that are not stable but such that $\xi= \xi^{(\alpha)}$ holds true for certain $\alpha \in (0, 1)$. 
    
 More generally, infinitely extensible offspring distributions are characterised by a branching mechanism, that is, by a function $\psi \colon  [0, \infty) \rightarrow \bR$ that is of the following 
L\'evy-Khintchine form 
\begin{equation}
\label{LKforminfdivoff}
\psi (\lambda) = {\bf a} \lambda + \frac{_1}{^2} {\bf b}
\lambda^2 + \int_{(0, \infty)} \left( e^{-\lambda x} -1+\lambda x \un_{\{ x<1\}} \right) \pi (dx) \; , 
\end{equation}
where ${\bf a} \in \bR$, ${\bf b}\in [0, \infty)$ and $\pi$ is a Borel measure on $(0, \infty)$ such that $\int_{(0, \infty)}
(1\wedge x^2) \pi (dx) < \infty$. A result similar to the theorem below  was proved in \cite[Theorem 4.2]{DuWi1} in the context of Bernoulli leaf percolation. The proof in the present framework is very similar, but since it demonstrates the crucial appearance of the branching mechanism $\psi$, we include a brief account of it here.
\begin{theorem}
\label{infdivoffstructure} Let $\xi$ be an infinitely extensible offspring distribution. Then, there exists a  function $\psi \colon  [0, \infty) \rightarrow \bR$ of the form {\rm(\ref{LKforminfdivoff})} such that $\psi^\prime (1)= 1$ and 
$$ \varphi_\xi (r)= r+ \psi (1-r) \, , \quad r \in [0, 1] \; .$$
Conversely, suppose that $\psi$ is of the L\'evy-Khintchine form {\rm(\ref{LKforminfdivoff})} and suppose that there exists $\lambda_0\in (0, \infty)$ such that $\psi (\lambda_0)>0$. Then, the function 
$$ r  \in [0, 1] \longmapsto r +  \frac{\psi (\lambda_0 (1-r))}{\lambda_0 \psi^\prime(\lambda_0) } $$
is the generating function of an infinitely extensible offspring distribution. 
\end{theorem}
\begin{pf} We first assume that $\xi$ is an infinitely extensible offspring distribution. Let $(\xi_n, \alpha_n)$, $n \in \bN$, be as in Definition \ref{infdivoffdef}. This implies that $\varphi_{\xi}$ 
can be extended analytically to the interval $I_n=( \frac{-\alpha_n}{1-\alpha_n}, 1)$. 
Moreover for all $ r \in I_n$, and all $k \geq 2$, 
$$ \varphi_\xi^{(k)} (r)= \frac{(1-\alpha_n)^{k-1}}{ 1-\varphi^\prime_{\xi_n} (\alpha_n)} \varphi^{(k)}_{\xi_n} 
\left( \alpha_n + (1-\alpha_n)r \right) $$ 
is positive. Since $\lim_{n \rightarrow \infty} \alpha_n = 1$, it make sense to set $\psi (\lambda)= \varphi_\xi (1-\lambda)-(1-\lambda)$, for all $\lambda \in [0, \infty)$. Then, $\psi$ is $C^\infty$ on $(0, \infty)$ and $(-1)^k \psi^{(k)} (\lambda) \geq 0$, for all $\lambda \in (0, \infty)$ and all $k \geq 2$. 
We then apply to $\psi^{(2)}$ Bernstein's theorem on completely monotone functions (see e.g.\ Feller \cite[Theorem XIII.7.2]{Fe}): there exists a Radon measure $\nu$ on $(0, \infty)$ and ${\bf b} \geq 0$ such that  
$ \psi^{(2)}( \lambda)= {\bf b} + \int_{(0,\infty)}  e^{-\lambda x} \nu (dx) $. 
Observe that $\psi(0)=0$ and $\psi^\prime (1)= 1 - \varphi^\prime_\xi (0)= 1 -\xi (1)= 1$. 
Then, we set $ \pi (dx)= x^{-2} \nu (dx)$ and ${\bf a}= 1-{\bf b}+ \int_{(0,\infty)} (e^{-x}-\un_{\{ x< 1\} } )x\pi(dx) $. 
 This easily entails that $\psi$ is of the form (\ref{LKforminfdivoff}). The proof of the converse result is straightforward: we leave the details to the reader.
\end{pf}

\subsubsection{Definition of growth processes. Characterization of their laws.} 
\begin{definition}
\label{growrealdef} ({\it Growth processes}) Let $(\Omega, \cG, \bP)$ be a probability space. For all $\lambda \in [0, \infty)$, let $\widetilde{\cF}_\lambda \colon  \Omega \rightarrow \bT$ be a GW($\xi^*_\lambda, c^*_\lambda ; \mu^*_\lambda$)-real forest as in Definition \ref{GWrealtreedef}. We say that $(\widetilde{\cF}_\lambda)_{\lambda\in[0,\infty)}$ is a growth process if for all $\lambda^\prime\ge\lambda$, there exists a hereditary property $A_{\lambda, \lambda^\prime} \subset \bT$ such that $\bP$-a.s.~for all $\lambda^\prime \geq \lambda$, $\widetilde{\cF}_\lambda\!=\!R_{A_{\lambda,\lambda^\prime}}(\widetilde{\cF}_{\lambda^\prime})$. \cq 
\end{definition}

\begin{remark}
\label{reggrow}
Note that if $(\widetilde{\cF}_\lambda)_{\lambda\in [0, \infty)}$ is a growth process, then by (\ref{herasegrow}), $\bP$-a.s.~for any $\lambda^\prime \geq \lambda$, $\widetilde{\cF}_\lambda \preceq \widetilde{\cF}_{\lambda^\prime}$. \cq  
\end{remark}
\begin{remark}
\label{regmodif} Suppose that the family $(\widetilde{\cF}_\lambda)_{\lambda\in [0, \infty)}$ satisfies 
the weaker consistency property that 
for all $\lambda^\prime \geq \lambda$, $\bP$-a.s.~$\widetilde{\cF}_\lambda\!=\!R_{A_{\lambda,\lambda^\prime}}(\widetilde{\cF}_{\lambda^\prime})$, then by (\ref{Aredembed}), we easily construct a modification of the process that is a growth process according to Definition \ref{growrealdef}. \cq  
\end{remark}
 
We next introduce specific one-parameter families of infinitely extensible offspring distributions that play a key r\^ole. To that end, we fix the following setting:  
\begin{enumerate}
\item[$-$] a branching mechanism $\psi\colon  [0, \infty) \rightarrow \bR$ of the form (\ref{LKforminfdivoff}), where we furthermore assume the following: 
\begin{equation}
\label{psiassum}
\textsl{allow deaths:} \quad \exists \lambda_0 \in (0, \infty)\colon \psi (\lambda_0) >0,  
\qquad\textsl{conservative:} \quad \int_{0+} \frac{dr}{(\psi (r))_-}= \infty \; ; 
\end{equation}
\item[$-$] a probability measure $\varrho$ on $[0, \infty)$ that is distinct from $\delta_0$. 
\end{enumerate}
The fact that $\psi$ is a convex function and that $\psi (0)= 0$ has several consequences. 
First, observe that $\psi$ has either one or two roots: we denote by $q$ the largest one, and since $\psi$ takes positive values eventually, we obtain $(q, \infty)= \{ \lambda \in [0, \infty) \colon  \psi (\lambda) >0\}$. 
Moreover, $\psi^{\prime} (0+)$ exists in $[-\infty, \infty)$. 
\begin{enumerate}
\item[$-$] If $\psi^{\prime} (0+)< 0$, we say that $\psi$ is {\it super-critical} (and $q>0$, necessarily). 
\item[$-$] If $\psi^{\prime} (0+) =0$, we 
say that $\psi$ is {\it critical} (and $q=0$, necessarily). 
\item[$-$]If $\psi^{\prime} (0+)> 0$, we say that $\psi$ is 
{\it sub-critical} (and $q= 0$, necessarily). 
\end{enumerate}
Next, for all $\lambda \in [q, \infty)$, we introduce two probability measures on $\bN$ that are denoted by $\xi_\lambda$ and $\mu_\lambda$ and we also define $c_\lambda \in (0, \infty)$, as follows: 
\begin{equation}
\label{xilambdef}
c_\lambda= \psi^\prime (\lambda) \, , \quad \varphi_{\xi_\lambda} (r)= r+ \frac{\psi \left( (1-r)\lambda  \right)}{\lambda \psi^\prime (\lambda)} \quad {\rm and} \quad  \varphi_{\mu_\lambda} (r)= \int_{[0, \infty)} e^{-\lambda y (1-r)} \varrho (dy) \; , \quad r\in [0, 1].
\end{equation}
More explicitly, we get $\xi_\lambda (0)= \frac{\psi (\lambda)}{\lambda \psi^\prime (\lambda) }$, $\xi_\lambda (1)= 0$, and for all $k\geq 2 $, 
\begin{equation}
\label{explicitxi}
 \xi_{\lambda} (k) = \frac{\lambda^{k-1} \! \left| \psi^{(k)} (\lambda)  \right|}{k! \,  \psi^\prime (\lambda)}= \frac{\lambda^{k-1}}{k! \,  \psi^\prime (\lambda)} \left( 
 {\bf b}  \un_{\{ k= 2\}} + \int_{(0, \infty)} x^{k} 
e^{-\lambda x} \, \pi (dx)  \right). 
\end{equation}
Recall (\ref{nonexpdiscr}): we easily check that $\xi_\lambda$ is conservative iff the second assumption in (\ref{psiassum}) holds true. 
The following theorem provides a useful classification of growth processes. 
\begin{theorem}
\label{structuregrowth} 
Let $(\widetilde{\cF}_\lambda)_{\lambda \in [0, \infty)}$ be a growth process as in 
Definition \ref{growrealdef}. By {\rm(\ref{herasegrow})}, $\bP$-a.s.~$\lambda \mapsto Z^+_0 (\widetilde{\cF}_\lambda)$ is non-decreasing and we set $Z= \lim_{\lambda \rightarrow \infty}  Z^+_0 (\widetilde{\cF}_\lambda)$.   
Then, only the two following cases occur. 

\begin{enumerate}
\item[${\rm (I)}$] If $\bP (Z <\infty) = 1$, there exists a proper and conservative offspring distribution $\xi$, there exists $c\in (0, \infty)$ and there exists a probability measure $\mu$ distinct from $\delta_0$, such that 
\begin{equation}
\label{dilimun}
\forall \, k \in \bN \, , \quad \lim_{\lambda \rightarrow \infty} \xi^*_{\lambda} (k)= \xi (k) \; , \quad   
\lim_{\lambda \rightarrow \infty} \mu^*_{\lambda} (k)= \mu (k) \quad {\rm and} \quad \lim_{\lambda \rightarrow \infty} c^*_{\lambda} = c . 
\end{equation}
Moreover, there exists a unique non-increasing function $\alpha\colon [0, \infty) \rightarrow D_\xi$ such that 
$$ \lim_{\lambda \rightarrow \infty} \alpha_{\lambda} = 0 \quad {\rm and} \quad \forall \lambda \in [0, \infty) \, , \quad   \xi^*_\lambda = \xi^{(\alpha_\lambda)} \; , \quad  c^*_\lambda = c^{(\alpha_\lambda)} \quad {\rm and} \quad \mu^*_\lambda = \mu^{ \{  \alpha_\lambda  \} } \; , $$
where for all $\alpha \in D_\xi$, $(\xi^{(\alpha)}, c^{(\alpha)}, \mu^{\{ \alpha \}})$ is defined by 
{\rm(\ref{reddlaww})}. 
\item[${\rm (II)}$] If $\bP (Z= \infty) >0$, there exists a triplet $(\psi , \varrho, \beta)$ that satisfies the following. 
\begin{enumerate}
\item[$\rm(i)$] $\psi$ is a function of the Lévy-Khintchine form {\rm(\ref{LKforminfdivoff})} that satisfies {\rm(\ref{psiassum})}.  
\item[$\rm(ii)$] $\varrho$ is a probability measure $[0, \infty)$ distinct from $\delta_0$, 
\item[$\rm(iii)$] $\beta\colon  [0, \infty) \rightarrow [q, \infty)$ is a non-decreasing function such that $\lim_{\lambda \rightarrow \infty} \beta(\lambda)= \infty$. Here, $q$ stands for the largest root of $\psi$.  
\item[$\rm(iv)$] For all $\lambda \in [0, \infty)$, 
$$ \xi^*_\lambda = \xi_{\beta (\lambda)} \, , \quad c^*_\lambda= c_{\beta (\lambda)} \quad {\rm and} \quad  \mu^*_\lambda= \mu_{\beta (\lambda)} $$
where for all $\lambda \in [q, \infty)$, $(\xi_\lambda, c_\lambda, \mu_\lambda) $ 
is defined by ${\rm (\ref{xilambdef})}$.  
\end{enumerate}

Moreover, if another triplet $( \psi^*, \varrho^* ,\beta^*)$ satisfies {\rm (i)--(iv)}, there exists a constant $\kappa \in (0, \infty)$ such that 
$\varrho^* (dx)= \varrho (dx/ \kappa)$, and for all $ \lambda \in [0, \infty)$,  $\psi^* (\lambda)= \psi (\kappa \lambda)/ \kappa$ and $\beta^*( \lambda)=\beta (\lambda) / \kappa $.\pagebreak[4] 
\end{enumerate}
\end{theorem}
\begin{pf} Let $\widetilde{\cT}_\lambda$ be a GW($\xi^*_\lambda, c^*_\lambda$)-real tree: its law is then $Q_{\xi^*_\lambda, c^*_\lambda}$. We keep the same notation as in Definition \ref{growrealdef}.  
For all $\lambda^\prime \geq \lambda$, we recall from (\ref{Aminus}) the definition of the hereditary property $A_{\lambda, \lambda^\prime}^-$ and we set 
$$ \alpha_{\lambda , \lambda^\prime} := \bP \left( \widetilde{\cT}_{\lambda^\prime} \notin A^-_{\lambda, \lambda^\prime} \right) = Q_{\xi^*_{\lambda^\prime} , c^*_{\lambda^\prime}} \left( \bT \backslash A^-_{\lambda, \lambda^\prime} \right) \; .$$
By Theorem \ref{realforreduc}, 
$((\xi^*_{\lambda^\prime})^{(\alpha_{\lambda, \lambda^\prime})} , (c^*_{\lambda^\prime})^{(\alpha_{\lambda, \lambda^\prime})}, 
(\mu^*_{\lambda^\prime})^{ \{ \alpha_{\lambda, \lambda^\prime} \}})= 
(\xi_\lambda^* , c_\lambda^*, \mu_\lambda^* )$. This first implies that $\lambda \mapsto c^*_\lambda$ is non-decreasing. Recall from the definition of growth processes that $\mu^*_0 (0) < 1$. This entails that 
$\mu^*_\lambda (0) < 1$, for all $\lambda \in [0, \infty)$, because 
$ \varphi_{\mu^*_\lambda} (0) \leq  \varphi_{\mu^*_\lambda} (\alpha_{0, \lambda})= \varphi_{\mu^*_0} (0) = \mu^*_0 (0) < 1 $.
Next, (\ref{propreinverse}) easily entails 
 \begin{equation}
\label{propinv}
\varphi_{\xi^*_{\lambda^\prime}} (r)= r+  \frac{c^*_\lambda}{c^*_{\lambda^\prime}} \cdot 
\left( (1\!- \!\alpha_{\lambda, \lambda^\prime}) \, \varphi_{\xi^*_\lambda} \left(\frac{_{r-\alpha_{\lambda, \lambda^\prime}}}{^{1-\alpha_{\lambda, \lambda^\prime}}}\right) -(r\!-\!\alpha_{\lambda, \lambda^\prime})  \right)
 \quad {\rm and} \quad \varphi_{\mu^*_{\lambda^\prime} } (r)= \varphi_{\mu^*_\lambda} \left( \frac{_{r-\alpha_{\lambda, \lambda^\prime}}}{^{1-\alpha_{\lambda, \lambda^\prime}}} 
\right) \; .
\end{equation}
Moreover, by (\ref{compaaa}), $(1-\alpha_{\lambda, \lambda^{ \prime}}) (1- \alpha_{\lambda^\prime, \lambda^{\prime \prime}})= 1- \alpha_{\lambda, \lambda^{\prime \prime}}$, for all $\lambda^{\prime \prime} \geq \lambda^\prime \geq \lambda$. This implies that $\alpha_{\lambda, \lambda^\prime}$ is non-decreasing in $\lambda^\prime$ and non-increasing in $\lambda$. Then, for all $\lambda \in [0, \infty)$, we set $\alpha_\lambda = \lim_{\lambda^\prime \rightarrow \infty } \alpha_{\lambda, \lambda^\prime}$. 

\noi $\ \,\bullet\,$ We first assume that $\bP (Z< \infty )= 1$. We denote by 
$\mu$ the law of $Z$. Since $\mu^*_\lambda$ is the law of 
$Z^+_{0} (\widetilde{\cF}_\lambda)$, by definition, we get $\lim_{\lambda \rightarrow \infty} 
\mu^*_{\lambda} (k)= \mu (k)$, for all $k\in \bN$.

We next claim that for all $\lambda\in[0,\infty)$, $\alpha_\lambda<1$. We argue by contradiction: let us suppose that 
$\alpha_\lambda = 1$. Then, for all $\lambda^{\prime \prime} \geq 
\lambda^{\prime} \geq \lambda$ and all $r\in [0, 1]$, we get 
$$ \varphi_{\mu^*_\lambda} (r) = 
\varphi_{\mu^*_{\lambda^{\prime \prime}} } \left( \alpha_{\lambda, \lambda^{\prime \prime}} (1-r)+r \right) \geq 
\varphi_{\mu^*_{\lambda^{\prime \prime}} } \left( \alpha_{\lambda, \lambda^{\prime }} (1-r)+r \right)  \; .$$  
But $\lim_{ \lambda^{\prime \prime} \rightarrow \infty} \varphi_{\mu^*_{\lambda^{\prime \prime}} } \left( \alpha_{\lambda, \lambda^{\prime }} (1\!-\!r)+r \right) = 
\varphi_{\mu } ( \alpha_{\lambda, \lambda^{\prime }} (1\!-\!r)+r) $ and 
$\lim_{ \lambda^{\prime } \rightarrow \infty} \varphi_{\mu } \left( \alpha_{\lambda, \lambda^{\prime }} (1\!-\!r)+r \right)= \varphi_\mu (1)= 1$, since we have supposed $\alpha_\lambda= 1$. This implies that $\mu^*_\lambda (0)= 1$, which is impossible as already proved. Note that this argument also entails that $\mu (0) < 1$. 

 An elementary compactness argument shows that there is a {\it sub-probability} measure $\xi$ on $\bN$ and a 
sequence $\lambda_n \in [\lambda , \infty)$, $n \in \bN$, increasing to $\infty$, such that $\lim_{n \rightarrow \infty} \xi_{\lambda_n}^* (k)= \xi (k)$, for all $k\in \bN$. 
This convergence entails $\lim_{n \rightarrow \infty} \varphi_{\xi_{\lambda_n}^*} (r)=
 \varphi_{\xi} (r)$, for all $r\in [0, 1)$. 
 Next observe that since  $\alpha_\lambda < 1$, (\ref{propinv}) implies that 
 $\varphi_{\xi^*_\lambda}$ can be extended analytically 
 to the interval $(\frac{-\alpha_\lambda}{1-\alpha_\lambda} , 1)$. 
Thus, for all $0<r<1$, in (\ref{propinv}), we can let $\lambda^\prime \rightarrow \infty$ (along the sequence $\lambda_n$) and there are two cases to consider:  
since $\lambda^\prime \mapsto c^*_{\lambda^\prime}$ is non-decreasing, if $c:=\lim_{\lambda^\prime \rightarrow \infty} 
c^*_{\lambda^\prime} < \infty$, then we get 
\begin{equation}
\label{obred}
 \varphi_\xi (r) = r+  \frac{c^*_\lambda}{c} \cdot 
\left( (1- \alpha_{\lambda}) \, \varphi_{\xi^*_\lambda} \left(\frac{r-\alpha_{\lambda}}{1-\alpha_{\lambda}}\right) -(r-\alpha_{\lambda})  \right), 
\end{equation}
and if $\lim_{\lambda^\prime \rightarrow \infty} 
c^*_{\lambda^\prime} = \infty$, then we get $\varphi_\xi (r) = r $. However, the last possibility would imply that $\xi (1)= 1$, which is impossible because $\xi$ is the limit of proper offspring distributions. Thus (\ref{obred}) holds true. Note that (\ref{obred}) entails that $\varphi_{\xi} (1)= 1$ and that 
$\xi$ is a proper offspring distribution. Moreover, $\xi$ is completely determined by (\ref{obred}). This implies (\ref{dilimun}).  
Also note that (\ref{obred}) easily entails that $\alpha_\lambda \in D_\xi $ and that  
that $(\xi^{(\alpha_\lambda)} ,  c^{(\alpha_\lambda)} , \mu^{\{ \alpha_\lambda \}} )= 
(\xi^*_\lambda , c^*_\lambda, \mu^*_\lambda)$. 

  The function $\lambda \in [0, \infty) \mapsto \alpha_\lambda \in D_\xi$ is clearly non-increasing. Let us prove uniqueness: suppose that 
$\varphi_{\mu^*_\lambda} ( r)= \varphi_\mu (\alpha_\lambda + (1-\alpha_\lambda) r)=\varphi_\mu (\gamma_\lambda + (1-\gamma_\lambda) r)$. Since $\mu (0) < 1$, $\varphi_\mu$ is strictly increasing, and we get $\alpha_\lambda = \gamma_\lambda$. 
We next set $\alpha = \lim_{\lambda \rightarrow \infty} \alpha_\lambda $. Recall that 
$\varphi_{\mu^*_\lambda} (0) =  \varphi_\mu (\alpha_\lambda) \geq \varphi_{\mu} (\alpha)$, which implies $\varphi_\mu (0) \geq \varphi_\mu (\alpha)$ as $\lambda \rightarrow \infty$. Thus, 
$\alpha= 0$, since $\varphi_\mu$ is strictly increasing. This completes the proof of Case (I). 

\medskip

\noi $\ \,\bullet\,$ We now assume that $\bP (Z = \infty) >0$. We claim that $\alpha_1= \lim_{\lambda \rightarrow \infty} \alpha_{1, \lambda}= 1$. We argue by contradiction: let us suppose that $\alpha_1< 1$. We fix $r <1$. Since $\mu^*_1= (\mu^*_\lambda )^{\{ \alpha_{1, \lambda} \} }$ 
for all $\lambda \in (1, \infty)$, we get  
$$ \varphi_{\mu^*_1} (r) \! = \!  \bE \left[ (\alpha_{1,\lambda } + (1\!-\!\alpha_{1, \lambda} ) 
r )^{Z^+_0 (\widetilde{\cF}_\lambda )} \right] \leq  \bE \left[ (\alpha_{1 } + (1\!-\!\alpha_{1} ) 
r )^{Z^+_0 (\widetilde{\cF}_\lambda )} \right]  \underset{{\lambda \rightarrow \infty}}{\longrightarrow}  \bE \left[ (\alpha_{1} + (1\!-\!\alpha_{1} ) 
r )^{Z} \un_{\{ Z < \infty \} } \right].$$
We now let $r$ go to $1$ and we get $1= \varphi_{\mu^*_1} (1)\leq \bP (Z< \infty)$, which contradicts the assumption $\bP (Z= \infty) >0$.

   Since $\lim_{\lambda \rightarrow \infty} \alpha_{1, \lambda}= 1$ and since $\xi^*_1= (\xi^*_\lambda)^{(\alpha_{1, \lambda})}$, $\xi^*_1$ is infinitely extensible. Theorem \ref{infdivoffstructure} implies that there exists $\psi_1$ of the form (\ref{LKforminfdivoff}) that satisfies 
$\psi_1^{\prime} (1) = 1$ and  $\varphi_{\xi^*_1} (r) = r+ \psi_1 (1-r)$, $ r \in [0 , 1] $. We  denote by $q$ the largest root of $\psi_1$. Since $\psi_1 (1)= \xi^*_1 (0) \geq 0$, we get $q\leq 1$. Then, observe that 
since $\xi^*_1$ is conservative, $\psi_1$ satisfies $\int_{0+} du/ (\psi_1(u))_-= \infty$. Thus, $\psi_1$ satisfies (\ref{psiassum}).

 We next prove that $\mu^*_1$ is a mixture of Poisson distributions: we set $g(r)= \varphi_{\mu^*_1} (1-r)$ and observe that $g(r)= \varphi_{\mu^*_\lambda} (1-(1-\alpha_{1, \lambda})r)$, for all $\lambda \in [1, \infty)$. Since $\lim_{\lambda \rightarrow \infty} \alpha_{1, \lambda}= 1$, $g$ can be analytically extended to $(0 , \infty) $, and we easily see that $g$ is completely monotone. Bernstein's theorem  entails that there exists a probability measure $\varrho$ on $[0, \infty)$ such that $\varphi_{\mu^*_1} (r)=
 \int e^{-y(1-r)} \varrho (dy) $, for all $r\in [0, 1]$. Note that $\varrho (\{0 \})=\mu^*_1 (0)  < 1$. Thus $\varrho$ is distinct from $\delta_0$. 

We now set $\psi = c^*_1 \psi_1$ and we define $\beta \colon  [0, \infty) \rightarrow [q, \infty)$ by setting 
$$ \textrm{$\beta (\lambda)= 1-\alpha_{\lambda, 1}$,\quad if $\lambda \in [0, 1)$}, \qquad  \beta (1)= 1, \qquad {\rm and} \quad  \textrm{$\beta (\lambda) = \frac{1}{1- \alpha_{1, \lambda}} $,\quad
 if $\lambda \in (1, \infty)$.}$$
Recall notation $(\xi_\lambda, c_\lambda , \mu_\lambda)$ from (\ref{xilambdef}). 
Then, for all $\lambda \in [0, \infty)$, we check that $ \xi^*_\lambda = \xi_{\beta (\lambda)}$ , $c^*_\lambda= c_{\beta (\lambda)}$ and $\mu^*_\lambda= \mu_{\beta (\lambda)}$ by applying 
(\ref{reddlaww}) to $\alpha = \alpha_{\lambda , 1}$ and $(\xi, c, \mu)= (\xi^*_1,  c^*_1, \mu^*_1)$ if $\lambda < 1$, and  by applying (\ref{propreinverse}) to $\alpha= \alpha_{1, \lambda}$ and 
$(\xi, c, \mu)= (\xi^*_\lambda,  c^*_\lambda, \mu^*_\lambda)$, if $\lambda >1$. 
We have proved that $(\psi, \varrho, \beta)$ satisfies (i)--(iv).  
The last point of the theorem is not difficult to check: we leave it to the reader. 
\end{pf}

\subsection{Convergence to L\'evy forests.}
\label{Levytreeconvsec}
In this section we prove under a necessary and sufficient condition that 
every growth process converges almost surely to a 
Lévy forest represented as a $\bT$-valued random variable. To that end, we first consider 
the convergence in law of the branching processes to apply the general criterion 
in Theorem \ref{tightness2}.

\paragraph{Continuous-state branching process.} Continuous-state branching processes (CSBPs for short) are the continuous analogue, in time and space, of Galton-Watson branching processes. They were introduced by Jirina \cite{Ji} and Lamperti \cite{La1, La2, La3}. We first recall standard results on CSBPs whose proof can be found in Silverstein \cite{Sil68} and Bingham \cite{Bi2}.

We fix a branching mechanism $\psi\colon  [0, \infty) \rightarrow \bR$ of the form (\ref{LKforminfdivoff}) and we assume that $\psi$ satisfies (\ref{psiassum}). We also fix a probability measure $\varrho $ on $[0, \infty)$ that is distinct from $\delta_0$. 
A $[0, \infty)$-valued Feller processes
$Z^\varrho =(Z^\varrho_a)_{ a \in [0, \infty) }$ defined on the probability space $(\Omega, \cG, \bP)$ is a 
continuous-state branching process with branching mechanism $\psi$ and initial distribution $\varrho $ (a CSBP($\psi, \varrho$) for short)
if its transition kernels are characterised by the following:  
$$ \bE \left[\left. e^{-\theta Z^\varrho_{b+a}} \; \right| \; Z^\varrho_b \right] = \exp \left(-Z^\varrho_b \, u(a, \theta) \right) \; , \quad a,b, \theta \in [0, \infty) \; , $$
where for all $\theta\in [0, \infty) $, the function $a\mapsto u(a, \theta )$ is the solution of the differential equation $\partial u (a, \theta)/\partial a  = - \psi (u(a, \theta))$, with $u(0, \theta)= \theta$. 
Under our assumptions on $\psi$, 
the differential equation satisfied by $u (\cdot , \theta)$ has a unique nonnegative solution that is defined 
on $[0, \infty)$ for all $\theta \in [0, \infty)$. In particular, note 
that the null function is the unique solution for $\theta= 0$.

    Let us rewrite the equation characterising $u$ in a more convenient way. Recall that $q$ stands for the largest root of $\psi$. If $\theta \in \{ 0, q\} $, then $u(a, \theta)= \theta $ for all $a\in [0,  \infty)$. Observe that if $ \theta >q$ (resp.~$\theta < q$) then $a \mapsto u(a, \theta)$ is a decreasing 
(resp.~increasing) function converging to $q$ as $a$ tends to infinity. An easy change of variables then entails 
\begin{equation}
\label{equaintegrale}
\forall \, a \in [0, \infty)\, , \quad  \forall \, \theta\in [0, \infty) \backslash \{ 0, q \} \, ,  \quad \int_{u(a, \theta)}^\theta \frac{dr}{\psi (r)} = a \; .
\end{equation}
We derive from the equation governing $u$ and from 
basic properties of Laplace transforms that for all $a\in [0, \infty)$, $Z^\varrho_a$ is integrable iff both $\psi^\prime (0+)$ and $\int_{[0, \infty)} y \varrho (dy)$ are finite quantities. In this case we get 
\begin{equation}
\label{expecZZ}
 \bE [Z^\varrho_a ] =e^{-\psi^\prime (0+) a } \int_{[0, \infty)} y \varrho (dy)  \; .
\end{equation}

  We next easily derive from (\ref{equaintegrale}) that 
$$ \bP \left( \lim_{a \rightarrow \infty}Z^\varrho_a = 0\right) +  \bP \left( \lim_{a \rightarrow \infty} Z^\varrho_a = \infty  \right)= 1 \quad {\rm and} \quad   \bP \left( \lim_{a \rightarrow \infty} Z^\varrho_a= \infty \right)= 1 - 
\int_{[0, \infty)} e^{-qx} \varrho (dx)  \;  .$$ 
Thus, $\bP \left( \lim_{a \rightarrow \infty} Z^\varrho_a= \infty \right) >0$ iff $q>0$. Namely, this  
happens only in the super-critical cases. 
It is easy to derive from (\ref{equaintegrale}) that 
\begin{equation}
\label{psextinct}
 \bP\left(\lim_{a \rightarrow \infty} Z^\varrho_a = 0\right)= \bP\left( \exists a \geq 0 \colon   Z^\varrho_a = 0\right)\quad\iff\quad\int^\infty \frac{dr}{\psi (r)} < \infty . 
\end{equation}
Condition (\ref{psextinct}) is called the {\it finite time extinction assumption}.  
Now observe that for all $a \in (0, \infty)$, $\theta \mapsto u(a, \theta)$ is an increasing function and 
under (\ref{psextinct}), 
\begin{equation} 
\label{vvvdefini}
\textrm{$v(a):=\lim_{\theta \rightarrow \infty} u(a, \theta) $ exists and satisfies} \quad 
\int_{v(a)}^\infty \frac{dr}{\psi (r)} = a \; , \quad a \in (0, \infty) \; .
\end{equation}
Observe that the function $v\colon  (0, \infty) \rightarrow (q, \infty)$ is decreasing and one-to-one.  
We next define the extinction time of $Z^\varrho$ by 
\begin{equation}
\label{extinctdef}
\cE (Z^\varrho)= \inf \left\{ a \in (0, \infty)  \colon   Z^\varrho_a = 0 \right\}  \; ,
\end{equation} 
with the convention $\inf \emptyset = \infty$. We then easily get  
\begin{equation}
\label{tspdextinct}
\bP \left( \, \cE (Z^\varrho) \leq a \; \right) = \int_{[0, \infty)} \varrho (dy) \, \exp \left( -v(a)y \right) \; , 
\quad a \in (0, \infty)\; .
\end{equation}
We refer to \cite{Bi2} for more details on CSBPs. 

\paragraph{Convergence in law to Lévy forests.} Let us fix a branching mechanism $\psi$ 
is of the form (\ref{LKforminfdivoff}) that furthermore satisfies (\ref{psiassum}).
We also fix 
$\varrho$, a probability measure on $[0, \infty)$ that is distinct from $\delta_0$. Recall that $q$ stands for the largest root of $\psi$. For all $\lambda \in [q, \infty)$, recall from (\ref{xilambdef}) the definition of 
$(\xi_\lambda, c_\lambda, \mu_\lambda)$. 

We denote by $\widetilde{\cF}_\lambda\colon \Omega \rightarrow \bT$ a GW($\xi_\lambda, c_\lambda; \mu_\lambda$)-real forest. {\it We do not need to assume that it is a growth process}: here, we are only interested in the convergence in law of these random trees when $\lambda$ goes to infinity.
We also denote by $\widetilde{\cT}_\lambda\colon  \Omega \rightarrow \bT$, a 
GW($\xi_\lambda, c_\lambda$)-real tree.
Namely the law of $\widetilde{\cT}_\lambda$ is $Q_{\xi_\lambda, c_\lambda}$ and the law of 
$\widetilde{\cF}_\lambda$ is $P_{\xi_\lambda, c_\lambda, \mu_\lambda}$, as in Definition \ref{GWrealtreedef}. 

We first consider the corresponding branching processes. 
As an easy consequence of (\ref{inteqharris}) and (\ref{equaintegrale}), for all $a, \theta\in [0, \infty)$, we get 
\begin{equation}
\label{niveautree}
w_\lambda (a, \theta):= \bE \left[ \exp (-\theta Z^+_{a} (\widetilde{\cT}_\lambda) )\right] = 
Q_{\xi_\lambda, c_\lambda} 
\left[ \exp (-\theta Z_a^+ ) \right]= 1-\frac{u\left( a, \lambda (1\!- \! e^{-\theta} ) \right)}{\lambda} 
\end{equation}
and (\ref{forestharris}) also implies 
\begin{equation}
\label{niveauforestt}
 \bE \left[ \exp ( -\theta  Z^+_{a} (\widetilde{\cF}_\lambda) )  \right]= P_{\xi_{\lambda}, c_\lambda, \mu_\lambda} 
\left[ \exp (-\theta Z^+_a ) \right] 
= \int_{[0, \infty)} \exp\left(- y \, u(a, \lambda (1\!-\!e^{-\theta})) \right) \varrho (dy) .
\end{equation}
This easily implies that for all $a\in [0, \infty)$, $\frac{1}{\lambda} Z^+_a (\widetilde{\cF}_\lambda) \rightarrow Z^\varrho_a $ in distribution on $[0, \infty)$. We next use a result due to 
Helland \cite[Theorem 6.1]{He} that shows that the following convergence actually holds in distribution 
in the space of cadlag functions $\bD ([0, \infty), \bR)$ equipped with Skorohod's  metric: 
\begin{equation}
\label{convbranch}\textstyle
\left(\frac{1}{\lambda}Z^+_{a} (\widetilde{\cF}_\lambda) \right)_{a \in [0, \infty)} \;  \longrightarrow \; 
\left(\, Z_a^\varrho \, \right)_{a \in [0, \infty)} \qquad\mbox{as $\lambda \rightarrow \infty$.} 
\end{equation}
In the following lemma, we compute the law of the $h$-leaf-length erased forest $R_h (\widetilde{\cF}_\lambda)$. 
\begin{lemma}
\label{heraselim} Let $\psi$ be a function of the form   
{\rm(\ref{LKforminfdivoff})} that satisfies (\ref{psiassum}). Let 
$\varrho$ be a probability measure on $[0, \infty)$, distinct from $\delta_0$. Let 
$\widetilde{\cF}_\lambda$ be as above. Then, the following holds true. 
\begin{enumerate}
\item[$\rm(i)$] For all $h \in (0, \infty)$, all $a\in [0, \infty)$ and all $\lambda \in [q, \infty)$, 
\begin{equation}
\label{halphah}
Q_{\xi_\lambda , c_\lambda} (\Gamma \leq  h ) =  1\!- \! \frac{{u(h, \lambda)} }{{\lambda}} 
\quad {\rm and} \quad R_h (\widetilde{\cF}_\lambda) \overset{\textrm{\rm(law)}}{=} \widetilde{\cF}_{u(h, \lambda)} \; .
\end{equation}
Moreover, 
\begin{equation}
\label{hefflambda}
\bE \left[ \exp (- \theta Z^{(h)}_{a} (\widetilde{ \cF}_\lambda) ) \right]
= \int_{[0, \infty)} \varrho (dy) \, \exp \left( - y \, u \left( a, u(h,\lambda) (1\!- \! e^{-\theta}) \right)  \right) .
\end{equation}

\item[$\rm(ii)$] The laws $P_{\xi_\lambda , c_\lambda , \mu_\lambda}$ of the GW-real forests 
$\widetilde{\cF}_\lambda$ are tight on $\bT$ as $\lambda\rightarrow \infty$ 
if and only if $\psi $ satisfies 
$$ \int^\infty \frac{dr}{ \psi (r)} < \infty \; .$$

\item[$\rm(iii)$] Assume that $\int^{\infty} dr/ \psi (r) < \infty$ and recall from {\rm(\ref{vvvdefini})} the definition of $v$. Then, 
$$ \left(  Z^{(h)}_{a} (\widetilde{ \cF}_\lambda) \, ;  \,  R_h (\widetilde{\cF}_\lambda) \right) \longrightarrow \left( Z^+_{a} (\widetilde{\cF}_{v(h)}) \, ; \, \widetilde{\cF}_{v(h)} \right)$$ 
weakly on $\bN \times \bT$ as $\lambda \rightarrow \infty$. Consequently, 
\begin{equation}
\label{hefflim}
 \lim_{\lambda \rightarrow \infty} 
\bE \left[ \exp (- \theta Z^{(h)}_{a} (\widetilde{ \cF}_\lambda) ) \right]
=  \int_{[0, \infty)} \varrho (dy) \exp \left(- y  \,  u \left( a, v(h)(1\!-\! e^{-\theta}) \right)  \right) .
\end{equation}
\end{enumerate}
\end{lemma}
\begin{pf} Recall that $A_h = \{ \Gamma \geq h\}$ is a hereditary property on $\bT$. Thus, $R_h (\widetilde{\cT}_\lambda)= R_{A_h} (\widetilde{\cT}_\lambda)$. We then see that 
$ \alpha := Q_{\xi_\lambda , c_\lambda} ( \bT \backslash A_h^-)= Q_{\xi_\lambda , c_\lambda} (\Gamma \leq  h )$. Then, recall that $Q_{\xi_\lambda , c_\lambda} (\Gamma \leq  h ) = \lim_{\theta \rightarrow \infty} w_\lambda (h , \theta)$. Thus, (\ref{niveautree}) entails the equality in (\ref{halphah}). 
Then, Theorem \ref{realforreduc} entails that $R_{A_h} (\widetilde{\cF}_\lambda)$ is a 
GW($\xi_\lambda^{(\alpha)}, c_\lambda^{(\alpha)}; \mu_\lambda^{\{ \alpha\}}$)-real forest, and we easily  see that 
$$(\xi_\lambda^{(\alpha)}, c_\lambda^{(\alpha)}, \mu_\lambda^{\{ \alpha\}})= (\xi_{u(h, \lambda) } , c_{u(h, \lambda) }, \mu_{u(h, \lambda) })\; , $$
which proves the point $\rm(i)$.
 
We next prove $\rm(ii)$. Theorem \ref{tightness} asserts that the laws of the $\widetilde{\cF}_\lambda$ are  are tight on $\bT$ as $\lambda \rightarrow \infty$ iff for each $a, h\in (0, \infty)$ the laws of the $Z^{(h)}_a(\widetilde{\cF}_\lambda)$, $\lambda \in [q, \infty)$, are tight on $\bN$ as  
$\lambda \rightarrow \infty$.  If $\int^\infty dr / \psi (r) = \infty$, (\ref{equaintegrale}) entails that 
$\lim_{\lambda \rightarrow \infty} u(h, \lambda)= \infty$, and since $\varrho \neq \delta_0$, 
(\ref{hefflambda}) implies that the laws of the $Z^{(h)}_a(\widetilde{\cF}_\lambda)$ \textsl{are not} 
tight on $\bN$ as $\lambda \rightarrow \infty$. If $\int^\infty dr / \psi (r) < \infty$, then $\lim_{\lambda \rightarrow \infty} u(h, \lambda)= v(h)< \infty $, and (\ref{hefflambda}) implies (\ref{hefflim}), which proves that the laws of the $Z^{(h)}_a(\widetilde{\cF}_\lambda)$ are 
tight on $\bN$, and the proof of $\rm(ii)$ is complete.

 We next observe $\lim_{\lambda \rightarrow \infty} (\xi_{u(h,\lambda)}, c_{u(h,\lambda)}, \mu_{u(h,\lambda)})= (\xi_{v(h)}, c_{v(h)}, \mu_{v(h)})$. Then, Lemma \ref{poubelle} applies, which entails $\rm(iii)$. 
\end{pf}

\begin{theorem}
\label{growthconvergence}
Let $\varrho$ be a Borel probability measure on $[0, \infty)$ distinct from $\delta_0$. 
Let $\psi\colon  [0, \infty) \rightarrow \bR$ be a branching mechanism of Lévy-Khintchine form {\rm(\ref{LKforminfdivoff})}. We assume that $\psi$ satisfies {\rm (\ref{psiassum})} and {\rm (\ref{psextinct})}. 
We denote by $q$ is largest root. 
Let $(\Omega, \cG, \bP)$ be a probability space. For all $\lambda  \in [q, \infty)$, let 
$\widetilde{\cF}_\lambda\colon  \Omega \rightarrow \bT$ be a ${\rm GW}(\xi_\lambda , c_\lambda; \mu_\lambda)$-real forest where $(\xi_\lambda , c_\lambda, \mu_\lambda)$ is given by {\rm (\ref{xilambdef})}. 
Then the following joint convergence holds true in distribution in $\bD ([0, \infty) , \bR) \times \bT${\rm :} 
\begin{equation}
\label{convjointe}\textstyle
\left( \; \left( \frac{1}{\lambda}Z^+_{a} (\widetilde{\cF}_\lambda) \, \right)_{ a \in [0, \infty )}   \; ; \;  \widetilde{\cF}_\lambda \;  \right) \; \; \underset{{\lambda \rightarrow \infty}}{\longrightarrow } \; \; \left(  \; 
( Z^\varrho_{a} )_{ a\in [0, \infty)} \; ; \; \widetilde{\cF } \; \right) \, , 
\end{equation}
where the limit is as follows. 
\begin{enumerate}
\item[\rm(i)] The limiting forest $\widetilde{\cF}$ is a $(\psi, \varrho)$-{\bf Lévy forest}. 
We denote by $\LevP_{\psi, \varrho}$ its distribution on $\bT$ which only depends on $\psi$ and $\varrho$. Furthermore, we have $\LevP_{\psi^*, \varrho^*} = \LevP_{\psi, \varrho}$ iff there exists $\kappa \in (0, \infty)$ such that  $\psi^* (\lambda)= \psi (\kappa \lambda)/ \kappa$ and $\varrho^*(dy)= \varrho (dy/ \kappa)$. 
\item[\rm(ii)] $R_h (\widetilde{\cF})$ is a {\rm GW($\xi_{v(h)}, c_{v(h)}; \mu_{v(h)}$)}-real forest. 
\item[\rm(iii)] The process $Z^\varrho$ is a {\rm CSBP}{\rm ($\psi, \varrho$)}.
\item[\rm(iv)] $ \lim_{h \rightarrow 0} \; \frac{1}{v(h) } Z^{(h)}_{a} (\widetilde{\cF})= Z^\varrho_a $, in probability, for all $a\in [0, \infty)$. 
\end{enumerate}
\end{theorem}
\begin{pf} By Lemma \ref{heraselim} $\rm(iii)$ and Corollary \ref{convcriter},
$\lim_{\lambda \rightarrow \infty} \widetilde{\cF}_\lambda = \widetilde{\cF}$ weakly in $\bT$ and 
$\rm(ii)$ holds true. Moreover, the last point of $\rm(i)$ is then a consequence of Theorem \ref{structuregrowth}.

  We next prove the joint convergence. To simplify notation, we set $Y^\lambda=(\frac{1}{\lambda}
 Z^{+}_a (\widetilde{\cF}_\lambda))_{a \in [0, \infty)}$. Recall from (\ref{convbranch}) that $Y^\lambda $ converges weakly in $\bD ([0, \infty), \bR)$ to a CSBP($\psi, \varrho$). 
Thus, the joint laws of the $(Y^\lambda, \widetilde{\cF}_\lambda)$ are tight in 
$\bD([0, \infty) ,  \bR)\!\times\! \bT$ as $\lambda \rightarrow \infty$.  
We want to prove that there is a unique limiting distribution by proving $\rm(iv)$ for a possible limit. To that end, let us assume that along a sequence $(\lambda_n)_{ n \in \bN}$ increasing to $\infty$, the following convergence holds in distribution in $\bD([0, \infty) ,  \bR)\times \bT$.
\begin{equation}
\label{subseqpaire}
\left( \;  Y^{\lambda_n}    \; ; \;  \widetilde{\cF}_{\lambda_n} \;  \right) \; \; 
\underset{{n \rightarrow \infty}}{\longrightarrow} \; \; \left(  \;  Y \; ; \; \widetilde{\cF}  \; \right). 
\end{equation}
By a slight abuse of notation, we assume that $Y$ and $\widetilde{\cF}$ are defined on $(\Omega, \cG, \bP)$. We fix $a\in [0, \infty)$ and $h\in (0, \infty)$. We know from Lemma \ref{heraselim} that $Z^{(h)}_a (\widetilde{\cF}_{\lambda_n})$ converges in law on $\bN$. Then the laws of $(Z^{(h)}_a (\widetilde{\cF}_{\lambda_n});   Y^{\lambda_n}    ; \widetilde{\cF}_{\lambda_n} )$ are tight on $\bN\times \bD([0, \infty) ,  \bR)\times \bT$. There is a increasing sequence of integers $(n(k))_{k\in \bN}$ such that the following convergence holds in distribution in $\bN \times \bD([0, \infty) ,  \bR)\times \bT$.
\begin{equation}
\label{subsubseq}
\left( \; Z^{(h)}_a (\widetilde{\cF}_{\lambda_{ n(k) }})  \; ; \;  Y^{\lambda_{n(k)}}    \; ; \;  \widetilde{\cF}_{\lambda_{n(k)}} \;  \right) \; \; 
\underset{{k \rightarrow \infty}}{\longrightarrow} \; \; \left(  \; X\; ;  \;  Y^\prime  \; ; \; \widetilde{\cF}^\prime  \; \right). 
\end{equation}
Again, to simplify notation, we assume that $X$, $Y^\prime $ and $\widetilde{\cF}^\prime $ 
are defined on $(\Omega, \cG, \bP)$. Clearly $(Y^\prime, \widetilde{\cF}^\prime)$  and $(Y, \widetilde{\cF})$ have the same law. The $\dgh$-continuity of $R_h$ implies that 
$\lim_{k\rightarrow \infty}(R_h (\widetilde{\cF}_{\lambda_{n(k)}} ) ;  \widetilde{\cF}_{\lambda_{n(k)}} ) =
 (R_h (\widetilde{\cF}^\prime) ,  \widetilde{\cF}^\prime )$ weakly on $\bT^2$. Moreover, $\widetilde{\cF}^\prime$ is a $(\psi, \varrho)$-Lévy forest and by $\rm(ii)$, $R_h (\widetilde{\cF}^\prime)$ has the same law as $\widetilde{\cF}_{v(h)}$. Lemma \ref{heraselim} asserts that, weakly on $\bN\times \bT$,  
 $\lim_{k\rightarrow \infty}(Z^{(h)}_a(\widetilde{\cF}_{\lambda_{n(k)}} ) ; R_h ( \widetilde{\cF}_{\lambda_{n(k)}} ))=
 (Z^+_a (\widetilde{\cF}_{v(h)}); \widetilde{\cF}_{v(h)} )$, which has the same law as 
 $(Z^{+}_a (R_h (\widetilde{\cF}^\prime)) ; R_h (\widetilde{\cF}^\prime))$. 
Thus $(X;  R_h (\widetilde{\cF}^\prime ))$ has the same law as $(Z^{+}_a(R_h(\widetilde{\cF}^\prime))  ;  R_h (\widetilde{\cF}^\prime ))$, which implies that 
\begin{equation} 
\label{zhidenti}
\textrm{$\bP$-a.s.} \quad X= Z^+_a (R_h(\widetilde{\cF}^\prime))= Z^{(h)}_a(\widetilde{\cF}^\prime) \; ,
\end{equation}
the last equality being the definition of $Z^{(h)}_a$. We next recall from (\ref{halphah}) that $Q_{\xi_\lambda , c_\lambda} (\Gamma > h)= \lambda^{-1}u(h, \lambda)$. We recall from (\ref{branchfor}) that conditionally given $Z^+_a (\widetilde{\cF}_\lambda)$, $\abv (a,  \widetilde{\cF}_\lambda)$ is distributed according to $Q_{\xi_\lambda , c_\lambda}^{\circledast Z^+_a (\widetilde{\cF}_\lambda)}$. Since $Z^{(h)}_a (\widetilde{\cF}_\lambda)= \langle \cM_0 (\abv (a,  \widetilde{\cF}_\lambda)) , \un_{\{ \Gamma >h \}} \rangle $, we know that conditionally given $Z^+_a (\widetilde{\cF}_\lambda)$, the variable $Z^{(h)}_a (\widetilde{\cF}_\lambda)$ has a binomial law with parameters $  \lambda^{-1}u(h, \lambda)$ and $Z^+_a (\widetilde{\cF}_\lambda)$. Thus, for all $\theta \in [0, \infty)$ and all continuous bounded functions $f$, we easily get 
$$ \bE \left[ e^{-\theta Z^{(h)}_a (\widetilde{\cF}_{\lambda_{n(k)}} )} f( \lambda_{n(k)}^{-1}Z^{+}_a (\widetilde{\cF}_{\lambda_{n(k)}})) \right] = \bE \left[ \left( 1-
\frac{{u(h, \lambda_{n(k)})}}{{\lambda_{n(k)}}}  \left(1\!-\!e^{-\theta} \right) \right)^{Z^{+}_a (\widetilde{\cF}_{\lambda_{n(k)}})}  f( \lambda_{n(k)}^{-1}Z^{+}_a (\widetilde{\cF}_{\lambda_{n(k)}})) \right]. $$
We then pass to the limit as $k\to\infty$ to get 
$$ \bE \left[ \exp ( - \theta X)
\,  f(Y^\prime_a)  \right] =\bE \left[ \exp \left(-v(h) 
 Y^\prime_a \left(1- e^{-\theta} \right)\right) \,  f(Y^\prime_a)  \right].$$
This entails that conditionally given $Y^\prime_a$, $X$ is a Poisson random variable with parameter 
$v(h)Y^\prime_a$. By (\ref{zhidenti}), for all $K, \varepsilon \in (0, \infty)$, we get 
$$\bP \left( \, \left| \frac{1}{{v(h)}}Z^{(h)}_{a} (\widetilde{\cF})- Y_a       \right| >   \varepsilon \right)
=  \bP \left( \, \left| \frac{1}{{v(h)}}Z^{(h)}_{a} (\widetilde{\cF}^\prime)- Y^\prime_a       \right| >   \varepsilon \right)
\leq   \frac{K }{\varepsilon^2 v(h) } + \bP (  Y^\prime_a  >K) . $$
This entails that for all $a \in [0, \infty)$, $ \lim_{h \rightarrow 0}\frac{1}{v(h) } Z^{(h)}_{a} (\widetilde{\cF})=Y_a$ in probability. Thus, we get the uniqueness of the limit for the joint laws, which completes the proof of the theorem.  
 \end{pf}

\paragraph{Almost sure convergence of growth processes.} We turn now to the main result of this section, which asserts almost sure convergence for growth processes and their branching processes. We only consider the most interesting case when the limiting tree is not a GW-tree. 
\begin{theorem}
\label{asconvgrowth} Let $(\Omega, \cG, \bP)$ be a probability space and let $\widetilde{\cF}_\lambda\colon  \Omega \rightarrow \bT$,  
$\lambda \in [0, \infty)$, be a growth process such that $\bP ( Z^+_0 (\widetilde{\cF}_{\lambda}) \rightarrow \infty) >0$. Let $(\psi , \varrho, \beta)$ be the triplet governing the growth process as specified in Theorem 
\ref{structuregrowth} ${\rm (II)}$. We furthermore assume $\int^\infty dr/ \psi (r)< \infty$. 
Then there exists a random CLCR real tree $\widetilde{\cF}\colon \Omega \rightarrow \bT$ such that
\begin{equation}
\label{aslimfo}
 \textrm{$\bP$-a.s.} \qquad \lim_{\lambda \rightarrow \infty} \dgh \left(\widetilde{\cF}_{\lambda} , \widetilde{\cF} \right) = 0 \; .
\end{equation} 
The random tree $\widetilde{\cF}$ is a $(\psi, \varrho)$-Lévy forest as defined in 
Theorem \ref{growthconvergence} and there exists a cadlag CSBP{\rm ($\psi, \varrho$)} denoted by $(Z^\varrho_a)_{a\in [0, \infty)}$ such that for all $a\in [0, \infty)$ 
\begin{equation}
\label{asapproxZ}  
\textrm{$\bP$-a.s.} \qquad  Z^\varrho_a = \lim_{\lambda \rightarrow \infty} \frac{1}{{\beta(\lambda)}} 
Z^{+}_a ( \widetilde{\cF}_{\lambda} )  =
\lim_{h\rightarrow 0+} \frac{1}{{v(h)}} Z^{(h)}_a ( \widetilde{\cF}) \; .
\end{equation}
Moreover, the same limits hold in $L^1$ if $\psi^{\prime}(0+)$ and $\int_{[0, \infty)} y \varrho (dy) $ are both finite. 
\end{theorem}
\begin{pf} We keep the notation of Theorem \ref{structuregrowth}: $\widetilde{\cF}_{\lambda}$ is a 
GW($\xi_{\beta(\lambda)}, c_{\beta (\lambda)}; \mu_{\beta( \lambda)}$)-real forest, where for all $\lambda \in [q, \infty)$, $(\xi_\lambda , c_\lambda, \mu_\lambda)$ is given by (\ref{xilambdef}). Lemma \ref{heraselim} asserts that for all $a\in [0, \infty)$ and all $h\in (0, \infty)$, $Z^{(h)}_a(\widetilde{\cF}_{\lambda})$ converges in law as $\lambda$ goes to $\infty$. By Remark \ref{reggrow}, 
Theorem \ref{tightness2} applies and there exists $\widetilde{\cF}$ such that (\ref{aslimfo}) holds. 
By Theorem \ref{growthconvergence}, $\widetilde{\cF}$ is a $(\psi, \varrho)$-Lévy forest and 
there exists a cadlag CSBP{\rm ($\psi, \varrho$)} denoted by $(Z^\varrho_a)_{a\in [0, \infty)}$ such that for all $a\in [0, \infty)$, $\lim_{h\rightarrow 0+} \frac{1}{{v(h)}} Z^{(h)}_a (\widetilde{\cF}_{\lambda})= Z^\varrho_a$ in probability 
and such that the following convergence holds weakly on $\bD ([0, \infty) , \bR) \times \bT$ 
\begin{equation}
\label{cconvjointe}\textstyle
\left( \; \left( \frac{1}{{\beta (\lambda)}}Z^+_{a} (\widetilde{\cF}_\lambda) \, \right)_{ a \in [0, \infty )}   \; ; \;  \widetilde{\cF}_\lambda \;  \right) \; \; \underset{{\lambda \rightarrow \infty}}{\longrightarrow } \; \; \left(  \; 
( Z^\varrho_{a} )_{ a\in [0, \infty)} \; ; \; \widetilde{\cF } \; \right). 
\end{equation}
  We fix $a\in [0, \infty)$ and we turn now to the a.s.~convergence of the branching processes. 
We use the following martingale argument. 
For all $\lambda \in [0, \infty)$, we denote by $\cH_\lambda$ the sigma-field generated by the variables $Z^+_a (\widetilde{\cF}_{\lambda^{\prime}})$, 
$\lambda^{\prime } \in [\lambda, \infty)$. {\it We then claim that $(\frac{1}{{\beta( \lambda)}}Z^+_{a} (\widetilde{\cF}_\lambda) )_{\lambda \in [0, \infty)}$ is a nonnegative backward martingale with respect to 
$(\cH_{\lambda})_{\lambda \in  [0, \infty)}$}. 

\smallskip

\noi
\textsl{Proof of the claim.} 
We first fix $\lambda, \lambda^\prime \in  [0, \infty)$,  such that $\lambda^\prime >\lambda$, and we denote by $\widetilde{\cT}_{\lambda^\prime}$ a 
GW($\xi_{\beta(\lambda^\prime) } , c_{\beta(\lambda^\prime) }  $)-real tree. We use the notation of Definition \ref{growrealdef} and we denote by $A_{\lambda, \lambda^\prime} \subseteq \bT$ the family of hereditary properties under which the growth process is consistent. We then set 
$$ \alpha_{\lambda, \lambda^\prime} := 
Q_{\xi_{\beta(\lambda^\prime) } , c_{\beta(\lambda^\prime) }  } \left( \bT \backslash A^-_{\lambda, \lambda^\prime}\right)=Q_{\xi_{\beta(\lambda^\prime)},c_{\beta(\lambda^\prime)}}\left(\bT\setminus A_{\lambda,\lambda^\prime}^-\right) \; .$$
The $A_{\lambda, \lambda^\prime}$-reduced offspring distribution is 
$\xi^{(\alpha_{\lambda, \lambda^\prime})}_{\beta (\lambda^\prime)}= \xi_{\beta (\lambda)}$, and a brief computation involving (\ref{reduclaw}) and (\ref{xilambdef}) entails 
\begin{equation}
\label{alphcompp}
\alpha_{\lambda, \lambda^\prime} = 1-\frac{{\beta (\lambda)}}{{\beta (\lambda^\prime)}} \; .
\end{equation}  
Let $G\colon \bT \rightarrow [0, \infty)$ be measurable and let $n^\prime \geq n$. We next compute 
$$ B:= \bE \left[ G\left( \abv (a, \widetilde{\cF}_{\lambda})  \right) \, 
 \un_{ \{  Z^+_a (\widetilde{\cF}_{\lambda^\prime } )  = n^\prime  \, ; \,  
Z^+_a (\widetilde{\cF}_{\lambda}) = n  \}} \right]  \; .$$
We recall from (\ref{branchfor}) that conditionally given the event $\{ Z^+_a (\widetilde{\cF}_{\lambda^\prime } )  
= n^\prime \}$, $\abv (a, \widetilde{\cF}_{\lambda^\prime } )$ is distributed as a forest of 
$n^\prime $ independent 
GW($\xi_{\beta (\lambda^\prime)} , c_{\beta (\lambda^\prime)}$)-real trees.  By (\ref{AbovAred}) 
in Lemma \ref{Amoinsprop}, $R_{A_{\lambda, \lambda^\prime}} 
(\abv (a, \widetilde{\cF}_{\lambda^\prime } ))= \abv (a, \widetilde{\cF}_{\lambda} )$. 
Then, by Theorem \ref{realforreduc}, conditionally given the event 
$\{ Z^+_a (\widetilde{\cF}_{\lambda^\prime } )  = n^\prime \}$, 
$ \abv (a, \widetilde{\cF}_{\lambda} )$ is a 
GW($\xi_{\beta (\lambda)}, c_{\beta (\lambda)}; \mu$)-real forest where $\mu$ stands for the binomial distribution with parameters $n^\prime$ and $1-\alpha_{\lambda, \lambda^\prime}$. This implies 
the following:
\begin{eqnarray*}
B &\!\!\!=\!\!\!&  \bP \left( Z^+_a (\widetilde{\cF}_{\lambda^\prime } )  = n^\prime   \right) P_{\xi_{\beta (\lambda)}, c_{\beta (\lambda)}, \mu } \left[ \un_{\{ Z^+_0 = n\}} G \, \right] \\
& \!\!\!=\!\!\!&  \bP \left( Z^+_a (\widetilde{\cF}_{\lambda^\prime } )  = n^\prime   \right) 
\binom{ n^\prime}{n} (1-\alpha^{}_{\lambda, \lambda^\prime})^{n} 
\alpha_{\lambda, \lambda^\prime}^{n^\prime -n} \, 
Q_{\xi_{\beta (\lambda)}, c_{\beta (\lambda)} }^{\circledast n} \left[ G\right] \\
& \!\!\!=\!\!\!& \frac{\bP \left( Z^+_a (\widetilde{\cF}_{\lambda^\prime } )  = n^\prime 
  \right) }{\bP \left( Z^+_a (\widetilde{\cF}_{\lambda } )  = n \,  \right) }
\binom{n^\prime}{n} (1-\alpha^{}_{\lambda, \lambda^\prime})^{n} 
\alpha_{\lambda, \lambda^\prime}^{n^\prime -n}
\, \bE \left[ G\left( \abv (a, \widetilde{\cF}_{\lambda})  \right) \, 
 \un_{\{  Z^+_a (\widetilde{\cF}_{\lambda}) = n  \}} \right] \; . 
\end{eqnarray*}

We next fix the real numbers $\lambda_p > \cdots > \lambda_0 \geq 0$, 
and the integers $n_p \geq \cdots \geq n_1 \geq 0 $, and for all $n_0 \in \{ 0, \ldots , n_1 \}$, we set 
$$u (n_0)=  \bP \left(  Z^+_a (\widetilde{\cF}_{\lambda_{p}}) \! =\!  n_{p} \,  ; \,   \ldots  \, ; \,     Z^+_a (\widetilde{\cF}_{\lambda_1 }) = n_1 \, ; \,  
Z^+_a (\widetilde{\cF}_{\lambda_0 }) \!= \! n_0 \right) . $$
We now apply the previous 
computation to $  G( \abv (a, \widetilde{\cF}_{\lambda_{p-1}}))= \un_{\{  Z^+_a (\widetilde{\cF}_{\lambda_{p-1}})  = n_{p-1} \, ; \, \ldots \, ; \, 
Z^+_a (\widetilde{\cF}_{\lambda_0 }) = n_0  \}} $, $ \lambda^\prime=\lambda_p $, $n^\prime= n_p$, $\lambda= \lambda_{p-1}$ and $ n= n_{p-1}$  
to get 
\begin{eqnarray*}
 u (n_0)   = \;   
\frac{{\bP \left( Z^+_a (\widetilde{\cF}_{\lambda_p} )  =   n_p \right)} }{{\bP \left( 
Z^+_a (\widetilde{\cF}_{\lambda_{p\!-\!1} } )  =  n_{p-1}   \right) }} & & \!\!\!\!\!\!\!\!\!\!\!  
\binom{n_p}{n_{p\!-\!1}} (1-\alpha^{}_{\lambda_{p\!-\!1}, \lambda_{p}})^{ n_{p \!-\! 1}} 
\alpha_{\lambda_{p\!-\!1}, \lambda_p }^{n_p  - n_{p \!-\! 1}} \;  \times \\
& & \qquad \qquad \; \bP \left( Z^+_a (\widetilde{\cF}_{\lambda_{p-1}})  =   n_{p-1} \, ; \,   \ldots  \, ; \,   
Z^+_a (\widetilde{\cF}_{\lambda_0 }) =  n_0 \right) \; .
\end{eqnarray*}
This entails 
$ u (n_0) = \bP ( Z^+_a (\widetilde{\cF}_{\lambda_p} )  =   n_p) 
\prod_{1\leq k\leq p} \binom{n_k}{n_{k\!-\!1}} (1-\alpha^{}_{\lambda_{k\!-\!1}, \lambda_{k}})^{n_{k \!-\! 1}} 
\alpha_{\lambda_{k\!-\!1}, \lambda_k }^{n_k  - n_{k \!-\! 1}} $. Thus we get 
\begin{equation}
\label{caution}
 u(n_0)= \binom{n_1}{n_{0}} (1-\alpha^{}_{\lambda_{0}, \lambda_{1}})^{n_0} 
\alpha_{\lambda_{0}, \lambda_1 }^{n_1  - n_{0}} \bP \left(  Z^+_a (\widetilde{\cF}_{\lambda_{p}}) \! =\!  n_{p} \,  ; \,   \ldots  \, ; \,     Z^+_a (\widetilde{\cF}_{\lambda_1 }) = n_1   \right) \, . 
\end{equation}
We now fix $\lambda^\prime  >\lambda$ and we denote by 
$\cP_{\lambda^\prime}$ the set of events of the form 
$$\{ Z^+_a (\widetilde{\cF}_{\lambda_{p}}) \! \geq \!  n_{p} \,  ; \,   \ldots  \, ; \,     Z^+_a (\widetilde{\cF}_{\lambda_1 }) \geq  n_1\}, \; \lambda_1, \ldots, \lambda_p \in [\lambda^\prime, \infty)\, , \; 
n_1, \ldots , n_p\in \bN\, , \; p\in \bN^*\; .$$ 
Clearly $\cP_{\lambda^\prime}$ contains $\Omega$, it is stable under intersection and it generates $\cH_{\lambda_1}$. Moreover, (\ref{caution}) 
with $\lambda= \lambda_0$ and $\lambda_1= \lambda^\prime$,  easily entails that for all $n \in \bN$, and all $A \in \cP_{\lambda^\prime} $, 
$$ \bE \left[ \un_{\{ Z^+_a (\widetilde{\cF}_\lambda) = n\}} \un_A \right]= 
\bE \left[ \binom{Z^+_a (\widetilde{\cF}_{\lambda^\prime})}{n} (1-\alpha^{}_{\lambda, \lambda^{\prime}})^{n} 
\alpha_{\lambda, \lambda^\prime }^{ Z^+_a (\widetilde{\cF}_{\lambda^\prime}) - n}  \un_A \right]. $$
A monotone class argument entails that the same equality holds for all $A\in \cH_{\lambda^\prime}$. Thus, for all $A \in \cH_{\lambda^\prime}$
\begin{eqnarray*}
 \bE \left[  Z^+_a (\widetilde{\cF}_\lambda)\un_A \right] &= & \sum_{n\geq 0} \bE \left[ n \binom{Z^+_a (\widetilde{\cF}_{\lambda^\prime})}{n} (1-\alpha^{}_{\lambda, \lambda^{\prime}})^{n} 
\alpha_{\lambda, \lambda^\prime }^{ Z^+_a (\widetilde{\cF}_{\lambda^\prime}) - n}  \un_A \right] \\
& =& \bE \left[  (1-\alpha_{\lambda, \lambda^\prime}) 
Z^+_a (\widetilde{\cF}_{\lambda^\prime})  \un_A  \right]   
= \frac{{\beta (\lambda)}}{{\beta (\lambda^\prime)}}  \bE \left[ Z^+_a (\widetilde{\cF}_{\lambda^\prime})  \un_A  \right] ,   
\end{eqnarray*}
by (\ref{alphcompp}), which immediately implies the claim. \cq

\medskip

The theorem of almost sure convergence of nonnegative backward martingale implies that 
for every sequence $(\lambda_n)_{n \in \bN}$ that increases to $\infty$, 
$\bP$-a.s.\ $\lim_{n \rightarrow \infty} \frac{1}{{\beta(\lambda_n)}} Z^{+}_a ( \widetilde{\cF}_{\lambda_n} ) $ exists. Then, the joint convergence (\ref{cconvjointe}) also implies that the limiting r.v.~is necessarily 
$\bP$-a.s.~equal to $Z^\varrho_a$.

  By (\ref{expecZZ}), if both $\psi^{\prime}(0+)$ and $\int_{[0, \infty)} y\varrho (dy) $ are finite, then $Z^\varrho_a$ is integrable. Standard results on backward martingales then entail that for all $\lambda \in [0, \infty)$, $\frac{1}{\beta (\lambda)} Z^+_a (\widetilde{\cF}_\lambda)$ is integrable and that 
$$ \lim_{\lambda \rightarrow \infty} \bE \left[  \, \Big| \frac{_1}{^{\beta (\lambda)}} Z^+_a (\widetilde{\cF}_\lambda) -Z^\varrho_a \Big| \,  \right] = 0 \; .$$

Since $\beta$ may have jumps, additional arguments are required to get the first equality in (\ref{asapproxZ}), which is proved so far only along a subsequence. 
To simplify notation, we first 
set $X_\lambda=  \frac{1}{{\beta(\lambda)}} Z^{+}_a ( \widetilde{\cF}_{\lambda} )$. 
Then, we fix $\varepsilon \in (0, \infty)$ and $\lambda^* \in (0, \infty) $ 
such that $\beta(\lambda^*) >0$. 
We recursively define an increasing sequence $(\lambda_n)_{n\in \bN}$ 
that tends to $\infty$ by fixing $\lambda_0 >\lambda^*$ and by setting  
$\lambda_{n+1}= \inf \{ \lambda>\lambda_n \colon  \log \beta (\lambda) > \varepsilon + \log 
\beta (\lambda_n+)\} $. 
Then observe that for all $\lambda, \lambda^\prime \in (\lambda_n , \lambda_{n+1})$, $e^{-\varepsilon} \leq \beta (\lambda^\prime) / \beta (\lambda) \leq  e^{\varepsilon}$.  
Since the growth process is $\preceq$-non-decreasing,  
(\ref{herasegrow}) entails that $\bP$-a.s.~for all $n \in \bN$, for all $\lambda,  \lambda^{ \prime},  \lambda^{\prime \prime} \in (\lambda_n , \lambda_{n+1})$ such that $ \lambda^{\prime \prime}\leq  \lambda\leq  \lambda^{\prime}$, 
$$ e^{-\varepsilon} X_{\lambda^{\prime \prime} } \leq  X_{\lambda} \leq e^{\varepsilon} X_{\lambda^{\prime } }. $$
For all $\lambda \in (0, \infty)$ we also define the sigma-field 
$\cH_{\lambda-}= \bigcap_{\lambda^\prime <\lambda } \cH_{\lambda^\prime}$ and we also denote by 
$\cH_{\lambda+}$ the sigma-field generated by $\bigcup_{\lambda^\prime > \lambda} \cH_\lambda$. 
We next fix a sequence $(h_p)_{p \in \bN}$ that decreases to $0$. 
By standard arguments for nonnegative martingales and for nonnegative backward martingales, we get the following: $\bP$-a.s.~for all $n \in \bN$,  
$$ \bE [X_{\lambda^*} \, | \, \cH_{\lambda_n\!-}]= \lim_{p \rightarrow \infty} 
X_{\lambda_n-h_p}=:X_{\lambda_n\! -} \quad \textrm{and} \quad 
\bE [X_{\lambda^*} \, | \, \cH_{\lambda_n\!+}]= \lim_{p \rightarrow \infty} 
X_{\lambda_n+h_p}=:X_{\lambda_n \!+} \; .$$ 
Now observe that $(X_{\lambda_n-})_{n \in \bN}$ and $(X_{\lambda_n+})_{n \in \bN}$ are 
positive backward martingales with respect to the backward filtrations 
$(\cH_{\lambda_n-})_{n \in \bN}$ and $(\cH_{\lambda_n+})_{n \in \bN}$. So, they both 
converge $\bP$-a.s. By Theorem \ref{tightness2}, $\widetilde{\cF}_{\lambda-h_p}$ converges $\bP$-a.s., as $p\rightarrow \infty$. We denote the limit by $\widetilde{\cF}_{\lambda-}$.  
By Lemma \ref{poubelle}, observe that $(X_{\lambda-h_p}, \widetilde{\cF}_{\lambda-h_p})$ converges in distribution to a GW($\xi_{\beta(\lambda-)} ,c_{\beta(\lambda-)}, \mu_{\beta(\lambda-)} $)-real forest as $p\rightarrow \infty$. Moreover, 
Theorem \ref{growthconvergence} applies to such laws, which implies a joint weak convergence that is similar to (\ref{cconvjointe}) with the left-limit process instead of the normal one. 
This easily entails that $\bP$-a.s.~$\lim_{n\rightarrow \infty} X_{\lambda_n-}= Z^\varrho_a$. A similar argument also implies that $\bP$-a.s.~$\lim_{n\rightarrow \infty} X_{\lambda_n+}= Z^\varrho_a$. 

This proves that for all $\varepsilon \in (0, \infty)$, there exists an event $\Omega_\varepsilon \in \cF$ of probability one on which the following holds true: $X_{\lambda_n}$, $X_{\lambda_n+}$ and 
$X_{\lambda_n-}$ tend to $Z^\varrho_a$ as $n\rightarrow \infty$, and 
$$\forall n \in \bN, \; \forall \lambda \in [\lambda_n , \lambda_{n\!+\!1}], \quad  e^{-\varepsilon} \min (X_{\lambda_n+} ,  X_{\lambda_n} , X_{\lambda_{n\!+\!1}} )
 \leq  X_{\lambda} \leq e^{\varepsilon} \max (X_{\lambda_{n\!+\!1}-} ,  X_{\lambda_n} , 
 X_{\lambda_{n+1}} ). 
 $$
Thus, on  $\Omega_\varepsilon$, 
$ e^{-\varepsilon}Z^\varrho_a \leq \liminf_{\lambda \rightarrow \infty} X_\lambda \leq 
\limsup_{\lambda \rightarrow \infty} X_\lambda\leq e^{\varepsilon}Z^\varrho_a$. This easily entails the first equality in (\ref{asapproxZ}).

It remains to prove the assertions concerning $Z^{(h)}_a (\widetilde{\cF})$. To that end, recall from 
(\ref{vvvdefini}) the definition of the function $v$ that is continuous 
decreasing from $(0, \infty)$ to $(q, \infty)$. We denote by $v^{-1}$ its inverse and for all 
$\lambda \in [0, \infty)$, we set $\widetilde{\cF}_{\lambda}^\prime= 
R_{v^{-1} (q+1+\lambda)} (\widetilde{\cF})$. It clearly defines a growth process that $\dgh$-converges to $\widetilde{\cF}$.
This process is governed by a triplet of the form 
$(\psi, \varrho, \beta^\prime)$ with $\beta^\prime (\lambda)= q+1+ \lambda$. 
Now observe that for  
all $a \in [0, \infty)$ and all $h$ such that $v^{-1} (1\!+\! q) > h >0$,    
$$
\frac{1}{{v(h)}} Z^{(h)}_a (\widetilde{\cF})=  \frac{1}{{\beta^\prime (v(h)\! -\! 1\!-\! q)}}
Z^{+}_a \left( \widetilde{\cF}^\prime_{v(h) - 1- q}\right) \; .$$
We now apply the first  equality in (\ref{asapproxZ}) to that specific growth process to obtain the second equality in (\ref{asapproxZ}). \end{pf}

As an application of the previous results on growth processes, we first state a simple characterisation of Lévy forests that is used to derive limit theorems of GW-forests to Lévy forests. 
\begin{theorem} 
\label{Levyfochar}
Let $(\Omega, \cG, \bP)$ be a probability space. Let $\widetilde{\cF} \colon  \Omega \rightarrow \bT$ be a random CLCR real tree.  
We assume that $\bP (\widetilde{\cF} \neq \pnt) >0$ and 
that for all sufficiently small $h \in (0, \infty)$, $R_h (\widetilde{\cF})$ is a GW-real forest.  
Then, the following holds true. 
\begin{enumerate}
\item[{\rm (a)}] Either $\bP (Z^+_0 (\widetilde{\cF} )< \infty)= 1$ and $\widetilde{\cF} $ is a GW-real forest . 
\item[{\rm (b)}] Or $\bP (Z^+_0 (\widetilde{\cF} )= \infty) >0$ and 
$\widetilde{\cF} $ is a L\'evy forest as in Theorem \ref{growthconvergence}. 
\end{enumerate}
\end{theorem}
\begin{pf} First note that $\bP (R_h (\widetilde{\cF}) \neq \pnt ) >0$, for all sufficiently small 
$h$.
Let $h_0$ be such that for all $h\in (0, h_0)$, $R_h (\widetilde{\cF})$ is a GW-real forest such that 
$\bP (R_h (\widetilde{\cF})\neq \pnt) >0$. Then, for all 
$\lambda \in [0, \infty)$, we set $\widetilde{\cF}_\lambda = R_{h_0 e^{-\lambda}} (\widetilde{\cF})$. Clearly, $(\widetilde{\cF}_\lambda)_{\lambda \in [0, \infty)}$ is a growth process. We then apply 
Theorem \ref{structuregrowth}. Namely,  
since $Z^+_0 (\widetilde{\cF}_\lambda) \rightarrow Z^+_0 (\widetilde{\cF})$, 
either $\bP (Z^+_0 (\widetilde{\cF}) < \infty  )= 1$ and the growth process is as in Theorem \ref{structuregrowth} (I): by Lemma \ref{poubelle}, $\widetilde{\cF}_\lambda$ converges in law 
to a GW-forest as $\lambda \rightarrow \infty$, which implies that $\widetilde{\cF}$ is a GW-real forest.  
Or $\bP (Z^+_0 (\widetilde{\cF}) = \infty  ) >0$, 
and the growth process is as in Theorem \ref{structuregrowth} (Il): the growth process is then governed by a triplet $(\psi, \varrho, \beta)$. Since the laws of the $\widetilde{\cF}_\lambda$ are tight in $\bT$ as $\lambda \rightarrow \infty$, Lemma \ref{heraselim} $\rm(ii)$ implies that $\int^\infty dr/ \psi (r) < \infty$, and Theorem \ref{asconvgrowth} implies that $\widetilde{\cF}$ is a Lévy forest.  
\end{pf}

Theorem \ref{Levyfochar} allows us to strengthen Lemma \ref{regenene}, as follows. Recall from Section \ref{GWrealtreesec} the definition of the sigma-field $\cB_{a+}$. 

\begin{theorem}
\label{vrairege} Let $Q$ be a Borel probability measure on $\bT$. We first assume that 
$Q(0 < D< \infty)>0 $. We also assume that for all $a\in[0,\infty)$, $Q(Z_a^+<\infty)=1$. Then, the following assertions are equivalent.
\begin{enumerate}
\item[$\rm(i)$] For every $a \in [0, \infty)$, the conditional distribution given $Z^+_a$ of $\abv (a, \cdot)$ under $Q$ is $Q^{\circledast Z^+_{a} } $.   
\item[$\rm(ii)$] For every $a \in [0, \infty)$, the conditional distribution given $\cB_{a+}$ of $\abv (a, \cdot)$ under $Q$ is $Q^{\circledast Z^+_{a} } $. 
\item[$\rm(iii)$] $Q$ is the distribution of a {\rm GW($\xi, c$)}-real tree, for a certain $c \in (0, \infty)$ and a certain proper conservative offspring distribution $\xi$. 
\end{enumerate}
\end{theorem}
\begin{pf} If we assume (iii), then $a\mapsto Z_a^+$ is cadlag $Q$-a.s., hence Lemma \ref{regenene} applies and we get (ii) and (i). It remains to 
prove $\rm(i) \Rightarrow \rm(iii)$ without assuming, as in Lemma \ref{regenene}, that $a\mapsto Z_a^+$ is cadlag $Q$-a.s.
To this end, first note that (i) implies $Q(Z_0^+=1)=1$, as in the proof of
Lemma \ref{regenene} in Appendix \ref{regenenepf}. Next, observe that $\Gamma\ge D$ on $\bT\setminus\{\pnt\}$, and recall from Lemma \ref{measDthetak} that $\lim_{h\rightarrow 0}D\circ R_h=D$. Thus, there exists $h_0\in(0,\infty)$ such that for all $h\in(0,h_0)$, $Q(0<D\circ R_h<\infty)>0$ and $q_h:=Q(\Gamma>h)>0$. Then it makes sense to define $Q_h$ as the law of $R_h$ under $Q(\,\cdot\,|\Gamma>h)$ and we have proved that $Q_h(0<D<\infty)>0$. Next, observe that $a\mapsto Z_a^+\circ R_h=Z_{a+}^{(h)}$ is cadlag $Q$-a.s., which implies that $a\mapsto Z_a^+$ is cadlag $Q_h$-a.s. 

Let us show that $Q_h$ satisfies (i). Let us fix $a\in[0,\infty)$, $n\in\bN\setminus\{0\}$ and a measurable function $G\colon\bT\rightarrow[0,\infty)$. The property (i) for $Q$, (\ref{AbovAred}) and arguments similar to those used in the proof of Theorem \ref{GWrealreduc} imply the following
\beq Q_h\left[\un_{\{ Z^+_a = n \} } G (\abv (a, \cdot) )\right]&\!\!\!\!=\!\!\!&q_h^{-1}Q\left[\un_{\{Z_a^+\circ R_h=n\}}G(\abv(a,R_h))\right]\\
&\!\!\!\!=\!\!\!&q_h^{-1}\sum_{N\ge n}Q\left[\un_{\{Z_a^+=N;\, Z_a^+\circ R_h=n\}}G\left(R_h(\abv(a,\cdot))\right)\right]\\
&\!\!\!\!=\!\!\!&q_h^{-1}\sum_{N\ge n}Q\left[\un_{\{Z_a^+=N\}}Q^{\circledast N}\left[\un_{\{Z_a^+\circ R_h=n\}}G(R_h)\right]\right]\\
&\!\!\!\!=\!\!\!&q_h^{-1}\sum_{N\ge n}Q(Z_a^+=N)\binom{N}{n}(1-q_h)^{N-n}q_h^nQ_h^{\circledast n} \left[ G \right]=Q_h(Z_a^+=n)Q_h^{\circledast n}\left[G\right].
\eeq
This implies that $Q_h$ satisfies (i). Then Lemma \ref{regenene} implies that $Q_h$ is the law of a GW-real tree. Thus, for any $h\in(0,h_0)$, $R_h$ under $Q$ is the law of a GW-real forest and since $Q(Z_0^+=1)=1$, Theorem \ref{Levyfochar} entails that $Q$ is the law of a GW-real tree, which completes the proof.
\end{pf}

\subsection{Invariance principles for GW-trees.}
\label{Invprincsec}
In this section we apply the tightness results of Section \ref{convcritsec} and the previous results on L\'evy forests to obtain limit theorems for discrete GW-real trees. We consider two asymptotic regimes: we first discuss convergence in distribution to GW-real trees and we next discuss convergence to L\'evy forests. 

\medskip

\noi{\bf Notation}. Unless the contrary is explicitly mentioned, all the random variables are defined on the same probability space $(\Omega, \cG, \bP)$. 
For all $p\in \bN$, 
\begin{enumerate}
\item[--] $\xi^p= (\xi^p (k))_{k \in \bN}$ is an offspring distribution such that $\xi^p (1) < 1$,  
\item[--] $\nu^p= (\nu^p (k))_{k \in \bZ}$ is given by $\nu^p(k)= \xi^p(k+1)$ if $k\geq -1$ 
and $\nu^p(k)= 0$ if $k<-1$,   
\item[--] $\mu^p = (\mu^p(k))_{k\in \bN}$ is a probability distribution on $\bN$ such that $\mu^p (0) < 1$, 
\item[--] $(Y^{p}_k )_{k \in \bN}$ is a discrete-time GW-Markov chain with initial distribution $\mu^p$ and offspring distribution $\xi^p$. 
\end{enumerate}

\noi 
{\it Discrete GW-real trees and forests.} We modify the lifetime part of Definition \ref{GWrealtreedef} and say that a
  $\bT$-valued random variable $\widetilde{\cT}_p$ is called a discrete GW($\xi^p$)-real tree if its distribution $Q$ satisfies
  $$Q\left[  \un_{\{{\bf k} = n \} } G (\vartheta) g( D)   \right]= \xi^p (n) Q^{\circledast n} [G] g(1)\; .$$
This is just a way to view graph trees as metric spaces by joining the vertices by intervals of unit length. 
Now let $\widetilde{\cT}_p (n)$, $n \geq 1$, be independent copies of $\widetilde{\cT}_p $. Let 
$N_p $ be an $\bN$-valued variable with law $\mu^p$, independent of 
$(\widetilde{\cT}_p (n))_{n\geq 1}$. We then set 
$$ \widetilde{\cF}_p:= \circledast_{1\leq n \leq N_p} \widetilde{\cT}_p (n) \; , $$
with the convention that $ \widetilde{\cF}_p= \pnt$ if $N_p= 0$. Then, $ \widetilde{\cF}_p$ is a random forest of $N_p$ independent discrete GW($\xi^p$)-real trees.

\smallskip

\noi 
{\it Tree-scaling.} 
The main purpose of this section is to obtain convergence in law of $\widetilde{\cF}_p $ when suitably rescaled. More precisely, let $\widetilde{T} \in \bT$ and $c \in (0, \infty)$; let $(T, d, \rho)$ be a representative of $\widetilde{T}$. Then, $(T, cd, \rho)$ is a CLCR real tree and its pointed isometry class only depends on $\widetilde{T}$ and $c$: we denote it by 
$c   \widetilde{T}$. Note that the function 
$(c , \widetilde{T}) \mapsto c \widetilde{T}$ is continuous. 
 
 \smallskip

\noi 
{\it Notation.} The symbol $*$ stands for the convolution product of measures on $\bR$. 
For every measure $\mathrm{m}$ on $\bR$, we set $\mathrm{m}^{*1}= \mathrm{m}$ and for all $n \geq 1$, we set $\mathrm{m}^{*(n+1)}= \mathrm{m}*\mathrm{m}^{*n}$. 
Also, $\lfloor \cdot \rfloor $ stands for the integer-part function.  

\paragraph{Convergence to Galton-Watson trees.} We first state a convergence result for GW-Markov chains that are rescaled in time {\it but not in space}. This result is quite close to Grimvall's \cite[Theorem 3.4]{Gr} that is actually a limit-theorem for GW-chains that are rescaled in time {\it and} space. 
Since it is not explicitly written in \cite{Gr}, we state it here as a lemma: its proof can be adapted in a straightforward way 
from that of \cite[Theorem 3.4]{Gr} that strongly relies on \cite[Theorem $2.2^\prime$]{Gr73} (whose proof also extends to our setting).  
\begin{lemma}
\label{discrbranchlim} Let $(\gamma_p)_{p\in \bN}$ be a positive sequence converging to $\infty$. Then, the following statements are equivalent. 
\begin{enumerate}
\item[\rm(i)] There exist a probability measure $\nu$ on $\bZ$ with $\nu(0)<1$ and a probability measure $\mu$ on $\bN$ with $\mu (0) < 1$, such that the following two limits hold in law as $p\rightarrow \infty$. 
$$ (\nu^p)^{*\lfloor \gamma_p\rfloor } \longrightarrow \nu  \quad {\rm and} \quad \mu^p \longrightarrow \mu . $$ 
\item[\rm(ii)] The one-dimensional marginal distributions of $(Y^{p}_{\lfloor \gamma_p a \rfloor })_{a\in [0, \infty)}$ converge to those of an $\bN$-valued process that is not constant.

\item[\rm(iii)] The process $(Y^{p}_{\lfloor \gamma_p a \rfloor })_{a\in [0, \infty)}$ converges weakly on 
$\bD ([0, \infty), \bR )$ to a continuous-time $\bN$-valued Galton-Watson branching process that is not constant. 
\end{enumerate}
\end{lemma}
\begin{remark}
\label{locrem}
Assume that $\rm(i)$, $\rm(ii)$ and $\rm(iii)$ hold true. Note that $\nu$ has to be an infinitely divisible distribution 
on $\bZ$. Then $\nu\neq\delta_0$ is the law of $X_1$, where $(X_t)_{t\in [0, \infty)}$ is a compound Poisson process with holding-time parameter $c$ and jump law $\pi$, a probability 
measure on $\bZ \backslash \{ 0 \}$. 
Denote by $(X^{p}_k)_{k\in \bN}$ a random walk with jump law 
$\nu^p $, and initial state $X^{p}_0= 0$. Standard arguments imply that 
$(X^{p}_{\lfloor \gamma_p t \rfloor })_{t\in [0, \infty)}$ converges weakly on $\bD ([0, \infty) , \bR)$ 
to $(X_t)_{t\in [0, \infty)}$. 
Since we deal with integer-valued processes, the joint law of the first jump time
and the size of the jump of $X^{p}$ converges to that of $X$. The assumptions on $\nu^p$ first imply that the support of $\pi$ is included in $\{ -1, 1, 2, \ldots \}$. If we set $\xi (k)= \pi (k- 1)$, $k\in \bN$, 
then it is easy to see that $\xi$ is the (proper and conservative) offspring distribution of the continuous-time $\bN$-valued GW-branching process mentioned in 
$\rm(iii)$ and that $c$ is its lifetime parameter. 
The previous joint convergence then entails 
that for all $k \in \bN$, 
\begin{equation}
\label{loclimi}
 \lim_{p \rightarrow \infty} \xi_{\#}^p(k) = \xi (k) \, , \quad   \lim_{p \rightarrow \infty} \mu^p (k)= \mu(k)  \, , \quad   
\lim_{p \rightarrow \infty} \gamma_p (1-\xi^p (1))= c \; , 
\end{equation}
where $\xi_{\#}^p (k)= \xi^p (k) /(1-\xi^p (1))$ if $k\neq 1$ and $\xi_{\#}^p (1)= 0$.  \cq 
\end{remark}

The previous result is used to derive the following limit theorem. Recall that 
$ \widetilde{\cF}_p$ stands for a random forest of $N_p$ independent discrete GW($\xi^p$)-real trees and that $N_p$ has law $\mu^p$.  
\begin{theorem}
\label{discrGWlim} Let $(\gamma_p)_{p\in \bN}$ be a positive sequence converging to $\infty$.\,We assume that the one-dimensional marginal distributions of $(Z^+_a( \frac{{1}}{{{\gamma_p}}} \widetilde{\cF}_p))_{a\in[0,\infty)}$ converge to those of an $\bN$-valued process that is not constant.
  Then, there exists a {\rm GW($\xi, c ; \mu$)}-real forest $\widetilde{\cF}$, 
as in Definition \ref{GWrealtreedef} with $\mu(0)<1$, such that  
$$ \textstyle\left(  \left(  Z^+_a \left( \frac{{1}}{{{\gamma_p}}}  \widetilde{\cF}_p \right) \right)_{a \in [0, \infty)}   \, ; \,   \frac{1}{{{\gamma_p}}} \widetilde{\cF}_p \,  \right)  \;  \overset{\textrm{\rm(law)}}{\underset{{p\rightarrow \infty}}{\longrightarrow}} \; \left( \left(  Z^+_a ( \widetilde{\cF})\right)_{a \in [0, \infty)}   \, ; \,    \widetilde{\cF} \, \right) \; ,$$
weakly on  $\bD( [0, \infty) , \bR) \times \bT$. 
Moreover, ${\rm(\ref{loclimi})}$ holds true. 
\end{theorem}
\begin{pf} First observe that $Y^{p}_k:=Z^+_{k\gamma_p} ( \frac{{1}}{{{\gamma_p}}} \widetilde{\cF}_p )$, $k\in \bN$, is a GW($\xi^p$)-Markov chain with initial distribution $\mu^p$. Thus, Lemma \ref{discrbranchlim} $\rm(i)$, $\rm(ii)$ and $\rm(iii)$ hold true.

We first prove convergence for single GW($\xi^p$)-discrete trees $\widetilde{\cT}_p$. This corresponds to the case where $\mu^p (1)= 1$, to which 
Lemma \ref{discrbranchlim} also applies. Thus, 
$(Z^+_a ( \frac{1}{{{\gamma_p}}}  \widetilde{\cT}_p) )_{a\in [0, \infty)}$
converges weakly on $\bD ([0, \infty) , \bR)$ as $p\rightarrow \infty$ to a continuous-time  
GW($\xi, c$)-branching process whose initial value is $1$. Moreover, $\xi$ is proper and conservative and (\ref{loclimi}) holds true. Since continuous-time GW-branching processes have no fixed-time discontinuity, convergence of finite-dimensional marginals holds true. 
This entails the tightness of the laws of $\frac{{1}}{{{\gamma_p}}}  \widetilde{\cT}_p$, $p\in\bN$, since for all $a , h\in [0, \infty)$, 
$ Z^{(h)}_a ( \frac{{1}}{{{\gamma_p}}}  \widetilde{\cT}_p) \leq  Z^+_a ( \frac{{1}}{{{\gamma_p}}}  \widetilde{\cT}_p )$, and by Theorem \ref{tightness}. 

We next want to prove that the real trees $\frac{{1}}{{{\gamma_p}}}  \widetilde{\cT}_p$ converge weakly to a GW($\xi, c$)-real tree by showing that every weak limit satisfies the branching property of 
Theorem \ref{vrairege}. 
To that end, observe first that the joint laws of $(  Z^+_a( \frac{{1}}{{{\gamma_p}}}  \widetilde{\cT}_p ))_{a \in [0, \infty)}$ and  $\frac{{1}}{{{\gamma_p}}} \widetilde{\cT}_p $, are tight on $\bD( [0, \infty) , \bR) \times \bT$. Let $(Y_a)_{a\in [0, \infty) }$ and $\widetilde{\cT}$ be such that 
\begin{equation}
\label{wsubsq}
 \textstyle \left(  \left(  Z^+_a \left( \frac{{1}}{{{\gamma_{p(k)}}}}  \widetilde{\cT}_{p(k)} \right) \right)_{a \in [0, \infty)}   \, ; \,   \frac{{1}}{{{\gamma_{p(k)}}}} \widetilde{\cT}_{p(k)} \,  \right)  \; \underset{{k\rightarrow \infty}}{\longrightarrow} \; \left( \left( Y_a \right)_{a \in [0, \infty)}   \, ; \,    \widetilde{\cT} \, \right) 
\end{equation} 
weakly on $\bD( [0, \infty) , \bR) \times \bT$ along the increasing sequence of positive integers 
$(p(k))_{k\in \bN}$. Without loss of generality, by the Skorohod representation theorem (and by a slight abuse of notation), we can assume that (\ref{wsubsq}) holds almost surely. To simplify notation we set $\widetilde{\fT}_k= \frac{{1}}{{{\gamma_{p(k)}}}}  \widetilde{\cT}_{p(k)} $. Thus, we assume that 
\begin{equation}
\label{subsqq}
 \left( \left(Z^+_a (\widetilde{\fT}_k)\right)_{a\in [0, \infty) } \, ; \,  \widetilde{\fT}_k \right)   \; \underset{{k\rightarrow \infty}}{\longrightarrow} \; \left( \left( Y_a \right)_{a \in [0, \infty)}   \, ; \,    \widetilde{\cT} \, \right) \quad \textrm{a.s. in $\bD( [0, \infty) , \bR) \times \bT$.}
 \end{equation} 
We first claim that 
\begin{equation}
\label{hardpoint}
\bP \left( Z_0^+ (\widetilde{\cT})= 1 \right) = 1 \quad {\rm and} \quad \bP \left( 
0 < D(\widetilde{\cT}) < \infty \right) = 1\; .
\end{equation}
\noi
\textit{Proof of (\ref{hardpoint}).} 
We first use the following standard result on Skorohod convergence for $\bN$-valued 
cadlag functions: by (\ref{subsqq}), 
the first jump time and the value at the first jump time of the processes 
$Z_{\cdot}^+ ( \widetilde{\fT}_k)$ converge a.s.~to the first jump time
and the value at the first jump time of the process $Y$. Observe that since
 $Z^+_0 (\widetilde{\fT}_k)= 1$, we get 
$$ D(\widetilde{\fT}_k)= \inf\{ a \in [0, \infty)\colon Z^+_a (\widetilde{\fT}_k) \neq 1 \} \quad {\rm and} \quad
 {\bf k} (\widetilde{\fT}_k)=   Z^+_{D(\widetilde{\fT}_k)} (\widetilde{\fT}_k)\;  . $$
Then, if we set $D^\prime=  \inf\{ a \in [0, \infty)\colon Y_a \neq 1 \} $ and ${\bf k}^\prime= Y_{D^\prime}$,  
the previous arguments and the continuity of $\abv$ entail that 
\begin{equation}
\label{ascvDk}
D^\prime =\lim_{k\rightarrow \infty}  D(\widetilde{\fT}_k) \, , \quad {\bf k}^\prime = \lim_{k\rightarrow \infty} {\bf k} (\widetilde{\fT}_k) \quad \textrm{and} \quad \abv(D^\prime , \widetilde{\cT} )= 
\textrm{$\dgh$--$\!\lim_{k\rightarrow \infty}$} \vartheta (\widetilde{\fT}_k) . 
\end{equation}
Since $Y$ is a GW($\xi, c $)-branching process, $D^\prime$ and ${\bf k}^\prime$ are independent, $D^\prime$ is exponentially distributed with mean $1/c$ and ${\bf k}^\prime$ has law $\xi$ that is proper and conservative. Thus, $0\! <  \! D^\prime\!<\!\infty$ and ${\bf k}^\prime \neq 1$ a.s.

  Let $(\fT_k, d_k, \rho_k)$ be any representative of $\widetilde{\fT}_k$ and let $(\cT, d, \rho)$ be any 
representative of $\widetilde{\cT}$. Then, $\bP$-a.s.~for any fixed $r\in (0, D^\prime)$, 
for all sufficiently large $k$ the closed ball $\overline{B}_{\fT_k} (\rho_k, r)$ rooted at $\rho_k$ is equivalent to the interval $[0, r]$ rooted at $0$ and since $\dgh(\fT_k, \cT) \rightarrow 0$, the closed ball 
$\overline{B}_{\cT} (\rho, r)$ rooted at $\rho$ is also equivalent to the interval $[0, r]$ rooted at $0$, which implies that $Z^+_0 (\cT)= 1$ and $D(\cT) \geq r$. This proves that 
\begin{equation}
\label{presque}
\textrm{$\bP$-a.s.} \qquad D(\widetilde{\cT}) \geq D^\prime >0 
\quad \textrm{and} \quad Z^+_0 (\widetilde{\cT}) = 1 \; .
\end{equation}
It remains to prove that $D(\widetilde{\cT}) \!<\! \infty$ a.s., and by (\ref{presque}), this boils down to proving that 
$\widetilde{\cT}$ is not the pointed isometry class of a half-line rooted at its finite end.   
 If ${\bf k}^\prime = 0$, then for any $r >D^\prime$, (\ref{ascvDk}) implies that for all sufficiently large $k$, ${\bf k} (\widetilde{\fT}_k)= 0$, which implies $\vartheta (\widetilde{\fT}_k)= \pnt$ and thus $\abv(D^\prime , \widetilde{\cT} )= \pnt$. Consequently, if ${\bf k}^\prime=0$, $\widetilde{\cT}$ is not the pointed isometry class of a half-line rooted at its finite end, and therefore $D(\widetilde{\cT} ) < \infty$.  

Since ${\bf k}^\prime\neq 1$ a.s., it only remains to 
consider the case where ${\bf k}^\prime \geq 2$. To that end, we first prove that 
\begin{equation} 
\label{ppresque} 
\textrm{$\bP$-a.s.~on $\{ {\bf k}^\prime \geq 2 \}$,} \quad \limsup_{k\rightarrow \infty} D(\vartheta \widetilde{\fT}_k ) >0 .  
\end{equation} 
Indeed, first note that for any $\varepsilon \in (0, \infty)$, 
$ \un_{\{  \limsup_{k\rightarrow \infty} D(\vartheta \widetilde{\fT}_k )< \varepsilon \, ; \, {\bf k}^\prime \geq 2 \}} \leq \liminf_{k\rightarrow \infty} \un_{\{  D(\vartheta \widetilde{\fT}_k ) \leq  \varepsilon  \, ; \, {\bf k} (\widetilde{\fT}_k)  \geq 2   \}} $. 
This inequality combined with Fatou lemma entails 
$$ \bP \big(  \limsup_{k\rightarrow \infty} D(\vartheta \widetilde{\fT}_k )< \varepsilon \, ; \, {\bf k}^\prime \geq 2 \big) \leq \liminf_{k\rightarrow \infty} \bP \big( D(\vartheta \widetilde{\fT}_k ) \leq  \varepsilon  \, ; \, {\bf k} (\widetilde{\fT}_k)  \geq 2\big) \; .$$
Now observe that $\bP \big( D(\vartheta \widetilde{\fT}_k ) \! \leq \!  \varepsilon  \, ; \, {\bf k} (\widetilde{\fT}_k) \! \geq \! 2\big) \!= \!  \bE \big[  \bP ( D( \widetilde{\fT}_k ) \! \leq \!  \varepsilon )^{ {\bf k} (\widetilde{\fT}_k) } \un_{\{{\bf k} (\widetilde{\fT}_k)  \geq 2 \}} \big]$, by the branching property for discrete Galton-Watson trees. This, combined with (\ref{ascvDk}) entails that 
$$ \lim_{k\rightarrow \infty} \bP \big( D(\vartheta \widetilde{\fT}_k ) \leq  \varepsilon  \, ; \, {\bf k} (\widetilde{\fT}_k) \geq 2\big)= \sum_{n\geq 2} \big( 1-e^{-c\varepsilon}\big)^n \xi (n)  
\quad \underset{\varepsilon \rightarrow 0}{-\!\!\!-\!\!\!-\!\!\!\longrightarrow} \; 0 \; , $$  
which implies (\ref{ppresque}). 

 Then, denote by $(\cT^\prime, d^\prime, \rho^\prime)$ a representative of $\abv (D^\prime, \widetilde{\cT})$. By (\ref{ppresque}), a.s.~if ${\bf k}^\prime \geq 2$, there exists $r >0$ such that $D(\vartheta \widetilde{\fT}_k) >r$ for infinitely many $k$, and (\ref{ascvDk}) implies that the closed ball 
 $\overline{B}_{\cT^\prime} (\rho^\prime , r)$ is equivalent to ${\bf k}^\prime$ copies of $[0, r]$ glued at $0$. It implies that $\widetilde{\cT}$ is not a half-line and thus,  $D(\widetilde{\cT} ) < \infty$. This completes the proof of the claim (\ref{hardpoint}). \cq 
 
\medskip

We next denote by $Q$ the law of $\widetilde{\cT}$. Then, (\ref{hardpoint}) can be rephrased as 
\begin{equation}
\label{hhardpnt}
Q(Z^+_0 = 1) = 1 \quad \textrm{and} \quad Q (\, 0 \! < \! D \! < \! \infty \, ) = 1 \; .
\end{equation}
Let us also denote by $Q_{p(k)}$ the law of $\widetilde{\fT}_k$.  
We fix $a\in [0, \infty)$ and for all $k\in \bN$, we set $a_k:= \lfloor  \gamma_{p(k)} a \rfloor / \gamma_{p(k)} \rightarrow a$, as $k\rightarrow \infty$.  
Let $G_1, G_2\colon  \bT \rightarrow [0,\infty)$ be continuous and bounded. First observe that 
$$ \bE \left[ G_1 \left(  \blw (a_k, \widetilde{\fT}_k ) \right)  G_2 \left(  \abv (a_k, \widetilde{\fT}_k )\right) \right]=  \bE \left[ G_1\!\left(  \blw (a_k, \widetilde{\fT}_k) \right)  Q_{p(k)}^{\circledast Z^+_{a_k} (\widetilde{\fT}_k) } \left[ G_2 \right] \right] \, , $$
by the branching property for discrete Galton-Watson trees. As $k \rightarrow \infty$, (\ref{subsqq}) and the continuity of $\blw$ and $\abv$ stated in Lemma \ref{abovebelow} imply that 
\begin{equation}
\label{subbr}
Q\! \left[ G_1 (  \blw (a  ,  \cdot   ) )  
G_2 (  \abv (a , \cdot   ) ) \right]\!=\!  \bE\! \left[ G_1\! \left(\!  \blw (a ,\widetilde{\cT} )\! \right)  
G_2\! \left( \! \abv (a, \widetilde{\cT}) \!\right)\! \right]\!=\!  \bE \! \left[ G_1 (  \blw (a, \widetilde{\cT} ) )  
Q^{\circledast Y_a  } [ G_2 ] \right] .\hspace{-0.1cm}
\end{equation}
This equality extends to all nonnegative measurable functions $G_1$ and $G_2$. 
Since $Q(Z^+_0 = 1) = 1$, (\ref{subbr}) with $G_1 \equiv 1$ and $G_2 (\cdot)= f( Z^+_0 (\cdot))$ first implies $Q [ f(Z^+_a )] = \bE [f(Y_a)]$. This first proves that $Q( Z^+_a <\infty)=1 $ for any 
$a\in [0, \infty)$, and it also entails for any measurable functions $G: \bT \rightarrow [0, \infty)$ and $f: \bN \rightarrow [0, \infty)$ that 
$$ Q\left[ f( Z^+_a) G (\abv(a \, , \, \cdot \, ) ) \right] = \bE \left[ f(Y_a) Q^{\circledast Y_a  } \left[ G \right] \right] = 
Q\left[ f(Z^+_a) Q^{\circledast Z^+_a} \left[ G \right]  \right]  . $$
Therefore, the conditional distribution given $Z^+_a$ of $\abv(a\, , \, \cdot \, )$ under $Q$ is $Q^{\circledast Z^+_a} $. 

 We thus have proved that Theorem \ref{vrairege} applies to $Q$ that is therefore 
the law of a GW($\xi^\prime , c^\prime $)-real tree. Since $Y_a$ under $\bP$ has the same law as 
$Z^+_a $ under $Q$, we easily get $(\xi^\prime , c^\prime)= (\xi, c)$.

Since $\widetilde{\cT}$ is a GW($\xi, c$)-real tree, $a\mapsto  Z^+_a (\widetilde{\cT})$ is cadlag 
and it has no fixed discontinuity. Thus, for any $a \in (0, \infty)$, this process is left-continuous at time $a$ a.s., and $Z^+_a (\widetilde{\cT})$ is equal a.s.~to a measurable functional $G_1$ of $\blw(a, \widetilde{\cT})$. The previous arguments and (\ref{subbr}) entail that for any $a\in [0, \infty)$, $Z^+_a (\widetilde{\cF})= Y_a$ a.s., and since both processes are cadlag, we get 
$Z^+_\cdot (\widetilde{\cT})= Y$ a.s. This proves the uniqueness of the limit in the joint convergence and it actually proves the theorem in the case of single trees, namely when $\mu^p= \delta_1$, $p \in \bN$. 

\medskip

Let us prove the general case. 
Denote by $Q_p$ the law of $\frac{1}{{\gamma_p}}\widetilde{\cT}_p$. 
Let $F\colon \bD([0, \infty) , \bR) \! \times \! \bT\rightarrow \bR$ be bounded and 
continuous for the product topology. 
We proved that $\int F( Z^+_{\cdot } (\widetilde{T}) ; \widetilde{T}) Q_p(d\widetilde{T}) \longrightarrow \int F( Z^+_{\cdot } (\widetilde{T}) ;  \widetilde{T})  Q_{\xi, c}(d\widetilde{T}) $. Now observe that
for any $k\in \bN$, 
$$\int_{\bT} F\left( Z^+_{\cdot } (\widetilde{T})\, ;\,  \widetilde{T}\right) \, Q^{\circledast k}_p(d\widetilde{T}) = 
\int_{\bT^k} F\left( Z^+_{\cdot } (\widetilde{T}_1) + \cdots + Z^+_{\cdot } (\widetilde{T}_k)  \, ;\,   \widetilde{T}_1\circledast \cdots \circledast \widetilde{T}_k \right) \, Q^{\otimes k}_p (d\widetilde{T}_1 \ldots d\widetilde{T}_k) $$
Since $k$ independent continuous-time GW-branching processes have distinct jump times a.s., we easily get 
\begin{eqnarray*}
 \lim_{p\rightarrow \infty}  \int_{\bT^k} F\left( Z^+_{\cdot } (\widetilde{T}_1) + \cdots + Z^+_{\cdot } (\widetilde{T}_k) \,  ; \,   \widetilde{T}_1\circledast \cdots \circledast \widetilde{T}_k \right) \, Q^{\otimes k}_p (d\widetilde{T}_1 \ldots d\widetilde{T}_k) =& & \\
& & \hspace{-70mm} =\int_{\bT^k} F\left( Z^+_{\cdot } (\widetilde{T}_1) + \cdots + Z^+_{\cdot } (\widetilde{T}_k)\,  ;  \,  \widetilde{T}_1\circledast \cdots \circledast \widetilde{T}_k \right) \, Q^{\otimes k}_{\xi, c} (d\widetilde{T}_1 \ldots d\widetilde{T}_k)\\
& & \hspace{-70mm} = \int_{\bT} F\left( Z^+_{\cdot } (\widetilde{T}) \, ; \,  \widetilde{T}\right) \, 
Q^{\circledast k}_{\xi, c}(d\widetilde{T}) 
\end{eqnarray*}
This implies the following: 
\begin{eqnarray*}
\lim_{p\rightarrow \infty} \bE \left[ F\left( Z^+_\cdot  \left(\frac{_{1}}{^{\gamma_p}} 
\widetilde{\cF}_p \right) \, ; \, \frac{_{1}}{^{\gamma_p}} 
\widetilde{\cF}_p \right)  \right] &=& \lim_{p\rightarrow \infty} \sum_{k\in \bN} \mu^p (k) \int_\bT 
F\left( Z^+_{\cdot } (\widetilde{T})\, ;\,  \widetilde{T}\right) \, Q^{\circledast k}_p(d\widetilde{T}) \\
& =& \sum_{k\in \bN} \mu (k) \int_\bT 
F\left( Z^+_{\cdot } (\widetilde{T})\, ;\,  \widetilde{T}\right) \, Q^{\circledast k}_{\xi, c}(d\widetilde{T}) \, ,
\end{eqnarray*}
which completes the proof of the theorem. 
\end{pf}

\paragraph{Convergence to L\'evy forests.} We consider now the convergence of GW-trees 
to L\'evy forests. In these cases, the profiles of the trees are rescaled in time {\it and} space. More precisely, 
we make the two following assumptions. 
\begin{enumerate}
\item[\textrm{(A1)}] There is a positive sequence $(\gamma_p)_{p\in \bN}$ 
converging to $\infty$, such that the process $(\frac{1}{p}Y^{p}_{\lfloor  \gamma_p a \rfloor})_{a\in [0,  \infty)}$ converges weakly on $\bD([0, \infty) , \bR)$ to a CSBP($\psi, \varrho$), 
where $\varrho (\{ 0 \}) <1$ and $\int_{0+} \frac{dr}{ (\psi (r))_-}= \infty$.  

\item[\textrm{(A2)}] Set $\cE_p= \inf \{ k\in \bN \colon  Y^{(p)}_k = 0 \} $ and denote by $\cE$ the extinction time of a CSBP($\psi, \varrho$).  We assume that $\int^\infty  \frac{dr}{\psi (r)} < \infty$ and that $\frac{1}{\gamma_p}\cE_p \longrightarrow \cE $, weakly on $[0, \infty]$, as $p\rightarrow \infty$. 
 \end{enumerate}
  
Grimvall \cite[Theorem 3.4]{Gr} asserts that \textrm{(A1)} is equivalent to the following weak convergence on $\bR$: 
$$   \textstyle\nu^p \left( \frac{{\, \cdot \,  }}{{p}}\right)^{*p\lfloor \gamma_p\rfloor } \;  
\underset{{p\rightarrow \infty}}{\longrightarrow} \; \nu \quad {\rm and} \quad 
\mu^p \left(\,  \frac{{\, \cdot \,  } }{{p}}\, \right)  \;  
\underset{{p\rightarrow \infty}}{\longrightarrow} \; \varrho \; , $$
where $\nu$ is an infinitely divisible spectrally positive law such that $\int_\bR  e^{-\theta x} \nu (dx)= e^{\psi (\theta)}$. Analytic necessary and sufficient conditions equivalent to such a convergence can be found for instance in \cite[Theorem II.3.2]{Jacod}.

If we assume (A1) and $\int^\infty  \frac{dr}{\psi (r)} < \infty$, then we can show that (A2) is equivalent to the following 
$$ \liminf_{p\rightarrow \infty} \left( \varphi_{\xi^p}^{\circ \lfloor \gamma_p \rfloor } (0) \right)^p > 0 \; , $$
where $ \varphi_{\xi^p}^{\circ n}$ stands for the $n$-th iterate of $\varphi_{\xi^p}$.  
We leave the details to the reader (see also the comments following \cite[Theorem 2.3.1]{DuLG}). Let us mention that (A1) implies (A2) when $\xi^p= \xi$, for all $p\in \bN$. In this case, 
$\psi$ is necessarily a $\gamma$-stable branching mechanism, namely $\psi(\lambda)= \lambda^\gamma$, $\gamma \in (1, 2]$ (see \cite[Theorem 2.3.2]{DuLG}).  

  We now state the main result of the section. To that end, 
recall from (\ref{equaintegrale}) and (\ref{vvvdefini}) the notation 
$u(a, \theta)$ and $v(a)$, recall from Theorem \ref{growthconvergence} the definition 
of L\'evy forests and recall 
that $\widetilde{\cF}_p$ stands for a random forest of $N_p$ independent discrete GW($\xi^p$)-real trees, where $N_p$ has law $\mu^p$.  
\begin{theorem}
\label{invprincth} Assume ${\rm (A1)}$ and ${\rm (A2)}$. Then, there exists a $(\psi, \varrho)$-L\'evy forest $\widetilde{\cF}$  and a 
{\rm CSBP}($\psi, \varrho$) denoted by $(Z^\varrho_a)_{a\in [0, \infty)}$ such that for all $a\in [0, \infty)$, 
$\bP$-a.s.\ $\lim_{h\rightarrow 0+} \frac{1}{v(h)}Z^{(h)}_a (\widetilde{\cF})= Z^\varrho_a$, and such that 
weakly on $\bD([0, \infty) , \bR) \times \bT$, 
\begin{equation}
\label{expliicv}\textstyle
 \left(  \left( \,  \frac{1}{p}Z^+_a \left( \frac{{1}}{{{\gamma_p}}}  \widetilde{\cF}_p \right) \right)_{a \in [0, \infty)}   \, ; \,   \frac{{1}}{{{\gamma_p}}} \widetilde{\cF}_p \,  \right)  \;  \overset{\textrm{\rm(law)}}{\underset{{p\rightarrow \infty}}{\longrightarrow}} \; \left( (  Z^\varrho_a)_{a \in [0, \infty)}   \, ; \,    \widetilde{\cF} \, \right) \; .
\end{equation} 
\end{theorem}
\begin{remark}
\label{LGDvsWD}
We first mention that a closely related result has been proved by quite different methods 
in \cite[Theorem 2.3.1 and Corollary 2.5.1]{DuLG}: this result only deals with 
critical or sub-critical GW-trees, however the convergence holds for the contour process, which is a stronger convergence. Note that in the super-critical cases, the contour process is not a 
well-suited approach. \cq 
\end{remark}
\begin{remark}
\label{miniass}
As noticed in \cite{DuLG}, (A1) and (A2) are in some sense the minimal assumptions under which 
the convergence (\ref{expliicv}) holds. Indeed, (A1) and $\int^\infty \frac{dr}{\psi (r)} <\infty$, do not necessarily imply (A2): see \cite[pp.\ 60-61]{DuLG} for a counterexample. 
Note that $(\psi, \varrho)$-L\'evy forests can only be defined as locally compact real trees  
if $\int^\infty \frac{dr}{\psi (r)} <\infty$. Moreover, if we assume (A1), $\int^\infty \frac{dr}{\psi (r)} <\infty$ and $ \frac{{1}}{{{\gamma_p}}} \widetilde{\cF}_p \rightarrow \widetilde{\cF}$, then, by the $\dgh$-continuity of $\Gamma$, the total heights converge too, which implies (A2). \cq 
\end{remark}
\begin{pf} Recall that $\widetilde{\cT}_p$ stands for a single discrete GW($\xi^p$)-real tree and observe that 
$Z^+_{k\gamma_p}( \frac{{1}}{{{\gamma_p}}} \widetilde{\cT}_p)$, is a discrete-time GW($\xi^p$)-Markov chain whose initial state is equal to $1$.
For all $a, \theta , h \in [0, \infty)$ and all $p\in \bN$, we set 
$$\textstyle u_p (a, \theta) = - p\log \bE [ \exp(-\frac{1}{p} \theta Z^+_{a} 
( \frac{{1}}{{{\gamma_p}}} \widetilde{\cT}_p )  )] \quad {\rm and} \quad  
v_p (h)= -p \log \bP (\,  \Gamma ( \frac{{1}}{{{\gamma_p}}} \widetilde{\cT}_p  ) < h \,   ) \; .$$
Next, observe that $Y^{p}_k:=Z^+_{k\gamma_p}( \frac{{1}}{{{\gamma_p}}} \widetilde{\cF}_p)$, $k\in \bN$, is a discrete-time GW($\xi^p$)-Markov chain with initial distribution $\mu^p$, to which (A1) and (A2) apply. Then, by (A1),  
$$\textstyle \varphi_{\mu^p} (\exp(-\frac{1}{p}u_p (a, \theta)) )  =  \bE\left[ \exp(-\frac{1}{p}\theta  Y^{(p)}_{\lfloor \gamma_p a \rfloor } )   \right] \quad  \underset{{p\rightarrow \infty}}{\longrightarrow} \quad \int_{[0, \infty)} 
\varrho (dy) \,  e^{-y u(a, \theta)} $$
and by (A2) and (\ref{tspdextinct}), we get 
$$\textstyle  \varphi_{\mu^p} (\exp(-\frac{1}{p}v_p (h))) =   
\bP (\,  \frac{1}{{\gamma_p}} \cE_p < h \, ) \quad   \underset{{p\rightarrow \infty}}{\longrightarrow} \quad 
 \int_{[0, \infty)} \varrho (dy) \, e^{-y v(h)} \; . $$
We next use the following basic result (known as the second theorem of Dini). 
\begin{enumerate}
\item[(R)] Let $f_p\colon  [0, \infty) \rightarrow 
[0, \infty)$, $p\in \bN$, be a sequence of monotonic functions converging pointwise to a continuous function 
$f$. Then the convergence is uniform on every compact interval of $[0, \infty)$. 
\end{enumerate}
We first note that the functions $f_p (\theta)= \varphi_{\mu^p} (e^{-\theta/p}) $ and $f(\theta)= \int_{[0, \infty)} \varrho (dy) e^{-\theta y}$ are strictly monotonic and continuous, so that the inverses $f_p^{-1}$ converge to $f^{-1}$ pointwise on the interval $(\varrho( \! \{0\} \! ) \, , \,\!  1]$. Hence, for all $a, \theta \in [0, \infty) $ and for all $h \in (0, \infty)$, we get $u_p (a, \theta) \rightarrow u(a, \theta)$ and $v_p (h) \rightarrow v(h)$. 
Next observe that $u_p$ is monotone in each component and that $v_p$ is non-increasing. 
Thus, by  (R), 
\begin{equation}  
\label{expoconv}
(a_p, \theta_p , h_p) \! \underset{p\rightarrow\infty}{\longrightarrow} \! (a, \theta, h) \!\in \! [0, \infty)^2\! \times\! (0, \infty) \Longrightarrow 
\; u_p (a_p , \theta_p)     \underset{{p\rightarrow \infty}}{\longrightarrow}  u(a, \theta) 
\quad {\rm and} \quad v_p (h_p)   \underset{{p\rightarrow \infty}}{\longrightarrow}  v(h). 
\end{equation}
We next fix $h\in (0, \infty)$ and we set $\widetilde{\cF}^h_p :=R_{\lfloor \gamma_p h \rfloor } 
(\widetilde{\cF}_p )$. 
It is easy to see that $\widetilde{\cF}^h_p$ is a forest of discrete GW-real trees as defined at the beginning 
of the section. We want to apply Theorem \ref{discrGWlim} to $\widetilde{\cF}^h_p$. 
To that end, first denote by $\xi^p_h$ its offspring distribution and by $\mu^p_h$ the law of the number of trees in this forest. We do not need to compute them explicitly to see that 
$\xi^p_h (1) < 1$ and $\mu^p_h (0) <1$. Next fix $a\in [0, \infty)$ and observe that conditionally given $Y^{p}_{\lfloor \gamma_p a \rfloor} \! =\!  Z^+_{\lfloor \gamma_p a \rfloor} (\widetilde{\cF}_p) \! =\!  n$, the law of 
$Z^+_{\lfloor \gamma_p a \rfloor} (\widetilde{\cF}^h_p)$ is binomial with parameters $n$ and 
$\bP(\Gamma (\widetilde{\cT}_p) \! \geq \! \lfloor \gamma_p h \rfloor)$. To simplify notation, we set 
$h_p=  \lfloor \gamma_p h \rfloor /  \gamma_p$, which tends to $h$ as $p\rightarrow \infty$, and 
we define $\phi_p (h, \theta) \in [0, \infty)$ by   
$$ \textstyle\exp(-\frac{1}{p} \phi_p ( h, \theta ) )= \bP\big( \Gamma (\widetilde{\cT}_p) \! <\!  \lfloor \gamma_p h \rfloor \big) + e^{-\theta} \bP\big( \Gamma (\widetilde{\cT}_p) \! \geq \! \lfloor \gamma_p h \rfloor \big) = 1\!-\! ( 1\!-\! e^{-\theta}) (1 \!-\!  \exp(-\frac{1}{p} v_p (h_p))) . $$
By  (\ref{expoconv}), $\phi_p (h, \theta) \rightarrow v(h)  (1-e^{-\theta}) $ and $\bE [\exp(- \theta Z^+_{a}(\frac{{1}}{{{\gamma_p}}} \widetilde{\cF}^h_p )) ] \longrightarrow  \int_{[0, \infty)} \varrho (dy) e^{-y u(a , v(h) (1-e^{-\theta}) \, ) } $. 
Theorem \ref{discrGWlim} applies and the following convergence holds weakly on 
$\bD ([0, \infty), \bR) \times \bT$: 
$$ \textstyle \left(  (Z^+_{a}(\frac{{1}}{{{\gamma_p}}} \widetilde{\cF}^h_p ))_{a\in [0, \infty)} ;  
\frac{{1}}{{{\gamma_p}}} \widetilde{\cF}^h_p  \right) \quad   \underset{{p\rightarrow \infty}}{\longrightarrow} \quad  \left( (Z^+_a(\widetilde{\cF}^h))_{a\in [0, \infty)} ; \widetilde{\cF}^h \right) \; .$$
Here, it is easy to see that 
$\widetilde{\cF}^h$ is a GW($ \xi_{v(h)} ,   c_{v(h)};  \mu_{v(h)}$)-real forest, where we recall that the one-parameter family of laws $(\xi_\lambda, c_\lambda, \mu_\lambda)$ is derived from $\psi$ and 
$\varrho$ by (\ref{xilambdef}). 

Now observe that $ \dgh( R_h (\frac{{1}}{{{\gamma_p}}} \widetilde{\cF}_p  ) , \frac{{1}}{{{\gamma_p}}} \widetilde{\cF}^h_p ) \leq \frac{1}{{\gamma_p}} $. Moreover, for each fixed $a$, 
there exists a random variable $\Delta_p\colon \Omega \rightarrow \bN$ such that $ \Delta_p \!\geq\!  
Z^{+}_{a}(\frac{{1}}{{{\gamma_p}}} \widetilde{\cF}^{h}_{p} )\!-\! Z^{{(h)}}_{{a}} \! (\frac{{1}}{{{\gamma_p}}} \widetilde{\cF}_p) \! \geq \!0$, and such that conditionally given 
$ Z^+_{\lfloor \gamma_p a \rfloor} (\frac{{1}}{{{\gamma_p}}}\widetilde{\cF}^h_p)= n$, the law of 
$\Delta_p$ is binomial with parameters $n$ and 
$\bP\big( h_p \!<\!  \Gamma (\frac{{1}}{{\gamma_p}}\widetilde{\cT}_p) \! \leq\!  h_p + \frac{{2}}{{\gamma_p}} \big)$, 
which tends to $0$ as $p\rightarrow \infty$, by (\ref{expoconv}). Since the laws of the random variables 
$ Z^+_{\lfloor \gamma_p a \rfloor} (\frac{{1}}{{{\gamma_p}}}\widetilde{\cF}_p^h)$ are tight, we get 
$\Delta_p \rightarrow 0$, in probability. Thus, for all $h \in (0, \infty)$ and all $a\in [0, \infty)$, 
$( Z^{(h)}_a (\frac{{1}}{{{\gamma_p}}} \widetilde{\cF}_p) \, ;  R_h (\frac{{1}}{{{\gamma_p}}} \widetilde{\cF}_p  ) ) $ converges to $(  Z^+_a (\widetilde{\cF}^h) \, ; \widetilde{\cF}^h )$ weakly on $\bN \times \bT$, as $p\rightarrow \infty$. 

  Let $\widetilde{\cF}$ be a $(\psi, \varrho)$-L\'evy forest. Theorem \ref{growthconvergence} $\rm(ii)$ asserts that $R_h (\widetilde{\cF})$ and $\widetilde{\cF}^h$ have the same law for all $h\in (0, \infty)$. Then, 
by Corollary \ref{convcriter}, $\frac{{1}}{{{\gamma_p}}} \widetilde{\cF}_p \rightarrow \widetilde{\cF}$ weakly on $\bT$ as $p\rightarrow \infty$. Thus, we have proved 
\begin{equation}
\label{prejointcv}
\textstyle \frac{{1}}{{{\gamma_p}}} \widetilde{\cF}_p \; 
\overset{\textrm{(law)}}{\underset{p\rightarrow \infty}{-\!\!\!-\!\!\!\longrightarrow}} \;  \widetilde{\cF} \qquad \textrm{and}  \qquad \big( \,  Z^{(h)}_a (\frac{{1}}{{{\gamma_p}}} \widetilde{\cF}_p) \, ;  R_h \big( \frac{{1}}{{{\gamma_p}}} \widetilde{\cF}_p   \big) \big) \;  \overset{\textrm{(law)}}{\underset{{p\rightarrow \infty}}{-\!\!\!-\!\!\!\longrightarrow}}  \; \big( \, Z^{(h)}_a (\widetilde{\cF}) \, ; \,  R_h \big(\widetilde{\cF} \big) \big) \; .
\end{equation}
Then, note that the laws of 
$(\, (\frac{1}{p} Z^+_a (\frac{{1}}{{{\gamma_p}}} \widetilde{\cF}_p))_{a\in [0, \infty)} \, ; \frac{{1}}{{{\gamma_p}}} \widetilde{\cF}_p )$, $p\in \bN$, are tight on $\bD ([0, \infty) , \bR) \times \bT$. Let 
$(p(k))_{k\in \bN}$ be an increasing sequence of integers such that 
\begin{equation}
\label{subsqone}
\textstyle \left(  \left( \,  \frac{1}{p(k)}Z^+_a \left( \frac{{1}}{{{\gamma_{p(k)}}}}  \widetilde{\cF}_{p(k)} \right) \right)_{b \in [0, \infty)}   \, ; \,   \frac{{1}}{{{\gamma_{p(k)}}}} \widetilde{\cF}_{p(k)} \,  \right)  \;  \overset{\textrm{\rm(law)}}{\underset{{k\rightarrow \infty}}{-\!\!\!-\!\!\!\longrightarrow}} \; \left( (  W_b)_{b \in [0, \infty)}   \, ; \,    \widetilde{\cF} \, \right) \;
\end{equation}
where $(W_b)_{b \in [0, \infty)}$ is a CSBP($\psi, \varrho$). By Theorem \ref{asconvgrowth}, 
the proof of the joint convergence will be complete if we show that for each $a\in [0, \infty)$,  
$Z^{(h)}_{a} (\widetilde{\cF}) \rightarrow W_a$ in probability as $h\rightarrow 0$. 
To that end, first note that by (\ref{prejointcv}), there exists an increasing sequence of integers $(k_\ell)_{\ell\in \bN}$ such that the following limit holds true weakly on $\bN^2 \times \bT$ 
\begin{equation}
\label{suuub}
\left( Z^{(h)}_{a} \left(\frac{{1}}{{{\gamma_{p(k_\ell)}}}} \widetilde{\cF}_{p(k_\ell)}\right) \, ;  
\frac{1}{{p(k_\ell)}}   Z^{+}_a \left(\frac{{1}}{{{\gamma_{p(k_\ell)}    }}} \widetilde{\cF}_{p(k_\ell)}\right) \, ;  \frac{{1}}{{{\gamma_{p(k_\ell)}}}} \widetilde{\cF}_{p(k_\ell)}\right)   \quad   \underset{{\ell \rightarrow \infty}}{\longrightarrow} \quad  \left(  
X \, ; W_a \, ;  \widetilde{\cF}\,  \right) . 
\end{equation}
The $\dgh$-continuity of $R_h$ and (\ref{prejointcv}) imply that $(X, R_h(\widetilde{\cF}))$ and 
$(Z^{(h)}_a (\widetilde{\cF}) , R_h(\widetilde{\cF}))$ have the same law, which implies that $X= Z^{(h)}_a (\widetilde{\cF}) $ a.s.  

Next, observe that the conditional law of $Z^{(h)}_a (\frac{1}{{\gamma_p}} \widetilde{\cF}_p )$ given $Z^+_a (\frac{1}{{\gamma_p}}\widetilde{\cF}_p )= n$, is binomial with parameters $n$ and 
$\bP( \Gamma (\widetilde{\cT}_p ) \geq \gamma_p h + \gamma_p a -\lfloor \gamma_p a \rfloor)$. 
Previous computations and (\ref{suuub}) imply that for all $\theta \in [0, \infty) $ and all bounded continuous functions $f$, 
$$\bE \left[ e^{-\theta Z^{(h)}_a (\widetilde{\cF})} f \left( W_a \right) \right] = \bE \left[  e^{-W_a v(h) (1-e^{-\theta})} f \left( W_a \right) \right]  . $$
Thus, the conditional law of $Z^{(h)}_a(\widetilde{\cF})$ given $W_a$ is a Poisson distribution with parameter 
$v(h) W_a$. This implies that for all $a \in [0, \infty)$ and all $h, K , \varepsilon  \in (0, \infty)$, 
$$\textstyle \bP(\, |  W_a-\frac{1}{{v(h)}} Z^{(h)}_a(\widetilde{\cF}) | > \varepsilon   )  \leq  \frac{K}{\varepsilon^2 v(h)} + \bP (W_a >K) \; , $$      
which implies the desired result. 
\end{pf}

\appendix

\section{Proofs of preliminary results on real trees.}\label{appendixA}
Let us first recall basic results on the Gromov-Hausdorff metric (in the context of real trees).  
Let $(T, d, \rho)$ and $(T^\prime, d^\prime, \rho^\prime)$ be two CLCR real trees and let $\varepsilon \in (0, \infty)$. A function $f\colon T\rightarrow T^\prime$ is a {\it pointed $\varepsilon $-isometry} if it satisfies the following 
conditions.
\begin{enumerate}
\item[\rm (a)]  $f(\rho)=\rho^\prime$. 
\item[\rm(b)] $ \textrm{dis} (f):=\sup \left\{ \lvert d( \s,s)- d^\prime(f(\s),f(s)) \rvert \, ;\;  
  \s , s\in T \right\}  < \varepsilon $. This quantity is called the {\it distortion} of $f$. 
\item[\rm(c)] $f(T)$ is a $\varepsilon $-net of $T^\prime$. Namely, every point of $T^\prime$ is at 
distance at most $\varepsilon$ of $f(T)$.
\end{enumerate}
The following lemma is a translation into our tree context of \cite[Corollary 7.3.28]{BuBu}. 
\begin{lemma}
\label{epsiiso} 
If $\dgh_{{\rm cpct}} (T, T') < \varepsilon$, then there exists a pointed $4\varepsilon$-isometry from 
$T$ to $T'$. If there exists a pointed $\varepsilon$-isometry from 
$T$ to $T'$, then $\dgh_{{\rm cpct}} (T, T') < 4\varepsilon$.
\end{lemma}
Recall (\ref{branchheight}) that gives the height of the branch point 
of two points $\s , s \in T$. 
If $f\colon T\rightarrow T^\prime$ is a pointed $\varepsilon $-isometry then, 
(\ref{branchheight}) implies that 
\begin{equation}
\label{isobranch}
\forall \s , s \in T \, , \quad \left| d( \rho , \s \wedge s) -d^\prime (\rho^\prime , f(\s) \wedge f(s)) \right| < 
\textstyle\frac{3}{2} \varepsilon \; .
\end{equation}

\subsection{Proof of Lemma \ref{abovebelow}.} 
\label{AbvBlwpfsec} 
Observe that for every CLCR real tree $(T, d, \rho)$, $\blw (a, T)$ is compact and 
that the ball of centre $\rho$ with radius $r$ of $\abv (a, T)$ is equal to $\abv ( a, \overline{B}_T (\rho , a+r))$. Thus, without loss of generality, we only need to consider compact real trees.

  For all $a, b\in (0, \infty)$, we easily see that $  \dgh_{{\rm cpct}} \left(\blw (a, T) , \blw (b, T) \right) \leq \lvert  a-b\rvert $. Assume that $(T^\prime , d^\prime, \rho^\prime) $ is a compact rooted real tree. 
The definition (\ref{dghcpct}) of $\dgh_{{\rm cpct}}$ easily entails that 
$\dgh_{{\rm cpct}} (\blw (a, T), \blw (a, T^\prime)) \leq 3 \dgh_{{\rm cpct}} (T, T^\prime)$ for all $a\in (0, \infty)$. Thus, 
$$   \dgh_{{\rm cpct}} \left(\blw (a, T) , \blw (b, T^\prime) \right) \leq \lvert  a-b\rvert +3 \dgh_{{\rm cpct}} (T, T^\prime) \; , $$
which entails the joint continuity for $\blw$. 

Next, define the following pseudo-metric $d_a$ on $T \times T$ by
\begin{equation}
\label{pseudoabove}
d_a (\sigma , s)= a \vee d(\rho, \sigma) +  a \vee d(\rho, s) -  
2 \left( a \vee  d(\rho, \sigma \wedge s ) \right) 
\end{equation}
and say that $\s \equiv s$ iff $d_a(\s, s)= 0$. Let 
$\rho_a$ be the equivalence class of $\rho$. Then, $(T/ \! \equiv_{d_a} , d_a, \rho_a)$ is isometric to $\abv (a, T)$.  Note that $0 \leq d(\s , s) -d_a (\s, s) \leq 2a $, which easily implies that the canonical projection from $T$ to $T/\equiv_{d_a}$ is a $2a$-pointed isometry and by Lemma \ref{epsiiso}, we get $\dgh_{{\rm cpct}} (T, \abv (a, T)) \leq 8 a$. Since $\abv (a_1+ a_2, T) = \abv (a_1 , \abv (a_2 , T))$, we easily get 
\begin{equation}
\label{Treefixed}
\forall a, b \in [0, \infty) \; , \quad \dgh_{{\rm cpct}} \left(  \abv (a , T) ,  \abv (b , T) \right)  \leq 8 
\lvert a-b \vert \; .
\end{equation}
Let $(T^\prime, d^\prime, \rho^\prime)$ be a compact rooted real tree and let $\varepsilon > \dgh_{{\rm cpct}} (T, T^\prime)$. Lemma \ref{epsiiso} implies that there exists a 
$4\varepsilon$-pointed isometry $f$ from $T$ to $T^\prime$ and 
(\ref{pseudoabove}) and (\ref{isobranch}) entail $| d_a  (\sigma, s )-d^\prime_a( f(\sigma), f(s))| \leq 20 \varepsilon$, for all $\s , s \in T$. 
An easy argument shows that $f$ induces a pointed $28\varepsilon$-isometry from $\abv(a,T)$ to 
$\abv (a, T^\prime)$. By Lemma \ref{epsiiso} and (\ref{Treefixed}) we get 
$$\forall a, b \in [0, \infty) \; , \quad \dgh_{{\rm cpct}} \left(  \abv (a, T) ,  \abv (b , T^\prime) \right)  \leq 112 \, \dgh_{{\rm cpct}} (T, T^\prime) + 8 \lvert a-b \vert , $$
which completes the proof of Lemma \ref{abovebelow} \cqfd

\subsection{Proof of Lemma \ref{measuraZht}.} 
\label{mesusec}
\paragraph{Proof of Lemma \ref{measuraZht} $(i)$.} 
First note that for all $a, h \in [0, \infty)$ and every CLCR real tree $T$, 
 $Z^{+}_a (T)= \langle  \cM_a (T)  \rangle$ and $Z^{(h)}_{a+} (T)=\langle \cM_a(R_h(T))\rangle $.
Next, observe that $\cM_a(T)= \cM_0(\abv (a, T))$. So we only need to prove that $\cM_0$ 
is measurable. The definition of $\dgh$ entails that we only need to prove that the restriction of 
$\cM_0$ to $(\bT_{{\rm cpct}} , \dgh_{{\rm cpct}})$ is measurable. To that end, we first 
prove the following claim. 

\noi
\underline{\textsl{Claim 1.}} For all  
bounded Lipschitz functions 
$F\colon \bT_{{\rm cpct}} \rightarrow [0, \infty)$ and $\phi\colon [0, \infty) \rightarrow [0, \infty)$, with $\phi$ vanishing in a neighbourhood of $0$, 
and for all sufficiently small $h \in (0, \infty)$, 
\begin{equation}
\label{claim1}
\tag{ \textsl{Claim 1}}
\langle (\cM_0 \circ R_h)(\, \cdot \, ) \,  , \,  F \cdot \phi \circ \Gamma  \rangle \colon  \bT_{{\rm cpct}} \longrightarrow [0, \infty) \; \, \textrm{is Borel-measurable.}
\end{equation}
We first prove that \textsl{Claim 1} implies the desired result. 
For all $\widetilde{T} \in \bT_{{\rm cpct}}$, 
set $G(\widetilde{T})= F(\widetilde{T}) \phi ( \Gamma (\widetilde{T}))$ and say that 
$G$ is $C$-Lipschitz. Let $h_0 \in (0, \infty)$, be such that $\phi (y)= 0$, for all $y\in (0, h_0)$. Observe that for all $h \in (0, h_0/2)$, 
$$ \forall \widetilde{T} \in \bT_{{\rm cpct}} \, , \quad \left|  \langle \cM_0(R_h (\widetilde{T} ) ), G  \rangle-  \langle \cM_0   (\widetilde{T} ) , G  \rangle \right| \leq Ch \, \# Z^{(h_0/2)}_{0} (\widetilde{T})  \longrightarrow 0\qquad\mbox{as $h\rightarrow 0$.} $$
By \textsl{Claim 1}, $\langle \cM_0   (\cdot ) , G  \rangle$ is measurable and a 
monotone class argument shows that $\widetilde{T} \mapsto \langle \cM_0   (\widetilde{T} ) , F\cdot 
\phi \circ \Gamma   \rangle $ is measurable for all bounded measurable 
$F\colon \bT_{{\rm cpct}} \rightarrow [0, \infty)$ and for all bounded Lipschitz 
$\phi\colon [0, \infty) \rightarrow [0, \infty)$ vanishing in a neighbourhood of 
$0$. Let $\phi_n \colon  [0, \infty) \rightarrow [ 0, \infty)$, $n \in \bN$, be a sequence of such functions 
such that $\phi_n \leq \phi_{n+1}$ and $\sup_{n \in \bN} \phi_n = \un_{(0, \infty)}$. 
By monotone convergence,  
$ \langle \cM_0   (\, \cdot \,  ) , F \rangle = \lim_{n\rightarrow \infty} 
\langle \cM_0   (\, \cdot \, ) , F \cdot \phi_n \circ \Gamma \rangle $, that is therefore 
measurable. Thus, $\cM_0$ is measurable and \textsl{Claim 1} entails Lemma \ref{measuraZht} $\rm(i)$. \cq 

\medskip

To prove \textsl{Claim 1}, we prove \underline{\textsl{Claim 2}} that is stated as follows. 
Fix $h_0 >h>0$. Fix $\phi\colon  [0, \infty) \rightarrow [0, \infty) $ a bounded Lipschitz function such that $\phi (y)= 0$, for all $y\in [0, h_0]$. 
Fix $F\colon  \bT_{{\rm cpct}} \rightarrow [0, \infty)$, a bounded Lipschitz function. Let $(T,d, \rho)$ be a compact rooted real tree. Then, 
for all $u\in [0, h)$, we set 
$$ \Psi_u (\widetilde{T})= \langle \cM_u \left( R_h (T) \right) , F\cdot \phi \circ \Gamma  \rangle \; .$$
Denote by  $(T_j, d, \sigma_j)$, $j\in J$, the subtrees of $R_h (T)$ above level $u$. Then $ \Psi_u (\widetilde{T})=\sum_{j\in J} F(\widetilde{T}_j) \phi (\Gamma (\widetilde{T}_j)) $. 
We then set 
$$ \Delta (T):= \left\{ \, d(\rho , \s) \; ; \, \sigma \in {\rm Br} (R_h (T)) \cup {\rm Lf} (R_h (T)) \,   \right\} . $$
Note that $\Delta (T)$ is a finite set. Then, we claim that 
\begin{equation}
\label{nearcont}
\tag{\textsl{Claim 2}}
\forall \widetilde{T}  \in \bT_{{\rm cpct}} \, , \; \forall u \in (0,h) \backslash \Delta (T) \, , \qquad 
\Psi_u (\widetilde{T}^\prime) 
\rightarrow \Psi_u (\widetilde{T})\qquad\mbox{as $\dgh_{{\rm cpct}} (\widetilde{T}^\prime, \widetilde{T}) \rightarrow 0$.} 
\end{equation}
Let us first prove that \textsl{Claim 2} implies \textsl{Claim 1}. 
To simplify notation, we write 
$G= F\cdot \phi \circ \Gamma $ and note that $G$ is $C$-Lipschitz. 
Suppose that $[u,v] \subseteq [0,h) \backslash \Delta (T)$. To each subtree $T_j$ 
above level $u$ in $R_h(T)$ corresponds a unique subtree of $R_h (T)$ above level $v$ that is simply the tree  
$T_j$ shortened at its root by a line of length $v-u$. Thus $|\Psi_u (\widetilde{T})-\Psi_v (\widetilde{T})| \leq C(v-u)\# \{ j\in J\colon  \Gamma (T_j) \geq h_0 \}$. 
This proves that $u \mapsto \Psi_u (\widetilde{T})$ is right-continuous on $[0,h)\backslash \Delta (T)$. 
For all $K \in (0, \infty)$ and for all $u\in (0, h)$, we set $\Phi_{u , K} (\widetilde{T})= 
\int_0^1 (K \wedge \Psi_{uv} (\widetilde{T})) dv $. Since $\Delta(T)$ is a finite set, it is Lebesgue negligible. Hence, \textsl{Claim 2} 
and dominated convergence imply that $\Phi_{u,K} \colon  \bT_{{\rm cpct}} \rightarrow [0, \infty)$ is $\dgh_{{\rm cpct}}$-continuous. Dominated convergence also implies that 
for all $\widetilde{T} \in \bT_{{\rm cpct}} $, $\lim_{u\rightarrow 0} \Phi_{u, K} (\widetilde{T})= K \wedge \Psi_0 (\widetilde{T})$. This entails that $\Psi_0 = \langle (\cM_0 \circ R_h)(\, \cdot \, )  , G  \rangle $ is Borel-measurable, which proves \textsl{Claim 1}. \cq  

\medskip

It remains to prove \textsl{Claim 2}. 
We use the previous notation $G$. We fix $u \in (0, h) \backslash \Delta (T)$. Since $\Delta (T)$ is finite, we fix $\vre \in (0, u/4)$ that can be chosen arbitrarily small and 
such that $[u-2\vre , u+2\vre] \subseteq (0, h) \backslash \Delta (T) $. Note that $\vre <h_0 /4$. Let $(T^\prime , d^\prime , \rho^\prime )$ be a compact rooted real tree such that $\dgh_{{\rm cpct}} ( R_h (T), R_h (T^\prime )) < \vre /4$. By Lemma \ref{epsiiso} there exists a pointed $\vre$-isometry $f\colon R_h(T)\rightarrow R_h(T^\prime)$. We next denote by $(T^\prime_k , d^\prime , \s^\prime_k)$, $k \in J^\prime $, the subtrees of $R_h (T^\prime )$ above level $u$, so that $\Psi_u (\widetilde{T}^\prime )= \sum_{k\in J^\prime } G(\widetilde{T}^\prime_k)$. Recall that $(T_j, d, \sigma_j)$, $j\in J$, stand for the subtrees of $R_h (T)$ above level $u$. We next set $J_\vre = \{ j\in J\colon  \Gamma (T_j) > 2\vre \}$ and we construct an {\it injective} function $\pi\colon  J_\vre \rightarrow J^\prime $ such that 
\begin{equation}
\label{septsixeps}
\forall j\in J_\vre \, , \quad \dgh_{{\rm cpct}} \left(T_j , T^\prime_{\pi (j)} \right) < 76 \vre \; .
\end{equation}
\textsl{Construction of $\pi$:} for each $j\in J_\vre$, we fix $\gamma_j\in T_j$ such that $d(\s_j, \gamma_j )= 2\vre$. Since $d(\rho, \s_j)= u$ and $\s_j\in \lgeo  \rho , \gamma_j \rgeo$, we get $d(\rho , \gamma_j)= u+ 2 \vre $ and $d^\prime (\rho^\prime, f( \gamma_j) ) > u+\vre $, and there exists $k\in J^\prime$ such that $f(\gamma_j) \in T^\prime_k$. We set $\pi (j)= k$. 

\textsl{$\pi$ is injective:} let $i\in J_\vre \backslash \{ j \}$. Then $\gamma_i \wedge \gamma_j = \s_i \wedge \s_j$ and (\ref{isobranch}) implies that $d^\prime (\rho^\prime, f(\gamma_i) \wedge f(\gamma_j) ) < d(\rho , \s_i \wedge \s_j) + 3\vre /2$. Since $[u-2\vre , u+2\vre] \subseteq (0, h) \backslash \Delta (T) $, we get 
$d(\rho , \s_i \wedge \s_j) < u-2\vre $. Consequently, $d^\prime (\rho^\prime, f(\gamma_i) \wedge f(\gamma_j) ) < u$ and $\pi (i) \neq \pi (j)$. 

We next prove that 
\begin{equation} 
\label{insideTj}
\forall j \in J_\vre \, , \; \forall \gamma \in T_j \; \textrm{ such that} \; d(\sigma_j , \gamma) \geq 2 \vre \; , \quad f(\gamma ) \in T^\prime_{\pi (j)} \; .
\end{equation}
Indeed, since  $[u-2\vre , u+2\vre] \subseteq (0, h) \backslash \Delta (T) $, $d(\rho , \gamma \wedge 
\gamma_j ) > u+ 2\vre $ and (\ref{isobranch}) entails $d^\prime(\rho^\prime , f(\gamma) \wedge f(\gamma_j) ) > u$, which implies (\ref{insideTj}).

  Next observe that for all $j\in J_\vre $, 
\begin{eqnarray}
\label{subrootdist}
d^\prime ( \s^\prime_{\pi (j)} , f(\s_j)) & \leq &  d^\prime (  \s^\prime_{\pi (j)} , f(\gamma_j)) + d^\prime ( f(\gamma_j) , f( \s_j)) \leq d^\prime ( \rho^\prime, f(\gamma_j)) -u + d( \s_j, \gamma_j) + \vre  \nonumber \\
& \leq & d( \rho , \gamma_j) + \vre -u + 2 \vre + \vre =  u+ 2 \vre +  \vre -u + 2 \vre + \vre    =6 \vre . 
\end{eqnarray}
 We next define $f_j \colon  T_j \rightarrow R_h (T^\prime)$ by setting $f_j (\s) = f(\s)$ if $d(\s_j , \s) \geq 
 2 \vre $, and $f(\s)= \s^\prime_{\pi (j)}$ if $d(\s_j , \s) < 2 \vre$. We deduce from (\ref{insideTj}) that 
$ f_j (T_j) \subseteq T^\prime_{\pi (j)}$. Next observe that if $d( \s , \s_j) < 2\vre$, then 
(\ref{subrootdist}) implies 
$$ d^\prime (f_j (\s), f(\s))= d^\prime (\s^\prime_{\pi (j)} , f(\s) ) \leq d^\prime (\s^\prime_{\pi (j)} , f( \s_j) ) + d^\prime ( f(\s_j), f(\s)) \leq 6 \vre + d(\s , \s_j) + \vre < 9 \vre . $$
Hence, $d^\prime (f_j(\s) , f (\s)) < 9 \vre $, for all $\s \in T_j$. Consequently, ${\rm dis} (f_j) \leq {\rm dis} (f) + 18 \vre < 19 \vre$. 
We next prove that $f_j (T_j)$ is a $4\vre$-net of $T^\prime_{\pi (j)}$: let $\s^\prime \in T^\prime_{\pi (j)}$ such that $d^\prime( \s^\prime_{\pi (j)} , \s^\prime) > 4 \vre$. Since $f$ is an $\varepsilon$-isometry, there exists $\s \in R_h (T)$ such that 
$d^\prime (f(\s) , \s^\prime) <\vre$. Thus,  $d( \rho ,\s) > d^\prime (\rho^\prime , f(\s) ) -\vre > d^\prime( \rho^\prime  ,  \s^\prime) -2\vre >u+2\vre$, and (\ref{insideTj}) implies that $\s \in T_j$ and $f_j (\s)= f( \s) $, which proves that $f_j (T_j)$ is a $4\vre$-net of $T^\prime_{\pi (j)}$. 
Thus, $f_j$ is a pointed $19\vre$-isometry from $(T_j , d, \s_j)$ to $(T^\prime_{\pi (j)}, d^\prime, \s^\prime_{\pi (j)} )$, which entails (\ref{septsixeps}) by Lemma \ref{epsiiso}. 

  We next prove that 
\begin{equation}
\label{hachez}
\left\{  k \in J^\prime \colon  \Gamma (T^\prime_k) \geq h_0  \right\} \subseteq \pi \left( \left\{ j \in J \colon  \Gamma (T_j) \geq h_0 /2 \right\} \right) \subseteq \pi (J_\vre ) \; .
\end{equation}
Let $\s^\prime \in T^\prime_k$ such that $d^\prime (\s^\prime_k , \s^\prime) = h_0$. There exists $\s 
\in R_h (T)$ such that $d^\prime(f(\s) , \s^\prime) < \vre$. We then get 
$$ d(\rho , \sigma ) -u> d^\prime (\rho^\prime , f(\s)) -\vre -u> d^\prime (\rho^\prime , \s^\prime) -2 \vre -u > h_0/2>2\vre. $$ 
Thus $\s \in T_j$ for a certain $j\in J_\vre$ and $f(\s ) \in T^\prime_{\pi (j)}$ by (\ref{insideTj}). Moreover, 
(\ref{branchheight}) easily entails that $d^\prime (\rho^\prime, f(\s) \wedge \s^\prime) > d^\prime (\rho^\prime , \s^\prime) -\vre  >u$. This implies that $f(\s)\in T^\prime_k$ and $\pi (j)= k$, which completes the proof of (\ref{hachez}). 

  Now recall that $G= F\cdot \phi \circ \Gamma $ and that $\phi (y)= 0$, if $y\in [0, h_0]$. Thus, (\ref{hachez}) implies that $\Psi_u (\widetilde{T})= \sum_{j\in J_\vre} G(\widetilde{T}_j)$ and $\Psi_u (\widetilde{T}^\prime)= \sum_{j\in J_\vre} 
G(\widetilde{T}^\prime_{\pi (j)})$. Recall that $G$ is $C$-Lipschitz. Then, (\ref{septsixeps}) and  (\ref{hachez}) imply 
\begin{equation}
\label{psifinal}
 \left| \Psi_u (\widetilde{T})-\Psi_u (\widetilde{T}^\prime) \right| \leq 76 C \vre \, \# \{ j\in J\colon  \Gamma (T_j) \geq h_0/2 \}. 
\end{equation}
To summarise, we have fixed 
$\widetilde{T} \in \bT_{\textrm{cpct}}$, $u \in (0, h) \backslash \Delta (T)$ and we have proved that (\ref{psifinal}) holds true 
for all sufficiently small $\varepsilon \in (0, \infty)$ and for all $\widetilde{T}^\prime \in \bT_{\textrm{cpct}}$ such that $\dgh_{\textrm{cpct}} ( R_h(\widetilde{T}), R_h(\widetilde{T}^\prime)) < \varepsilon /4$,   
Since $R_h$ is $\dgh_{{\rm cpct}}$-continuous, this entails \textsl{Claim 2} 
and the proof of Lemma \ref{measuraZht} $\rm(i)$ is complete. \cqfd

\paragraph{Proof of Lemma \ref{measuraZht} $\rm(ii)$.} Since $(\bT, \dgh)$ is a Polish space, the 
Borel Isomorphism Theorem implies that there exists a one-to-one Borel-measurable function 
$\phi\colon \bT \rightarrow \bR$ such that its inverse $\phi^{-1}\colon \bR \rightarrow \bT$ is also Borel-measurable. Let $M = \sum_{1\leq k\leq n} \delta_{\widetilde{T}_k} $ be in 
$\ccM (\bT)$, where 
$\phi(\widetilde{T}_1) \leq \cdots \leq \phi(\widetilde{T}_n)$. Then, for all $k \in \bN$, we define $\Lambda_k (M)$ as follows: 
$$ \textrm{$\Lambda_k (M)= \pnt $, if $k= 0$ or if $k>n$} \quad {\rm and} \quad 
\textrm{$\Lambda_k (M)= \widetilde{T}_k $,  if $1\leq k \leq n$.}$$
We next set 
$\ccM_f (\bT)= \{ M \in \ccM (\bT)\colon  \langle M \rangle < \infty \}$ that is clearly an element of 
the sigma-field $\cG_{\ccM (\bT)}$. 
\begin{lemma}
\label{extract}
For all $k\in \bN$, $\Lambda_k \colon  \ccM_f (\bT) \rightarrow \bT$ is measurable. 
\end{lemma}
\begin{pf} We only need to prove that $\phi\circ \Lambda_k$ is measurable. 
We set  
$A_{x, k}:=\{ M \in \ccM_f (\bT) \colon  \phi (\Lambda_k (M)) \leq x \} $ for all $x\in \bR$. Now, we observe that 
$A_{x, 0}= \emptyset $ if $\phi (\pnt) >x$, that $A_{x, 0} = \ccM_f (\bT) $, if $\phi (\pnt) \leq x$, and that 
for all $k \geq 1$, 
$$ A_{x, k}= \left( \{ \phi (\pnt) \leq x \} \cap \{ \langle M \rangle < k\} \right) \cup 
\left\{ M\in \ccM_f (\bT) \colon  \langle M \,  , \,  \un_{(\!-\infty, x]} \! \circ\!  \phi \,  \rangle \geq k  \right\} \; \in \; \cG_{\ccM (\bT)} \; ,$$
which implies the desired result
\end{pf}
Recall from (\ref{pastedef}) the definition of 
$\paste$. It is easy to check that 
$(\widetilde{T} , \widetilde{T}^\prime) \in \bT^2 \mapsto \widetilde{T} \circledast \widetilde{T}^\prime  \in \bT$ is continuous. Thus, this implies that 
\begin{equation}
\label{pastefini}
M \in \ccM_f (\bT) \; \longmapsto  \; \paste (M)= \circledast_{k\in \bN} \Lambda_k (M) \quad \textrm{is measurable.}
\end{equation}
Let $h \in (0, \infty)$. For all $M= \sum_{i\in I} \delta_{\widetilde{T}_i}$, we set $\Xi_h (M)=\sum_{i\in I} \delta_{R_h (\widetilde{T}_i)} $. Clearly $\langle  \Xi_h (M) \rangle <\infty$ and for all measurable $F\colon \bT \rightarrow [0, \infty)$, $\langle \Xi_h (M) , F \rangle = \langle M, F \circ R_h \rangle$. 
This implies that 
$\Xi_h \colon  \ccM(\bT) \rightarrow \ccM_f (\bT)$ is measurable. This result combined with (\ref{pastefini}) implies that $\paste \circ \Xi_h $ is measurable.  
Now, observe that for all $h \in (0, \infty)$ and all $M \in \ccM (\bT)$, 
$\dgh \left( \paste (M) , \paste (\Xi_h (M) ) \right) \leq h$, 
which implies Lemma \ref{measuraZht} $\rm(ii)$. \cqfd

\subsection{Proof of Lemma \ref{measDthetak}.}  
\label{measDthetakpf}
By (\ref{firstabove}), Lemma \ref{abovebelow} and Lemma \ref{measuraZht} $\rm(i)$, we only need to prove that $D$ is measurable. Let $(T,d,\rho)$ be a CLCR real tree. We first claim the following. 
\begin{equation}
\label{limRhD}
\lim_{h \to 0} D (R_h(T))= D(T) \; .
\end{equation}
Recall that $D(T)= \infty$ iff $T$ has no leaf and no branch point, namely iff $T$ is either a point tree or a finite number of half-lines pasted at their finite endpoint. In these cases, (\ref{limRhD}) obviously holds true. Let us assume that $D(T)$ is finite. First note that if $\s \in {\rm Br} (T)$, then $\s \in  {\rm Br} (R_h(T))$ for all sufficiently small $h\in (0, \infty)$. Let $\s \in {\rm Lf} (T)$ be such that 
$\lgeo \rho, \s \rgeo \cap {\rm Br} (T)= \emptyset$. Then for all $h \in (0, d(\rho, \s))$, there exists $\s^\prime \in \lgeo \rho, \s \rgeo $ such that $d(\s^\prime , \s )= h$. Thus, $\s^\prime \in {\rm Lf} (R_h(T))$ and $ d(\rho, \s)=d(\rho, \s^\prime)   +h$. Thus, for all $\s \in {\rm Br} (T) \cup {\rm Lf} (T)$, $d(\rho, \s) \geq \limsup_{h\rightarrow 0} D(R_h(T))$, which implies that $D(T) \geq \limsup_{h\rightarrow 0} D(R_h(T))$.

Conversely, 
observe that ${\rm Br} (R_h(T)) \subseteq {\rm Br} (T)$. Next, if $\s^\prime \in {\rm Lf} (R_h (T))$, then there exists $\s \in {\rm Lf} (T)$ such that $\s^\prime \in \lgeo \rho, \s \rgeo   $ and $d(\rho, \s^\prime )= d(\rho, \s )-h$. Therefore, for all $\s^\prime \in  {\rm Br} (R_h(T)) \cup {\rm Lf} (R_h(T))$, $d(\rho, \s^\prime ) \geq 
D(T) -h$. Thus, $\liminf_{h\rightarrow 0} D(R_h(T)) \geq D(T)$, which completes the proof of 
(\ref{limRhD}).

For  all $h \in (0, \infty)$, we next set 
$J_h (T)= \inf \{ a \in [0, \infty) \colon  Z_{0}^{(h)} (T)\neq Z_{a}^{(h)} (T)    \} $, with the convention that $\inf \emptyset = \infty$.  
By Lemma \ref{measuraZht} $\rm(i)$, $ Z_{a}^{(h)} $ is measurable and since $a\mapsto Z^{(h)}_a$ is caglad, the function $J_h \colon  \bT \rightarrow [0, \infty]$ is measurable.  
Now observe that $J_h (T) \geq D(R_h (T)) >0$. If $J_h (T) >D(R_h (T)) $ then, the lowest leaves of $R_h(T)$ are at the same distance from the root as the lowest branch points and there exists 
$\vre \in (0, \infty)$, such that for all $h^\prime \in (h, h + \vre ) $, 
$ D (R_{h^\prime} (T)) = J_{h^\prime} (T)= D (R_h (T))-h^\prime+h$. This implies that $\liminf_{h\in \bQ \cap (0, \infty) \to 0  } J_h (T) =\liminf_{h\in \bQ \cap (0, \infty) \to 0  } D(R_h(T))$, 
which implies the measurability of $D$ by (\ref{limRhD}). \cqfd

\subsection{Proof of Lemma \ref{edgechar}.}
\label{edgecharpf} 
Recall that for all $n \in \bN$, $\vartheta_n \colon  \bT \rightarrow \bT$ and $D_n \colon  \bT \rightarrow [0, \infty]$ are defined as $ \vartheta_{n+1} = \vartheta \circ \vartheta_n$ and $D_{n} = D \circ \vartheta_n$. We set 
$A= \{\widetilde{T} \in \bT \colon  \sum_{n \in \bN} D_n (\widetilde{T}) = \infty \}$, which is a Borel set of $\bT$ by 
Lemma \ref{measDthetak}. 

  Let $(T,d, \rho)$ be a CLCR real tree. Suppose that 
$\widetilde{T} \in \bT_{{\rm edge}}$. If $T$ has no leaf and no branch point, then 
$D_n (T)= \infty$, for all $n \in \bN$, and it belongs to $A$. 
Next assume that there is $\s \in {\rm Br} (T) \cup {\rm Lf} (T)$. 
Then, there exists 
$n_0 \in \bN$ such that $\sum_{0\leq k\leq n_0} D_k (T)= d(\rho, \s)$. If $d(\rho, \s)= \Gamma (T)$, then
for all $n>n_0$, $D_n(T)= \infty$. 
Let us assume that $\Gamma (T)= \infty$:  for all $r \in (0, \infty)$, set 
$n(r)= \# (  \overline{B}_T (\rho, r) \cap (  {\rm Br} (T) \cup {\rm Lf} (T))) $, which is finite 
since $\widetilde{T} \in \bT_{\rm edge}$; the previous arguments imply that $\sum_{0\leq k \leq n(r)}  
D_k(T)>r$, which implies that $\sum_{n \in \bN} D_n (T)= \infty$. This proves that $
\bT_{{\rm edge}} \subseteq A$. 

 Conversely, assume that $\widetilde{T} \in A $. If $D(T)= 0$, then $\vartheta_n \widetilde{T}=  \widetilde{T}$ and $D_n (T)= 0$, for all $n \in \bN$, which contradicts the assumption. If $D(T)= \infty$, then 
 $\widetilde{T}\in \bT_{{\rm edge}}$ (and recall that ${\bf k} (\widetilde{T})= 0$, by convention). 
Let us assume that 
$D(T) \in (0, \infty)$, then $\blw (T, D(T))$ is equivalent to a finite number of copies of the interval $[0, D(T)]$ pasted at $0$. If ${\bf k} (T)= \infty$, $\vartheta \widetilde{T}$ has infinitely many trees pasted at its root and the local compactness implies that there are leaves arbitrarily close to its root, which implies that $D_1 (T)= 0$; therefore $D_n (T)= 0$, for all $n \geq 1$, which contradicts the assumption. 
Thus, if $D(T) \in (0, \infty)$, then 
 ${\bf k} (T) <\infty $. These arguments and a simple recursion imply first 
that for all $n \in \bN$, $\vartheta_n \widetilde{T} \in A$ and ${\bf k} (\vartheta_n \widetilde{T}) < \infty$, and that for all $n \in \bN$ such that $R_n := D_0 (T) + \ldots + D_n (T) < \infty$, we get 
$$ \left\{ d(\rho, \s) \, ; \, \s \in {\rm Br} (T) \cup {\rm Lf} (T)\colon d(\rho, \s ) \leq R_n  \right\}= \{ R_0< R_1 < \cdots < R_n \} \; , $$ 
Thus, if $R_n < \infty$, we get 
$${\rm n} (\rho, T) + \sum {\rm n} (\s, T) \leq {\rm n} (\rho, T)+  (1+ {\bf k} (\widetilde{T}))+ (1+ {\bf k} (\vartheta_1 \widetilde{T}))+ \cdots + (1+{\bf k} (\vartheta_n \widetilde{T})) < \infty \, , $$
where the first sum is over the branch points $\s \in {\rm Br} (T)$ such that $d(\rho, \s) \leq R_n$. 
This easily entails that $(T,d,\rho)$ satisfies (\ref{degsumfin}) in the definition of real trees with edge lengths. \cqfd  

\section{Proofs of the preliminary results on GW-trees.}
\label{GWproofsec}

To prove the lemmas of Section \ref{GWsec} about GW-real trees, it is useful to be able to push forward distributions between a space of discrete combinatorial trees and $\bT_{\rm edge}$. 

\paragraph{Discrete trees with marks.} The discrete combinatorial trees that we consider are rooted, ordered and locally finite. We
use Ulam's coding (see Neveu \cite{Ne}) that allows to view such trees as subsets of the set of finite integer words 
$$\bU=\bigcup_{n\ge 0}(\bN^*)^n \; , $$
where $\bN^*$ is the set of positive integers. Here $(\bN^*)^0$ stands for $\{ \varnothing \}$, where $\varnothing$ is the empty word.  
Before recalling the formal definition of discrete trees in this context, let us set some notation: the concatenation of the two words $u=(a_1,\ldots,a_m)$ and $v=(b_1,\ldots,b_n)$ in $\bU$ is denoted by $w=u*v=(a_1,\ldots,a_m,b_1,\ldots,b_n)$. Note that $\varnothing *u=u=u*\varnothing$. A single-symbol word shall be denoted by $(j)$, where $j\in \bN^*$. 
The length of $u\in (\bN^*)^n$ is denoted by $|u|=n$, with the convention $|\varnothing |= 0$.

  For all $u\in \bU \backslash \{ \varnothing \}$, there exists $v\in \bU$ such that $u= v* (j)$ for a certain $j\in \bN^*$. Note that  $\lvert v \rvert = \lvert u \rvert -1$. We then call 
$v$ the {\it parent} of $u$ and we denote it by $\overleftarrow{u}$. We can view $\bU$ as a graph whose set of vertices is $\bU$ and whose set of edges is $ \{ \{\overleftarrow{u},u \} \, ; \, u \in \bU \backslash \{ \varnothing \} \}$, then we denote by 
$\lgeo u,v\rgeo$ the shortest path (with respect to the graph distance) between $u$ and $v$. We also set 
$ \rgeo u,v\rgeo := \lgeo u,v\rgeo \backslash \{ u\} $ and we define similarly  
$\lgeo u,v\lgeo$ and $\rgeo u,v\lgeo$. For $u,v\in\bU$, the 
last common ancestor of $u$ and $v$ is denoted by $u\wedge v$: we recall that 
$ \lgeo \varnothing ,  v \wedge u \rgeo = \lgeo \varnothing   , u \rgeo \cap  \lgeo \varnothing   , v  \rgeo$. 
\begin{definition}
\label{treeNeveudef} A non-empty subset $\ft\subset \bU$ is called a tree iff it satisfies the following conditions for all $u \in \ft$. 
\begin{enumerate}
\item[(a)] If $u\in \ft$ is different from $\varnothing$, then $\overleftarrow{u}\in \ft$. 
\item[(b)] There exists $k_u (\ft)\in\bN $, such that $\{ v\in \ft  \colon   \overleftarrow{v}=u\}=\{ u*(1), \ldots, u*(\, k_u(\ft)\, )\}$ if $k_u(\ft) \geq 1$ and $\{ v\in \ft  \colon   \overleftarrow{v}=u\}= \emptyset $ if 
$k_u(\ft)= 0$. 
\end{enumerate}
Note that (a) entails $\varnothing \in \ft$. We view $\varnothing$ as the progenitor of the population whose family tree is $\ft$. Then, $k_u (\ft)$ stands for the number of children of $u\in \ft$. 
We denote by $\bT^{{\rm discr}}$ {\it the set of all ordered rooted discrete trees}. \cq  
\end{definition}
We view $[0, \infty]$ as the compactification of $[0, \infty)$ and we denote by $\Delta$ a metric that generates this topology. 
We call $\fT= (\ft ; \bx)$ a 
{\it marked tree} if $\ft\in\bT^{\rm discr}$ and if $\bx=(x_u , u\in\ft)$, with $x_u\in[0,\infty]$, for all $u \in \ft$. We then denote by 
$\bT_{[0,\infty]}^{\rm discr}:= \bigsqcup_{\ft \in \bT^{{\rm discr}} }\left(\{ \ft \} \times [0,\infty]^{\ft}\right)$ the set of marked trees. We equip $\bT_{[0,\infty]}^{\rm discr}$ with the $\sigma$-algebra $\cG_{[0,\infty]}$ generated by the subsets 
\begin{equation}\label{elemset} A_{u,y}:=\{(\ft;\bx)\in\bT^{\rm discr}_{[0,\infty]} \colon  u\in \ft  \, , \, x_u > y \}  \, , \quad u\in\bU \, , \; y \in [0, \infty] \; .\end{equation}

\paragraph{Connection with real trees.} A discrete tree with \textit{finite} marks 
clearly corresponds to a real tree. For technical reason, we associate a real tree to $[0, \infty]$-marked discrete trees with an 
obvious restriction due to possibly infinite lifetime marks. More precisely, let 
$\fT= (\ft ; \bx) \in \bT^{{\rm discr}}_{[0, \infty]}$. 
For all $u \in \ft$, we introduce 
the following notation:  
\begin{equation}
\label{deathdef}
 \zeta_u = \sum_{v \in \lgeo \varnothing , u \rgeo} x_v \; , \quad u \in \ft . 
\end{equation} 
We can think of $\zeta_u $ as the {\it death-time of $u$} and of $\zeta_{\overleftarrow{u}}$ as 
{\it the birth-time of $u$}, with the convention $\zeta_{\overleftarrow{\varnothing}}= 0 $. For an obvious reason, we have to assume the following: 
\begin{equation} 
\label{inftypb}
\forall u \in \ft \, , \; \forall v \in \lgeo \varnothing , u \, \lgeo \, \, , \quad \zeta_v < \infty \; .  
\end{equation}
We associate with $\fT$ a rooted real tree denoted by $\Trr( \fT)=(T, d, \rho)$ 
as follows. We first set 
$$\rho =(\varnothing , 0) \quad {\rm and} \quad T= \{ \rho \} \cup \{  (u, s ) \, ; s\in (0,x_u] \cap (0, \infty), 
\, u\in \ft\} \; .$$ 
We then define a distance $d$ on $T \times T $ as follows: for all $\sigma =(u,s)\in T\setminus
\{\rho\}$, we set  $d(\rho, \sigma )=s+\sum_{v\in \lgeo \varnothing , u \lgeo } x_v $, which is finite by (\ref{inftypb}). Let $\sigma^\prime=(u^\prime,s^\prime)\in T\setminus
\{\rho\}$. We then set   
 \begin{displaymath}
d(\sigma, \sigma^\prime)= \left\{ \begin{array}{ll} 
d(\rho , \sigma)+ d(\rho , \sigma^\prime) -2 \sum_{v\in \lgeo \varnothing , u
  \wedge u^\prime \rgeo} x_v \, , & \textrm{ if $u\wedge u^\prime \notin \{ u, u^\prime  \}$,} \\
 |d(\rho , \sigma)- d(\rho , \sigma^\prime)  | \,  ,  &\textrm{otherwise.}
\end{array} \right.   
\end{displaymath}
It is easy to check that $\Trr (\fT):=(T,d,\rho)$ is a
rooted real tree. However, note that $T$ may neither be a real tree with edge lengths nor 
locally compact. We then introduce 
\begin{equation}
\label{domaintree}
 \bT^{{\rm discr}}_{{\rm edge}} =
\left\{ \fT= (\ft; \bx)\in \bT^{{\rm discr}}_{[0, \infty]}\colon \textrm{$\fT$ satisfies (\ref{inftypb}) and} \; 
 \forall a \in [0, \infty) \, , \; \# \{ u\in \ft \colon  \zeta_u \leq a \} < \infty \,   \right\} \; .
 \end{equation}
It is easy to check that if $\fT \in  \bT^{{\rm discr}}_{{\rm edge}} $, then $\Trr( \fT)$ is a real tree with edge lengths as in Definition \ref{realedgetrees}. Next observe that $  \bT^{{\rm discr}}_{{\rm edge}}$ belongs to the sigma-field $\cG_{[0, \infty]}$. Then, for all $\fT \in  \bT^{{\rm discr}}_{{\rm edge}} $, we denote by $\TR (\fT)$ the pointed isometry class of the real tree with edge lengths $\Trr (\fT)$. 
\begin{lemma}   
\label{TRmeas} $\TR\colon \bT^{{\rm discr}}_{{\rm edge}} \longrightarrow \bT_{{\rm edge}}$ is measurable.    
\end{lemma}
\begin{pf} 
Let $\ft \in  \bT^{{\rm discr}}$ be finite and set $U_{\ft}= \big\{ \bx \in [0, \infty]^\ft : (\ft; \bx) \in  \bT^{{\rm discr}}_{{\rm edge}}\big\}$, which is an open subset of $[0, \infty]^\ft$ equipped with the product topology. 
First note that $\bx \in U_\ft  \mapsto  \TR (\ft; \bx) $ is $\dgh$-continuous. For any $n \in \bN$, set 
$E_n = \bigsqcup\, ( \{ \ft \} \! \times \! [0, \infty]^\ft )$, where the disjoint union is taken over the set of discrete trees $\ft$ such that $|u| \leq n$ for all $u\in \ft$. Clearly, $E_n$ is a Polish space when it is equipped 
with the distance $d_n$ that is defined for any $\fT= (\ft ; \bx), \, \fT^\prime= (\ft^\prime ; \bx^\prime)$ in $E_n$ by 
$$ \textrm{$d_n (\fT , \fT^\prime )= 1$\quad if\quad $\ft \neq \ft^\prime$} \qquad \textrm{and} \qquad 
\textrm{$d_n (\fT , \fT^\prime )= \sum_{0\leq m<\#\ft} 2^{-m-1} (1\wedge 
\Delta (x_{u_m} , x^\prime_{u_m} ))$\quad if\quad $\ft= \ft^\prime$,} $$ 
where $u_0\!=\! \varnothing \!< \!u_1 \!<\! \cdots \!< \! u_{\# \ft -1}$ is the sequence of vertices of $\ft$ listed in the lexicographical order. 
Next, for all $\fT= (\ft ; \bx)\in \bT_{[0, \infty]}^{\rm discr}$ and all $n \in \bN$, 
we set $\ft_{| n}= \{ u\in \ft \colon  |u| \leq n\}$ and $\fT_{| n}= ( \ft_{| n} ; \bx_{| n } =( x_u, u\in \ft_{| n })) \in E_n$. 
We then define a metric $d$ on $ \bT_{[0, \infty]}^{\rm discr}$ by setting 
$  d(\fT, \fT^\prime) = \sum_{n\geq 0} 2^{-n} d_n (\fT_{| n } , \fT^\prime_{| n} ) $. 
By standard arguments, $(\bT_{[0, \infty]}^{\rm discr} , d)$ is a Polish space. The previous arguments 
entail that for any fixed $n \in \bN$, $\fT \in \bT^{{\rm discr}}_{{\rm edge}} \mapsto \TR (\fT_{|n})$ is $\dgh$-continuous. Moreover, for any fixed $\fT \in  \bT^{{\rm discr}}_{{\rm edge}}$, we easily get $\dgh (\TR (\fT), \TR(\fT_{| n}) ) \longrightarrow 0$ as $n\rightarrow \infty$. This implies that $\TR\colon   \bT^{{\rm discr}}_{{\rm edge}} \longrightarrow \bT_{{\rm edge}}$ is measurable with respect to the $d$-Borel sigma-field on $ \bT_{[0, \infty]}^{\rm discr}$, which turns out to be  $\cG_{[0, \infty]}$. 
\end{pf}

In the following lemma we prove that $\TR$ has a measurable section. This result 
is used in the proof of Lemma \ref{uniqGW}. 
\begin{lemma}
\label{section} There is $S\colon \bT_{\rm edge} \rightarrow \bT^{\rm discr}_{\rm edge}$  measurable 
such that $\TR (S(\widetilde{T}))= \widetilde{T}$, for all $\widetilde{T} \in  \bT_{\rm edge}$. 
\end{lemma}
\begin{pf} Recall Lemma \ref{extract} and its notation $\Lambda_k$. First note that for all 
$\widetilde{T} \!\in \! \bT_{\rm edge}$, $\langle \cM_0 (\widetilde{T} )\rangle \!= \!
Z^+_0 (\widetilde{T} ) \!< \!\infty$. Lemma \ref{extract} and Lemma \ref{measuraZht} (i) allow to define 
for all $k \!\in\! \bN^*$ a measurable function $\Phi_k\colon \bT_{\rm edge} \rightarrow  \bT_{\rm edge} $ by setting $\Phi_k (\widetilde{T} )=\Lambda_k (\cM_0 (\widetilde{T} ))$.
Then, for all words $u\in \bU$ we recursively define 
$\phi_u\colon  \bT_{\rm edge} \rightarrow \bT_{\rm edge} $ by setting for all $k \in \bN^*$, 
$\phi_{(k)}(\widetilde{T} ) = \Phi_k (\widetilde{T} )$, and $\phi_{u*(k)} = 
\Lambda_k(  \cM_0 ( \vartheta( \phi_u(\widetilde{T}))))$. 

If $ Z^+_0 (\widetilde{T} )\!=\!0$, then $\widetilde{T}\!=\!\pnt$ and we set $S(\widetilde{T} )\!=\!\{(\varnothing,0)\}$ that is the progenitor with zero lifetime and no children. Let us assume that $n:=Z^+_0 (\widetilde{T} ) \geq 1$. Then, for all $k\in \{ 1, \ldots, n\}$, we set $ \ft_k = \{ u\in \bU  \colon   \phi_{(k)*u} (\widetilde{T} ) \neq \pnt\}$ and for all 
$u \in \ft_k$, we set $x_u (k)= D(\phi_{(k)*u} (\widetilde{T} ) ) $.   
It is easy to check that $\fT_k :=(\ft_k; \bx(k)) \in \bT^{\rm discr}_{\rm edge}$. Then, we set $S(\widetilde{T} )=\{((k)*u,x_u(k));u\in\ft_k,k\in\{1,\ldots,n\}\}\in \bT^{\rm discr}_{\rm edge}$ and we easily check that 
$\TR (S(\widetilde{T}))= \widetilde{T}$. Also, $S$ is clearly measurable since $\phi_u$ and $D$ are measurable. 
\end{pf}

\subsection{Proof of Lemma \ref{uniqGW}.}
\label{uniqGWpf}

Lemma \ref{uniqGW} looks obvious, but it is not. Since it is the
point of entrance of GW-laws into the space $\bT$, we proceed with care in several steps.  

\medskip

\noi
\textbf{Step 1: existence.} Here we use a construction in $\bT^{{\rm discr}}_{[0, \infty)}$. 
For all $u\in \bU$, we define the {\it $u$-shift} $\theta_u$, by setting for all $w= u*v$, $\theta_u w=v$. For every subset $A \subseteq \bU$, we also define $\theta_u A $ as the (possibly empty) set of words $v\in \bU$ such that $u*v \in A$. For all $u\in\ft$, we set $\theta_u\fT=(\theta_u\ft;\theta_u \bx)$, where 
$\theta_u \bx=(x_{u*v},v\in\theta_u\ft) $, and we slightly abuse notation by writing $u\in\fT$ instead of $u\in\ft$ and $k_u(\fT)$ instead of $k_u(\ft)$.

 Let us fix an offspring distribution $\xi$ and $c\in (0, \infty)$. Let $(\Omega, \cG, \bP)$ be a probability on which is defined  a family $( N(u), \bx_u)_{u\in \bU}$ 
of i.i.d.~$\bN\! \times [0, \infty)$-valued r.v.~with law $\xi \otimes ce^{-cx}dx$. We then set 
$$ \ftau= \{ \varnothing \} \cup \bigcup_{n\geq 1} \big\{ u\!=\! (a_1, \ldots , a_n)\!\in \!(\bN^*)^n\colon  \forall k\! \in \! \{ 1, \ldots , n\}, \; \, a_k \! \leq \! N_{(a_1, \ldots, a_{k-1}) }  \big\}   \, , $$
with the convention that $(a_1, \ldots, a_{k-1})\!=\!\varnothing$ if $k\!=\!1$. We set $\bx\!=\!(x_u, u\in \ftau)$ and $\fTau\!=\!(\ftau; \bx)$. Clearly, $\fTau\!\in\!\bT^{{\rm discr}}_{[0, \infty]} $. 
Recall from (\ref{elemset}) the definition of the elementary sets $A_{u, y}$. We immediately see 
that $\{ \fTau \in A_{u,y} \}\!\in\!\cG $, which entail that $\fTau\colon \Omega \rightarrow  \bT^{{\rm discr}}_{[0, \infty]}$ is $(\cG, \cG_{[0, \infty]})$-measurable. Moreover, we easily check that $\fTau$ satisfies the following two properties.  
\begin{enumerate}
\item[(a)] The law of  $k_\varnothing(\ftau)$ is $\xi$, the law of $\bbx_\varnothing$ is $\eta(dx)=ce^{-cx}dx$, and $k_\varnothing(\ftau)$ and $\bbx_\varnothing $ are independent. 
\item[(b)] For all $k \in \bN^*$ such that $\xi (k) >0$, the subtrees $(\, \theta_{(j)}(\fTau) ; 1\le j\le k\, )$ under $\bP ( \, \cdot \,  | \, k_\varnothing (\ftau)=k)$ are i.i.d.\ copies of $\fTau$ under $\bP$, and they are independent of $\bbx_\varnothing$. 
\end{enumerate}
Furthermore, if $\fTau^\prime$ also satisfies (a) and (b), we check that for all 
$n \! \geq \!  1$, all $u_1, \ldots , u_n \in \bU$, all $y_1, \ldots , y_n\in [0, \infty]$, 
$ \bP ( \fTau \in A_{u_1, y_1} \cap \cdots \cap A_{u_n, y_n })\! =\!  \bP ( \fTau^\prime  \in A_{u_1, y_1} \cap \cdots \cap A_{u_n, y_n })$ and a monotone class argument entails that $\fTau$ and $\fTau^\prime$ have the same law, 
which we call the GW($\xi, c$)-distribution on
$\bT^{\rm discr}_{[0,\infty]}$ and which is therefore characterised by (a) and (b). We refer to Neveu \cite{Ne}, for more details.  

For any $a \in [0, \infty)$, we then set $Z^+_a (\fTau)= \# \{ u\in \ftau : \zeta_{\overleftarrow{u}} \! \leq \! a \! <\! \zeta_u \} \in \bN \cup \{ \infty\}$, that is the number of individuals that are alive at time $a$. Properties (a) and (b) imply that up to a possible explosion in finite time, the process $a\mapsto Z^+_a (\fTau)$ is a continuous-time $\bN$-valued Markov chain whose matrix-generator $(q_{i,j})_{i,j\in \bN}$ is given by 
$q_{i,i} =-ci$, $q_{i,j}= 0$ if $j<i-1$, $q_{i, i-1}= c i\xi (0)$ and $q_{i,j}= c i\xi (j-i+1)$ if $j >i$. Standard analytical computations imply that a.s.~explosion does not occur iff $\xi$ is conservative as defined in (\ref{nonexpdiscr}): see \cite[Section III.3]{AthNey} for more details. Then, if $\xi$ is conservative, $a\! \in [0, \infty) \mapsto Z^+_a (\fTau)$ is $\bN$-valued and cadlag a.s. Thus, 
\begin{equation}
\label{concons}
 \textrm{$\xi$ conservative} \quad \Longrightarrow \quad \textrm{$\bP$-a.s.} \quad \fTau \in \bT_{{\rm edge}}^{{\rm discr}} \; .
\end{equation}   
Let us furthermore assume that $\xi$ is proper. Then, 
the law on $\bT$ of $\TR(\fTau)$ satisfies (a) in Definition \ref{GWrealtreedef}. 
Note that this law is concentrated on $\bT_{\rm edge}$. This proves that for every proper conservative offspring distribution $\xi$ and every $c\in (0, \infty)$, there exists at least one probability measure on 
$\bT_{\rm edge}$ that satisfies (\ref{defiGWcalcul}) in Definition \ref{GWrealtreedef}.  \cq 

\medskip

\noi
\textbf{Step 2.} Let $Q$ be as in Definition \ref{GWrealtreedef}. We claim that $Q(\bT_{\rm edge} )= 1$. 

\smallskip

\noi
\textsl{Proof.} On the auxiliary probability space $(\Omega, \cG, \bP)$, we consider   
an $\bN$-valued Markov process $Z\!=\!(Z_t)_{t\in [0, \infty)}$ with initial state $Z_0\!=\! 1$ 
and with matrix-generator $(q_{i,j})_{i,j\in \bN}$ as defined above.  
We then set $J_0\!=\! \inf \{ t\!>\!0\colon  Z_t \!\neq \!Z_{0}\}$ and for all $n\! \in \!\bN$, 
$J_{n+1}\!=\! \inf \{ t \!>\! J_n\colon  Z_{t} \!\neq \!Z_{J_n} \}$, with the convention that $\inf \emptyset = \infty$. Then, $(J_n)_{n \in \bN}$ are the jump times of $Z$ and if $Z$ is absorbed at time $J_n$, $J_p= \infty$, for all $p>n$. Since $\xi$ is assumed to be conservative, $J_n \rightarrow \infty$ a.s.~as $n\rightarrow \infty$. 

Let $Q$ be a law on $\bT$ as in Definition \ref{GWrealtreedef}. It is easy to prove recursively that 
$D_0+ \cdots + D_n$ under $Q$ has the same law as $J_n$ under $\bP$. This implies that 
$Q$-a.s.\ $\sum_{n\geq 0} D_n = \infty$, which implies the claim by Lemma \ref{edgechar}. \cq

\medskip

\noi
\textbf{Step 3.} Let $Q$ be as in Definition \ref{GWrealtreedef} and suppose that $Q$ satisfies (\ref{defiGWcalcul}) 
with $(\xi, c)$ and $(\xi^\prime, c^\prime)$. Then $(\xi, c)=(\xi^\prime, c^\prime)$ follows straight from (\ref{defiGWcalcul}). \cq

\medskip

\noi
\textbf{Step 4.} Conversely: let $\xi$ be proper and conservative, let $c\in (0, \infty)$. Suppose that $Q$ and $Q^\prime$ satisfy (\ref{defiGWcalcul}) in Definition \ref{GWrealtreedef}. Then, we claim that $Q= Q^\prime$.   

\smallskip

\noi
\textit{Proof.} Recall the function $S$ from Lemma \ref{section}. In the definition of $S$, 
the vertices have been ordered in a way that causes a lack of exchangeability. This is why we introduce a shuffling kernel as follows. Let $\fT \in  \bT^{\rm discr}_{\rm edge}$. Denote by $K(\fT, d\fT^\prime)$ the law of the 
discrete marked tree obtained by permuting independently and uniformly the siblings (with their corresponding lifetime marks). It is easy to check that $K$ is a measurable kernel. 
Then by Lemma \ref{section}, it is easy to check that the two laws 
$P:=\int Q (d \widetilde{T}) K( S(\widetilde{T}), d\fT )$ and $P^\prime:=\int
Q^\prime (d \widetilde{T}) K( S(\widetilde{T}), d\fT )$ satisfy (a) and (b) of the definition of a discrete Galton-Watson distribution. As already mentioned there is a unique 
GW($\xi, c$) law on $\bT_{[0, \infty]}^{\rm discr}$. Therefore, $P= P^\prime$. 
But now observe that $Q$ is the law of $\TR$ under $P$ and that $Q^\prime$ is the law of $\TR$ under $P^\prime$. Thus, $Q= Q^\prime$, which entails the desired result. This finishes the proof of Lemma \ref{uniqGW}. \cqfd

\subsection{Proof of Lemma \ref{regenene}.}
\label{regenenepf}

\paragraph{Basic computations.} Before proving Lemma \ref{regenene}, let us prove some basic facts. Let us fix a proper conservative offspring distribution $\xi$ and $c\in (0, \infty)$. Recall that $\bT^{\rm discr}$ stands for the set of (ordered rooted) discrete trees with no mark as in Definition \ref{treeNeveudef}. We then set 
$\bT^{{\rm discr}}_f= \{ \ft \in \bT^{\rm discr}\colon \# \ft < \infty\}$. Let $\ft \in \bT^{{\rm discr}}_f$. 
Recall that for all $u\in \ft$, $k_u (\ft)$ stands for the number of children of $u$.  We denote by ${\rm Lf} (\ft)= \{ u \in \ft \colon  k_u(\ft) = 0\}$ the set of leaves of $\ft$.  
For each subset $ S \subseteq {\rm Lf} (\ft)$, we define the following weight 
$ w_\xi (\ft,S)= \prod_{u \in \ft \backslash S} \xi \left( k_u (\ft) \right)$. 
Let $\bx = (x_u)_{u\in\ft} \in [0,\infty)^{\ft } $, so that $\fT= (\ft ; \bx)$ is a $[0,\infty)$-marked tree. Recall from (\ref{deathdef}) notation $\zeta_u$ and $\zeta_{\overset{\leftarrow}{u}}$ for resp.\ the death time and the birth time of $u\in \ft$.  We also set $L (\bx)= \sum_{u\in \ft} x_u$.  
For all $a \in (0, \infty) $, we define 
$$D_{\ft,S,a}= \left\{ \bx = (x_u)_{u\in \ft} \in [0,\infty)^{\ft}\colon  \textrm{$\zeta_u < a$ if $u\in  \ft \backslash S$ and $x_u = a -\zeta_{\overset{\leftarrow}{u}}$ if $u\in S$} \right\} . $$
We also introduce the following finite measure on $[0, \infty)^{\ft}$:
$$M^{ c}_{\ft, S, a} (d\bx ):= \un_{D_{\ft,S,a}} (\bx) c^{\# \ft -\# S}
 e^{-cL(\bx) } \prod_{u\in S}\delta_{a-\zeta_{\overset{\leftarrow}{u}}}(dx_u)\prod_{u\in \ft \backslash S} dx_u \; , $$
 where $\delta_b (dy)$ stands for the Dirac mass at $b \in [0, \infty)$. 
 
 \medskip
 
  Let $\fTau\!=\! (\ftau ; \bbx)$ have the GW($\xi, c$)-distribution on $\bT^{\rm discr}_{[0,\infty]}$, as defined in the Section \ref{uniqGWpf} (and recall that $\fTau$ satisfies (a) and (b)). 
We then set ${\bf S}\!=\! \{ u \!\in\! \ftau\colon  \zeta_{\overset{\leftarrow}{u}} \!<\! a \!\leq \! \zeta_u \}$. 
We list ${\bf S}$ in the lexicographical order and write ${\bf S} \!= \!\{ u_1 \!< \! \cdots \!<\! u_{\# {\bf S}} \} $. 
The forest of discrete trees above level $a$ is then given by 
$\fFau^{a}= \left( \fTau_\ell \right)_{1 \leq \ell \leq \# {\bf S} }$, where $ \fTau_\ell= \left( 
\theta_{u_\ell} \ftau \, ; \,  \bx^\prime \!=\! (x^\prime_v , v\! \in \! \theta_{u_\ell} \ftau ) \right)$ with $x^\prime_\varnothing \!=\! \zeta_{u_\ell} \! -\! a$ and $x^\prime_v = x_{u_\ell * v}$, for any $v\in \theta_{u_\ell} \ftau$ distinct from $\varnothing$. Here, we use the convention that $\fFau^a$ is a cemetery point 
$\partial $ if ${\bf S}\! =\! \emptyset$.
The tree below level $a$ is then given by $\fTau_a\! =\!  (\ftau_a ; \bbx^a)$, where $ \ftau_{a } \!= \! \{ u \! \in \ftau \colon  \zeta_{\overset{\leftarrow}{u}} \! < \! a \}$, and where $\bbx^a_u\!= \!\bbx_u $, if 
$u \! \in\!  \ftau_{a} \backslash {\bf S}$ and $\bbx^a_u \!=\! a\!-\!\zeta_{\overset{\leftarrow}{u}} $ if $u\! \in\!  {\bf S}$. 
We next denote by $\Pi_{n} (d\fF )$ the law of a finite sequence (namely, a forest) of 
$n$ independent GW($\xi, c$)-discrete trees (with the convention that $\Pi_{0}$ is the Dirac mass on the cemetery point $\partial$). For all measurable functions $G_1\colon \{ \partial \} \sqcup \bigsqcup_{n\geq 1}  
(\bT^{{\rm discr}}_{[0, \infty]})^n \rightarrow [0, \infty)$ and  
$G_2\colon \bT^{{\rm discr}}_{[0, \infty]} \rightarrow [0, \infty)$, 
we easily get for all $\ft \in \bT^{{\rm discr}}_f$, for all $S \subseteq {\rm Lf} (\ft)$ and for all $a\in (0, \infty)$, 
\begin{equation}
\label{fondameasy}
\bE \left[  G_1 (\fFau^a ) G_2 (\fTau_{a} ) \un_{\{ \ftau_a = \ft ;\, {\bf S}= S \}}\right]
=  w_\xi (\ft, S) \,  \Pi_{\# S} [ G_1] 
\int_{D_{\ft, S, a}}  
M^{ c}_{\ft, S, a} (d\bx) \, G_2 ( \ft ; \bx ) \;  .
\end{equation}
We also set 
\begin{equation}
\label{pxitdef}
 p^a_{\xi ,c} (\ft) := \bE [ \un_{\{ \ftau_{a} = \ft \}}] = \sum_{S \subseteq {\rm Lf} (\ft)} w_{\xi} (\ft, S)\langle  
 M^{ c}_{\ft, S, a} \rangle \, ,  
\end{equation}
where $\langle  M^{c}_{\ft, S, a} \rangle $ stands for the (finite) mass of the measure 
$  M^{c}_{\ft, S, a}  $. Recall from (\ref{concons}) that since $\xi$ is proper and conservative, $\fTau \in 
\bT_{{\rm edge}}^{{\rm discr}}$ a.s.~and thus, 
\begin{equation}
\label{sumpxit}
\sum_{\ft \in  \bT^{{\rm discr}}_f }  p^a_{\xi, c} (\ft) = 1\; .
\end{equation}
Now set $\widetilde{\cT}= \TR (\fTau)$ that is a GW($\xi, c$)-real tree, whose law on $\bT$ 
is  
$Q_{\xi, c}$. 
Note that $\abv (a, \widetilde{\cT})= \TR(\fFau^a )$, that 
$\blw (a, \widetilde{\cT})=\TR (\fTau_{a})$, and that $Z^+_a (\widetilde{\cT})= \#{\bf S}$. 
Thus, (\ref{fondameasy}) implies that for all measurable functions $F, G \colon  \bT \rightarrow \bR_+$, for all $n \in \bN$ and all $a \in (0, \infty) $, 
\begin{equation}
\label{realfondameasy}
Q_{\xi , c} \left[ F\left( \abv (a, \cdot )\right) \, G \left( \blw (a, \cdot) \right) \, \un_{\{ Z^+_a = n\} } \right]= Q_{\xi , c} \left[ Q^{\circledast n}_{\xi ,c} [F]  \, G \left( \blw (a, \cdot) \right) \, \un_{\{ Z^+_a = n\} }  \right] .
\end{equation}
Note that for any $a\!\in\! [0, \infty)$, there is no $u\!\in \!\ftau$, 
such that $a\!=\! \zeta_u$, a.s. Thus, $b\mapsto (Z^+_b (\fTau);  \blw(b, \TR(\fTau)))$ 
is continuous at time $a$, a.s. Namely, 
\begin{equation}
\label{stochcontQ}
\forall a\in [0, \infty), \; \textrm{$Q_{\xi,c}$-a.s.} \quad 
b\mapsto\big( Z^+_b,  \blw (b, \, \cdot\, )\big) \; \textrm{is continuous at time $a$}.
\end{equation}
This implies in particular that (\ref{realfondameasy}) holds true with $a= 0$. 

\paragraph{Proof of Lemma \ref{regenene}.} Note that (\ref{realfondameasy}) and 
(\ref{stochcontQ}) prove that $\rm(iii) \Longrightarrow \rm(ii)$. 
The implication $\rm(ii) \Longrightarrow \rm(i)$ is obvious. 
It only remains to prove $\rm(i) \Longrightarrow \rm(iii)$. 
So, we assume that $Q$ is a probability measure on $\bT$ that satisfies $\rm(i)$, and that is such that 
$Q(0< D< \infty) >0$ and such that 
$Z^+$ is $Q$-a.s.\ cadlag.

  We first prove that $Q(Z^+_0= 1)= 1$. 
Indeed, by $\rm(i)$ with $a= 0$, $Q(Z^+_0= n)=  Q(Z^+_0= n)Q^{\circledast n}(Z^+_0= n)$, for all $n \geq 1$. 
Suppose that there exists $n \geq 2$ such that $Q(Z^+_0 = n) >0$, then 
$Q^{\circledast n}(Z^+_0= n)= 1$, which implies that $ Q^{\circledast n}(Z^+_0= k)= 0$ if $k\neq n$. 
But $Q^{\circledast n}(Z^+_0= n^2) \geq 
 Q(Z^+_0= n)^n >0$, which is impossible. 
Thus, $Q(Z^+_0= 0)+ Q(Z^+_0= 1)= 1$. 
Since $Q(D = \infty) < 1$,  $Q( Z^+_0= 0) < 1$, which implies that $Q(Z^+_0= 1) >0$. 
Now observe that $\rm(i)$ entails 
$Q(Z^+_0= 1) =Q(Z^+_0= 1) ^2$, so we get $Q(Z^+_0= 1) = 1$.

This implies that for all $a\in [0, D)$, $Z^+_a = 1\neq Z^+_D $. 
Namely, $D=  \inf \{ a\in [0, \infty) \colon  Z^+_a \neq Z^+_0\}$.
We next fix $b\in (0, \infty)$, $p\in \bN$, and $G\colon  \bT \rightarrow \bR$, bounded and measurable. 
By $\rm(i)$ with $a=2^{-p}b$, we get  
$$Q \! \left[  G( \abv (b, \cdot)) \un_{ 
\bigcap_{{1\leq k \leq 2^p}} \{  Z^+_{k2^{-p}b} = 1\} } \right] \!= \! Q \! \left(  Z^+_{2^{-p}b} \! = 1 \right)  
Q \! \left[ G( \abv (b\!-\! b2^{-p}, \cdot)) \un_{ 
\bigcap_{{1\leq k \leq 2^p-1}} \{ Z^+_{k2^{-p}b} = 1\} } \right].
$$
An easy inductive argument implies that 
$$ Q \! \left[ G( \abv (b, \cdot)) \un_{ 
\bigcap_{{1\leq k \leq 2^p}} \{ Z^+_{k2^{-p}b} = 1\} } \right] \! =\!  Q[G] Q \!\left( Z^+_{b2^{-p}} \!= 1 \right)^{2^p } \!\! = \!  Q[G] Q \! \left(  \bigcap_{^{1\leq k \leq 2^p}} \{Z^+_{k2^{-p}b} \!= \!  1\} \right) . $$
Since $a\mapsto Z_a^+$ is $Q$-a.s.\ cadlag, $\lim_{p\rightarrow\infty} \un_{ \bigcap_{{1\leq k \leq 2^p}} \{  Z^+_{k2^{-p}b} = 1\} } = \un_{\{ D >b\}}$ and a simple 
argument and the assumption $Q(0\!<\! D\!<\! \infty) \! >\!0$ imply that there is $c\in (0, \infty)$ such that $Q(D>b) =e^{-cb}$. Namely, 
$D$ under $Q$ has an exponential law with mean $1/c$ and we get  
\begin{equation}
\label{oublii}
Q\left[ G( \abv (b, \cdot)) \un_{\{ D > b \} } \right]= Q[G] \, Q( D >b) \; .
\end{equation}  
Recall that $\lceil \cdot \rceil$ stands for the ceiling function and set 
$D_p= 2^{-p} \lceil 2^p D \rceil $ that decreases to $D$ as $p\rightarrow \infty$. 
Fix $n\in \bN \backslash \{ 1\}$, $F\colon  \bT \rightarrow \bR $ and $f\colon  [0, \infty) \rightarrow \bR$, bounded and continuous. We then set 
$$ A_{k,p}\! =\!  Q \big[F (\abv (k2^{-p} , \cdot )) f(k2^{-p}) \un_{\{ D_p =2^{-p} k ; Z^{+}_{D_p} = n \}} 
\big]  $$
and
$$ A_p = \sum_{{k \in \bN}}  A_{k,p} \; =   Q  \big[ 
F( \abv (D_p , \cdot) ) f(D_p) \un_{\{ Z^+_{D_p} = n\}}\big] . 
$$
Then, $\lim_{p\rightarrow\infty} A_p= Q\left[ F( \vartheta) f( D)\un_{\{ {\bf k} = n \}} \right]$, by Lemma \ref{abovebelow} and since $Z^+$ is assumed to be $Q$-a.s.\ cadlag. 
We next apply (\ref{oublii}) with $b= (k-1)2^{-p}$, and then $\rm(i)$ with $a=2^{-p}$ to get   
$$  A_{k,p} = Q^{\circledast n} [F]  Q\left( Z^+_{2^{-p}} = n \right)  e^{-c(k-1)2^{-p}}   f(k2^{-p}).$$
By taking $F$ and $f$ equal to $1$, we get $Q(  D_p =2^{-p} k ; Z^{+}_{D_p} = n )= Q\left( Z^+_{2^{-p}} = n \right)  e^{-c(k-1)2^{-p}} $. Summing over $k \geq 1$ entails 
$Q\left( Z^+_{2^{-p}} = n \right)= (1-e^{-c2^{-p}})Q(   Z^{+}_{D_p} = n ) $. Thus 
 \begin{eqnarray*}
 A_{k,p} & =&  Q^{\circledast n} \big[ F \big]  Q\big( Z^+_{D_p} = n \big) (1-e^{-c2^{-p}}) e^{-c(k-1)2^{-p}}   f(k2^{-p}) \\
& =& Q^{\circledast n} \big[ F \big]  Q\big( Z^+_{D_p} = n \big)   Q\left( D_p = k2^{-p} \right)
f(k2^{-p}) .
 \end{eqnarray*}
Summing over $k\! \geq \!1$ entails $A_p \!=\! Q^{\circledast n} [F]  Q( Z^+_{D_p} \!=\! n) Q[ f(D_p)] $, which implies $Q[ F( \vartheta) f( D)\un_{\{ {\bf k} = n \}}] = 
Q^{\circledast n} [F] \xi (n) \int_{0}^\infty f(x) ce^{-cx}dx $, 
where $\xi (n)= Q (Z^+_{D}= n)$, for all $n\in \bN$. Namely, 
$Q$ is the law of a GW($\xi , c$)-real tree.  \cqfd

\subsection{Proof of Lemma \ref{poubelle}.}
\label{poubellepf}
The statement for GW-forests is easily derived from the analogous result for single GW-real trees. We only need to prove that for all $a \geq b \geq 0$, all $p\in \bN$ and all bounded measurable functions $F\colon \bT \rightarrow \bR$:
\begin{equation}
\label{limibelow}
\lim_{n\rightarrow\infty}  Q_{\xi_n , c_n} \left[ F \left( \blw (a, \, \cdot \,  \right)  \right]= 
Q_{\xi_\infty , c_\infty} \left[ F \left( \blw (a, \, \cdot \,  \right)  \right]\; .
\end{equation}
Recall (\ref{pxitdef}) and (\ref{sumpxit}). Thus, for all $n \in \bN \cup \{ \infty \}$, we have 
$$ B_n:=  Q_{\xi_n , c_n} [ F ( \blw (a, \, \cdot \,  )]= \sum_{\ft \in \bT^{{\rm discr}}_f } 
\sum_{S \subseteq {\rm Lf} (\ft)} w_{\xi_n} (\ft, S) \int_{D_{\ft, S, a}}  
M^{c_n}_{\ft, S, a} (d\bx) \, F \left( \TR( \ft ; \bx ) \right) \;  . $$
Observe that for all $\ft \in  \bT^{{\rm discr}}_f $ and all $S \subseteq {\rm Lf} (\ft) $, 
$\lim_{n \to \infty} w_{\xi_n} (\ft, S) = w_{\xi_\infty} (\ft, S) $, 
\begin{equation}
\label{pexilimi}
\lim_{n \to \infty} \int M^{ c_n}_{\ft, S, a} (d\bx) \, 
F \left( \TR( \ft ; \bx ) \right)=  \int M^{ c_\infty}_{\ft, S, a} (d\bx) \, 
F \left( \TR ( \ft ; \bx ) \right)
\end{equation}
and $\lim_{n \to \infty} p^a_{\xi_n ,c_n} (\ft)=  p^a_{\xi_\infty , c_\infty} (\ft )$. 
Thus, for all \textsl{finite} subsets  $A \subset  \bT^{{\rm discr}}_f $, we get 
$$\limsup_{n \rightarrow \infty}\,  \lvert B_n -B_\infty \rvert 
\leq  
\lVert F \rVert 
\, p^a_{\xi_\infty , c_\infty} \!\left( \bT^{{\rm discr}}_f \backslash A \right) 
+
 \lVert F \rVert \lim_{n\rightarrow \infty} 
p^a_{\xi_n, c_n} \,  \left( \bT^{{\rm discr}}_f \backslash A \right) 
=
2 \lVert F \rVert \left(1\!-\! p^a_{\xi_\infty, c_\infty} \! (A) \right) , $$
which implies (\ref{limibelow}) because $p^a_{\xi_\infty, c_\infty}$ is a probability function on the countable 
set $\bT^{{\rm discr}}_f$. 
\cqfd

\subsection{Proof of Lemma \ref{measAred}.} 
\label{RAproofsec}
Let $A \subseteq \bT$ be hereditary and let $h \in (0, \infty)$. We set $A_h= A\cap \{ \Gamma \geq h\}$. Note that $A_h$ is hereditary. Let $\phi$ stand for 
a Borel isomorphism from $\bT$ onto $\bR$. Namely, $\phi\colon \bT \rightarrow \bR$ is one-to-one and $\phi$ as well as $\phi^{-1}\colon \bR \rightarrow \bT$ are Borel-measurable. 
For all $k \in \bN^*$, let $L_k\colon \bT \rightarrow \bT$ be defined as follows.  
For all $\widetilde{T} \in \bT$, set $\cM_h (\widetilde{T})= \sum_{i\in I} \delta_{\widetilde{T}_i}$ and 
$n=\langle \cM_h (\widetilde{T}) , \un_{A_h} \rangle$; then, for all $k >n$, $L_k (\widetilde{T})= \pnt$ and if $n \geq 1$, 
$$ \sum_{i\in I} \un_{A_h} (\widetilde{T}_i)\delta_{\widetilde{T}_i}= \sum_{1\leq k\leq n} \delta_{L_k (\widetilde{T})} \quad {\rm with} \quad \phi (L_1 (\widetilde{T})) 
\leq \cdots  \leq \phi (L_n (\widetilde{T})) \; .$$
$L_k (\pnt)= \pnt$.  
We argue as in Lemma \ref{extract} to prove that $L_k$ is measurable.

Then for every word $u\in \bU$, we define a measurable function 
$L_u\colon \bT \rightarrow \bT$ such that $ L_\varnothing$ is the identity map on $\bT$, $L_{(k)}= L_k$, for all $k\in \bN^*$ and $L_v \circ L_u = L_{u*v}$, for all $u,v\in \bU$. 

We next fix a CLCR real tree $(T, d, \rho)$ and we define a $[0, \infty)$-marked discrete tree 
$(\ft; \bx) \in \bT^{{\rm discr}}_{[0, \infty)}$ by setting 
$$ \ft= \{\varnothing \} \cup \left\{ u \in \bU \backslash \{ \varnothing  \} \colon  L_u (\widetilde{T}) \neq \pnt \right\} \quad 
{\rm and}  \quad \forall u \in \ft \, , \quad x_u = h \; .$$
We then set $(T^\prime , d^\prime, \rho^\prime) = {\Trr} (\ft; \bx)$. Recall that each 
edge of $T^\prime$ has length $h$ and corresponds to one vertex in $\ft$. Recall that formally, 
$\rho^\prime = (\varnothing , 0 )$ and 
$T^\prime= \{ \rho^\prime \} \cup \{ (u, s) ; s \in (0, h] , \, u \in \ft \}$. 
We denote by 
$\Phi_h (\widetilde{T})$ the pointed isometry class of $T^\prime$. Recall from Lemma \ref{TRmeas} that 
$\TR$ is measurable. Since the functions $L_u$, $u\in \bU$, are measurable, 
$\Phi_h\colon  \bT \rightarrow \bT$ is then measurable.

     We now prove that $T^\prime$ and $R_A(T)$ are close with respect to $\dgh$. To that end, for all $\ell \in \bN$, we set 
$$ S_\ell = \{ \rho \} \cup \left\{ \s \in T \colon  d(\rho, \s)= \ell h \quad {\rm and} \quad 
\langle \cM_0 (\widetilde{\theta}_\s T) , \un_{A_h} \rangle \geq 1 \right\} \quad {\rm and} \quad S:= \bigcup_{\ell \in \bN} S_\ell \; .$$
It is easy to see that $S$ is a $(2h)$-net of $R_A(T)$. Indeed, let $\s \in R_A (T)$. There exists $\ell \in \bN$ such that $\ell h \leq d(\rho, \s ) < (\ell +1) h$. If $\ell = 0$, then $d(\rho, \s) \leq h$, which entails $d(\s , S) \leq h$. Assume that $\ell \geq 1$ and let $\s^\prime \in \lgeo \rho, \s \rgeo$ be such that $d(\rho, \s^\prime )= (\ell-1) h$. Then, $d(\s, \s^\prime) \in [h, 2h)$. Denote by $T_*^\circ$ the connected component of $\theta_{\s^\prime} T \backslash \{ \s^\prime \} $ that contains $\s$ and set $T_*= T^\circ_* \cup \{ \s^\prime \}$. Note that $T_*$ is an atom of $\cM_0 (\widetilde{\theta}_{\s^\prime} T)$ and that $\Gamma (T_*) \geq h$. Since $\s \in T_*^\circ \cap R_A (T)$, there exists $\s^{\prime \prime} \in T^\circ_*$ such that $\widetilde{\theta}_{\s^{\prime \prime}} T \in A$. Since $\theta_{\s^{\prime \prime}} T= \theta_{\s^{\prime \prime}} T_* $, we then get $\widetilde{T}_* \in A$. This implies that $\s^\prime \in 
S_{\ell-1}$. Thus, $d(\s , S) \leq 2h$, which proves that $S$ is a $(2h)$-net of $R_A(T)$.

  From the definition of the functions $L_u$ and of the tree $\ft$, we easily check that there is 
a function $\jmath\colon  \ft \rightarrow S$ that satisfies the following property: $\jmath$ is surjective, $\jmath (\varnothing) = \rho$, and for all $u\in \ft \backslash \{ \varnothing \}$, 
$L_u (\widetilde{T})$ is an atom of $ \cM_0 (\widetilde{\theta}_{\jmath (u)} T)$ and $d(\rho , \jmath (u))= |u| h $. We now define $f\colon T^\prime \rightarrow S$ by setting $f((u,s))= \jmath (u)$, for all $(u,s) \in T^\prime$. We easily see that 
$$ \forall \, (u,s), (u^\prime, s^\prime) \in T^\prime \, , \quad \left| d^\prime ((u,s),(u^\prime, s))-d(\jmath (u), \jmath (u^\prime)) \right| \leq 2h \; ,$$
which implies that the distortion of $f$ is less than 
$2h$. Since $f(T^\prime )= S$ is a $(2h)$-net of $R_A (T)$,
 $f$ is a pointed $(2h)$-isometry and Lemma \ref{epsiiso} implies that $\dgh( T^\prime , R_A (T)) \leq 8 h$.  This proves that for all $h \in (0, \infty)$, there exists a measurable function $\Phi_h\colon  \bT \rightarrow \bT$, such that 
$$ \forall \, \widetilde{T} \in \bT \, , \quad \dgh \left( \Phi_h (\widetilde{T}) , R_A (\widetilde{T}) \right) \leq 8h \; , $$
which completes the proof of Lemma \ref{measAred}.

%


\end{document}